\documentclass[preprint,3p,11pt]{elsarticle}
\usepackage{url}
\usepackage{amssymb}
\usepackage{amsmath}
\usepackage{stmaryrd}
\usepackage{siunitx}
\usepackage{commath}
\usepackage{subfig}
\usepackage{enumerate}
\usepackage[hidelinks]{hyperref}
\hypersetup{
  colorlinks   = true, 
  urlcolor     = red, 
  linkcolor    = red, 
  citecolor   = red 
}
\usepackage{cleveref}
\usepackage{multirow}
\usepackage{algorithm}
\usepackage{algpseudocode}

\makeatletter
\def\ps@pprintTitle{%
 \let\@oddhead\@empty
 \let\@evenhead\@empty
 \def\@oddfoot{}%
 \let\@evenfoot\@oddfoot}
\makeatother

\newdefinition{remark}{Remark}
\newproof{proof}{Proof}
\begin{document}
\begin{frontmatter}
  \title{A conforming sliding mesh technique for an
    embedded-hybridized discontinuous Galerkin discretization for
    fluid-rigid body interaction}
  \author[TLH]{Tam\'as L. Horv\'ath\corref{cor1}\fnref{label1}}
  \ead{thorvath@oakland.edu}
  \fntext[label1]{\url{https://orcid.org/0000-0001-5294-5362}}
  \address[TLH]{Department of Mathematics and Statistics, Oakland
    University, U.S.A.}
    
  \author[SR]{Sander Rhebergen\fnref{label2}}
  \ead{srheberg@uwaterloo.ca}
  \fntext[label2]{\url{https://orcid.org/0000-0001-6036-0356}}
  \address[SR]{Department of Applied Mathematics, University of
    Waterloo, Canada}  
  \begin{abstract}
    In (J. Comput. Phys., 417, 109577, 2020) we introduced a
    space-time embedded-hybridizable discontinuous Galerkin method for
    the solution of the incompressible Navier--Stokes equations on
    time-dependent domains of which the motion of the domain is
    prescribed. This discretization is exactly mass conserving,
    locally momentum conserving, and energy-stable. In this manuscript
    we extend this discretization to fluid-rigid body interaction
    problems in which the motion of the fluid domain is not known a
    priori. To account for large rotational motion of the rigid body,
    we present a novel conforming space-time sliding mesh
    technique. We demonstrate the performance of the discretization on
    various numerical examples.
  \end{abstract}
  \begin{keyword}
    Navier--Stokes \sep embedded-hybridized \sep discontinuous
    Galerkin \sep space-time \sep fluid-rigid body interaction \sep sliding mesh.
  \end{keyword}
\end{frontmatter}
\section{Introduction}
\label{sec:introduction}

In this paper we consider fluid-rigid body interaction problems in
which the motion of the rigid body is governed by a system of ordinary
differential equations (ODEs) and the fluid problem is governed by the
incompressible Navier--Stokes equations. In such problems aerodynamic
forces exerted by the fluid at the fluid-rigid body interface results
in movement of the rigid body which results in deformation of the
fluid domain. Problems of this kind are challenging to solve as the
fluid domain, on which the Navier--Stokes equations are defined, is
not known a priori. Nevertheless, these problems are of interest in
many engineering applications, for example, in simulations of flow
around a turbine or modeling the interaction between wind and a
bridge.

A popular approach to solve the incompressible Navier--Stokes
equations on deforming domains is by using an Arbitrary
Eulerian--Lagrangian (ALE) approach in which the deforming domain is
mapped to a fixed computational domain. The ALE method has
successfully been applied to fluid-rigid body interaction problems
(e.g., \cite{Robertson:2003,Dettmer:2006,Saksono:2007}), a wide
variety of fluid-structure interaction problems (e.g.,
\cite{Canic:2021,Bukac:2013,Bukac:2019,Neunteufel:2021}), and
free-surface problems (e.g., \cite{Fu:2020,Labeur:2009}).

The space-time method, used in the paper, is an alternative to the ALE
approach. In a space-time approach a time-dependent partial
differential equation (PDE) in $d$ spatial dimensions is first written
as a $(d+1)$-dimensional problem in space-time. One can now naturally
discretize a PDE on a deforming domain by noting that a time-dependent
domain is `fixed' in space-time.

The space-time problem can be discretized in several ways. For
example, in the context of incompressible flows, space-time finite
element methods using function spaces that are continuous in space and
discontinuous in time have been considered in \cite{Hauke:1994,
  Johnson:1994, Ndri:2001, Ndri:2002} and discontinuous function
spaces in both space and time have been considered in
\cite{Rhebergen:2013b,Tavelli:2015, Tavelli:2016, Vegt:2008}. These
methods have successfully been applied to PDEs on deforming domain
problems and fluid-structure interaction problems, see, for example,
\cite{Tezduyar:1992a,Tezduyar:1992b,Tezduyar:2006,Hubner:2004}.

In this paper we consider function spaces that are discontinuous in
space and time. However, unlike space-time DG methods
\cite{Rhebergen:2013b,Tavelli:2015, Tavelli:2016, Vegt:2008}, we
consider here a space-time embedded-hybridizable discontinuous
Galerkin (EHDG) discretization \cite{Horvath:2020}. This
discretization is a space-time extension of the hybridizable
discontinuous Galerkin (HDG) methods introduced previously in
\cite{Rhebergen:2017, Rhebergen:2018, Rhebergen:2020} for
incompressible flows, and a special case of the space-time HDG
discretization (see \cref{ss:discretization_NS}) that was introduced
in \cite{Horvath:2019} and analyzed in \cite{Kirk:2021}. Appealing
properties of this discretization include that it is exactly mass
conserving, locally momentum conserving, and energy-stable. We
furthermore remark that the HDG method was first introduced in
\cite{Cockburn:2009a} to reduce the number of globally coupled
degrees-of-freedom of usual DG methods. This is achieved by static
condensation.

The main focus of this manuscript is the extension of the space-time
EHDG method to fluid-rigid body interaction problems in which the
rigid body movement results in large deformation of the fluid
domain. Large rotational movement of the rigid body may result in a
decrease in the quality of the fluid domain mesh, or even mesh
entanglement, if the mesh construction is not managed
correctly. Therefore, to be able to account for large rotational
movement of the rigid body, we consider a sliding mesh approach
\cite{Ferrer:2012, Cottrell:2009, Anagnostou:1990}. For this, we
partition the computational fluid mesh into three parts: a rotating
mesh, a static mesh, and an annulus mesh. The rigid body is immersed
in the rotating mesh. This rotating mesh rotates with the same angular
velocity as the body and there is no relative motion between the
vertices in this mesh compared to the body. The static mesh, on the
other hand, remains fixed and does not deform with the rigid body
motion. The annulus mesh is a transition zone between the static mesh
and the rotating mesh (see \cref{fig:mixed_mesh}). We use edge
swapping in the annulus mesh to guarantee a good quality conforming
space-time mesh independent of the rotational movement of the rigid
body. This is different from \cite{Ferrer:2012} in which the annulus
mesh results in hanging nodes. We remark that space-time edge swapping
was used also in \cite{Wang:2015}. However, their approach resulted in
new mesh connectivity during the simulation. In the current work, we
will construct the space-time mesh from building blocks that can be
calculated a priori thereby reducing the computational cost that would
normally be associated with building new mesh connectivity.

Finally, we remark that the fluid-rigid body interaction problem can
be solved either by a staggered or a monolithic approach. In the case
of a staggered approach, separate flow (PDE) and structural (ODE)
solvers are applied alternatingly until convergence is reached. Such
an approach is used, for example, in
\cite{Robertson:2003,Calderer:2010,Farhat:2006,Kadapa:2017} and will
be used also in this article. In the monolithic approach the coupled
ODE and PDE systems are solved simultaneously. We refer to
\cite{Dettmer:2006,Fu:2020,Hubner:2004,Rugonyi:2001,Tezduyar:2007} for
more details on the monolithic approach.

The rest of this paper is organized as follows. In
\Cref{s:problem,s:discretization} we introduce, respectively, the
fluid-rigid body interaction problem and its discretization. We
discuss the conforming sliding mesh technique in \Cref{s:sliding_grid}
and we present numerical examples in \Cref{s:examples}. Finally,
conclusions are drawn in \Cref{s:conclusion}.

\section{The fluid-rigid body interaction problem}
\label{s:problem}
We are interested in the solution of fluid-rigid body interaction in
two spatial dimensions and on the time interval $t \in [0, T]$. This
problem is governed by the incompressible Navier--Stokes equations on
a fluid domain coupled to ordinary differential equations that
describe translational and rotational displacement of a single rigid
body. The fluid domain $\Omega(t) \subset \mathbb{R}^2$ is time
dependent due to the motion of the rigid body. It is therefore
convenient to present the Navier--Stokes equations on the space-time
domain,
\begin{equation}
  \mathcal{E} := \cbr{(t, \boldsymbol{x}) \, | \, 0 < t < T,\
    \boldsymbol{x} \in \Omega(t)} \subset \mathbb{R}^3.
\end{equation}
The Navier--Stokes equations for the velocity field
$\boldsymbol{u}:\mathcal{E} \to \mathbb{R}^2$ and the kinematic
pressure field $p:\mathcal{E} \to \mathbb{R}$ are then formulated as:
\begin{equation}
  \label{eq:navierstokes}
  \partial_t\boldsymbol{u} + \nabla \cdot \boldsymbol{\sigma} = \boldsymbol{f},
  \qquad
  \nabla \cdot \boldsymbol{u} = 0,
  \qquad \text{in } \mathcal{E},
\end{equation}
where $\boldsymbol{f}:\mathcal{E} \to \mathbb{R}^2$ is a forcing term
and $\boldsymbol{\sigma}$ is the momentum flux given by
\begin{equation}
  \label{eq:momflux}
  \boldsymbol{\sigma} = \boldsymbol{u} \otimes \boldsymbol{u}
  + p\boldsymbol{I} - 2\nu\boldsymbol{\varepsilon}(\boldsymbol{u}),
\end{equation}
where $\nu > 0$ is the kinematic viscosity, $\boldsymbol{I}$ is the
identity tensor, and
$\boldsymbol{\varepsilon}(\boldsymbol{u}) = (\nabla\boldsymbol{u} +
\nabla\boldsymbol{u}^T)/2$ is the symmetric part of the velocity
gradient.

Denote the space-time outward unit normal to the boundary
$\partial\mathcal{E}$ of the fluid space-time domain $\mathcal{E}$ by
$(n_t, \boldsymbol{n}) \in \mathbb{R}^{3}$, with temporal component
$n_t \in \mathbb{R}$ and spatial component
$\boldsymbol{n} \in \mathbb{R}^2$. This boundary $\partial\mathcal{E}$
is partitioned into five complementary subsets,
$\partial\mathcal{E}^D$ (the Dirichlet boundary),
$\partial\mathcal{E}^N$ (the Neumann boundary on which
$n_t + \boldsymbol{u}\cdot\boldsymbol{n} \ge 0$),
$\partial\mathcal{E}^F$ (the free-slip boundary on which we always
assume $n_t=0$), $\Omega(0)$ (the fluid spatial domain at initial time
$t=0$), and $\Omega(T)$ (the fluid spatial domain at final time
$t=T$). Note that since $n_t=0$ on $\partial\mathcal{E}^F$ we may
write $\partial\mathcal{E}^F = [0,T] \times \partial\Omega^F$, where
$\partial\Omega^F$ is the part of the boundary of $\Omega(t)$ on which
a free slip boundary is imposed. We then denote by $\boldsymbol{\tau}$
the unit tangent vector to $\partial\Omega^F$. We now prescribe the
following boundary conditions:
\begin{subequations}
  \begin{align}
    \label{eq:navierstokes-bc_a}
    \boldsymbol{u}
    &= \boldsymbol{w} && \text{on } \partial\mathcal{E}^{D},
    \\
    \label{eq:navierstokes-bc_c}
    (p \boldsymbol{I} - 2 \nu \boldsymbol{\varepsilon}(\boldsymbol{u})) \boldsymbol{n}
    &= \boldsymbol{g} && \text{on } \partial\mathcal{E}^{N},
    \\
    \label{eq:navierstokes-bc_fs1}
    \boldsymbol{u}\cdot\boldsymbol{n}
    &= 0 && \text{on } \partial\mathcal{E}^{F},
    \\
    \label{eq:navierstokes-bc_fs2}
    \del{(p \boldsymbol{I} - 2 \nu \boldsymbol{\varepsilon}(\boldsymbol{u})) \boldsymbol{n}}\cdot \boldsymbol{\tau}
    &= 0 && \text{on } \partial\mathcal{E}^{F},
    \\
    \label{eq:navierstokes-bc_e}
    \boldsymbol{u}(0, \boldsymbol{x})
    &= \boldsymbol{u}_0(\boldsymbol{x}) && \text{in } \Omega(0),    
  \end{align}
  \label{eq:navierstokes-bc}
\end{subequations}
where $\boldsymbol{g} : \partial\mathcal{E}^N \to \mathbb{R}^2$ is the
prescribed diffusive part of the momentum flux and
$\boldsymbol{u}_0 :\Omega(0) \to \mathbb{R}^2$ is the solenoidal
initial condition. Furthermore, we partition the Dirichlet boundary
into the two complementary sets $\partial\mathcal{E}^B$ and
$\partial\mathcal{E}^W = \partial\mathcal{E}^D \backslash
\partial\mathcal{E}^B$, where $\partial\mathcal{E}^B$ is the boundary
of the rigid body. Then
$\boldsymbol{w} : \partial \mathcal{E}^D \to \mathbb{R}^2$ is
prescribed on $\partial\mathcal{E}^W$ while on $\partial\mathcal{E}^B$
it is the unknown velocity of the interface $\boldsymbol{u}_B$ between
the flow and the body.

\begin{figure}[tbp]
  \begin{center}
    \includegraphics[width=.6\linewidth]{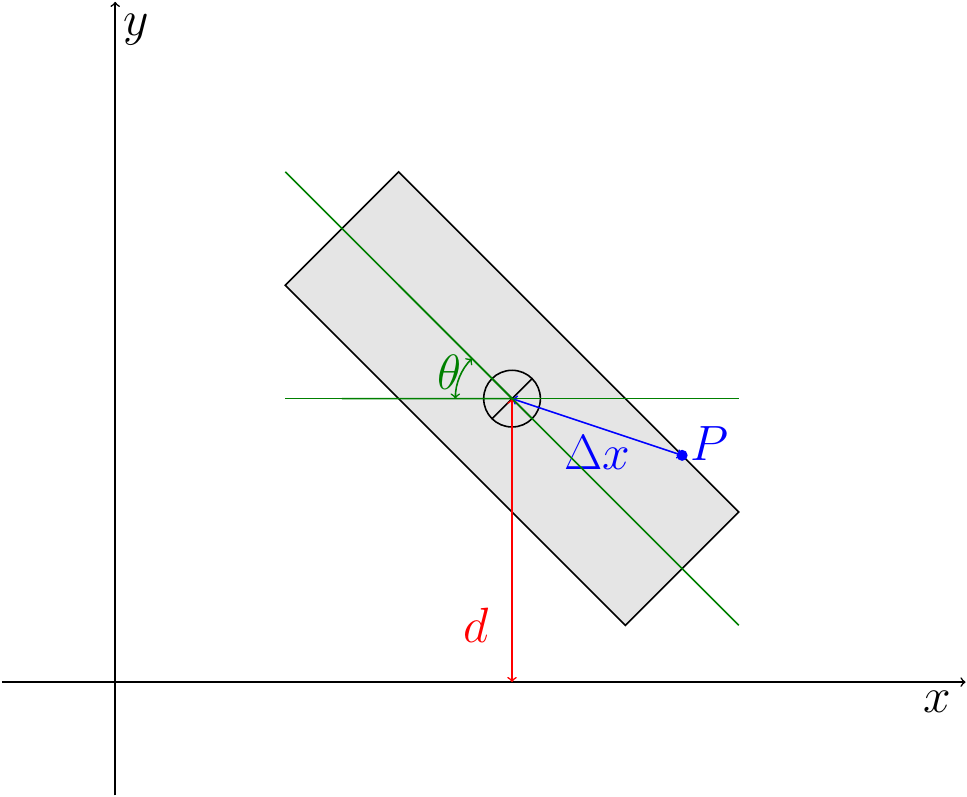} 
    \caption{Clarification of notation for the equations describing
      the motion of a rigid body. The relative position $\Delta x$ of
      a point $P$ on the surface of the rigid body, the rigid body
      rotation $\theta$ are described with respect to the center of
      gravity of the rectangular rigid body, which is denoted here by
      $\otimes$. }
    \label{fig:rb_notation}
  \end{center}
\end{figure}

To determine the interface velocity $\boldsymbol{u}_B$ we require the
force and moment exerted by the rigid body on the fluid. These are
given by, respectively,
\begin{equation}
  \label{eq:forcemoment}
  \boldsymbol{F}(t) = \int_{\partial\mathcal{E}^{B}(t)}\rho
  \del{p\boldsymbol{I}-2\nu\boldsymbol{\varepsilon}(\boldsymbol{u})}\boldsymbol{n} \dif s,
  \quad
  M(t) = \int_{\partial\mathcal{E}^{B}(t)}\Delta\boldsymbol{x} \times
  \sbr{\rho\del{p\boldsymbol{I}-2\nu\boldsymbol{\varepsilon}(\boldsymbol{u})}\boldsymbol{n}} \dif s,
\end{equation}
where $\partial\mathcal{E}^B(t)$ is the fluid-rigid body interface at
time $t$, $\rho$ is a constant density, $\boldsymbol{n}$ is the
outward normal unit vector of the boundary, and $\Delta\boldsymbol{x}$
is the relative position of the surface of the rigid body with respect
to its center of gravity. We consider the case where the rigid body
translation is displacement only in the $y$-direction. We denote this
displacement by $d$. Denoting the rigid body rotation with respect to
the rigid body's center of gravity by $\theta$ (see
\cref{fig:rb_notation} for further clarification of the notation), the
linear motion of the rigid body is described by
\begin{subequations}
  \label{eq:linearMotion}
  \begin{align}
    \label{eq:linearMotion_y}
    m \ddot{d} + c_y \dot{d} + k_y d &= F_y,
    \\
    \label{eq:linearMotion_theta}
    I_{\theta}\ddot{\theta} + c_{\theta}\dot{\theta} + k_{\theta}\theta &= M,    
  \end{align}
\end{subequations}
where $m$ is the mass per unit length of the body, $I_{\theta}$ is the
mass moment of inertia, $c_y$ and $c_{\theta}$ are the damping
constants, and $k_y$ and $k_{\theta}$ are the stiffness
constants. Furthermore, $F_y$ is the $y$-component of the force vector
$\boldsymbol{F}$. The velocity of the interface between the flow and
the body, $\boldsymbol{u}_B$, can then be determined from the solution
to \cref{eq:linearMotion}. 

\section{The discretization}
\label{s:discretization}

\subsection{Incompressible Navier--Stokes}
\label{ss:discretization_NS}

To discretize the fluid problem we first discretize the space-time
domain $\mathcal{E}$. For this we partition the time interval $[0, T]$
into time levels $0 = t^0 < t^1 < \hdots < t^N = T$ and define the
$n$th time interval as $I^n = (t^n, t^{n+1})$ for
$n = 0, \hdots, N-1$. The length of time interval $I^n$ is denoted by
$\Delta t^n = t^{n+1} - t^n$. A space-time slab is then defined as
\begin{equation}
  \mathcal{E}^n = \cbr[0]{(t, \boldsymbol{x}) \in \mathcal{E} \, | \, t \in I^n},
  \qquad n = 0, 1, \hdots, N-1.
\end{equation}
It will be useful to partition the boundary $\partial\mathcal{E}^n$ of
a space-time slab $\mathcal{E}^n$ into the complementary subsets
$\Omega(t^n)$, $\Omega(t^{n+1})$, and
$\mathcal{Q}^n := \cbr[0]{ (t, \boldsymbol{x}) \in \partial\mathcal{E}
  \, | \, t \in I^n}$. The discretization of $\mathcal{E}^n$, for
$n=0, \hdots, N-1$ is denoted by
$\mathcal{E}^n_h = \cup_{k=1}^{M}\mathcal{K}_k$, where $\mathcal{K}_k$
are space-time simplices. We defer to \cref{s:sliding_grid} for a
detailed discussion of creating the mesh. For now we mention that we
consider only conforming triangulations. The boundary
$\partial\mathcal{K}$ of a space-time element
$\mathcal{K} \in \mathcal{E}_h^n$ is formed either by
$\mathcal{Q}_{\mathcal{K}}$, $\mathcal{Q}_{\mathcal{K}} \cup K^n$, or
$\mathcal{Q}_{\mathcal{K}} \cup K^{n+1}$, where
$\mathcal{Q}_{\mathcal{K}} := \cbr[0]{ (t, \boldsymbol{x}) \in
  \partial\mathcal{K} \, | \, t \in I_n}$,
$K^n := \cbr[0]{ (t, \boldsymbol{x}) \in \partial\mathcal{K} \, | \, t
  = t^n}$, and
$K^{n+1} := \cbr[0]{ (t, \boldsymbol{x}) \in \partial\mathcal{K} \, |
  \, t = t^{n+1}}$. The outward unit space-time normal vector on the
boundary $\partial\mathcal{K}$ of $\mathcal{K} \in \mathcal{E}_h^n$ is
denoted by $\del[0]{n_t, \boldsymbol{n}} \in \mathbb{R}^{3}$ where
$n_t$ is the temporal part of the normal and $\boldsymbol{n}$ the
spatial part. Next, we define an interior facet in the time interval
$I^n$ by
$F := \partial\mathcal{Q}_{\mathcal{K}^+} \cap
\partial\mathcal{Q}_{\mathcal{K}^-}$, where $\mathcal{K}^-$ and
$\mathcal{K}^+$ are adjacent elements in $\mathcal{E}_h^n$, while a
boundary facet in the time interval $I^n$ is a facet of
$\mathcal{Q}_{\mathcal{K}}$ that lies on the boundary. The set and
union of all facets in the time interval $I^n$ are denoted by,
respectively, $\mathcal{F}^n$ and $\Gamma^n$.

In each space-time slab $\mathcal{E}_h^n$ the space-time
embedded-hybridized discontinuous Galerkin method results in
approximations $\boldsymbol{u}_h \in \boldsymbol{V}_h^n$ to the
velocity $\boldsymbol{u}$,
$\bar{\boldsymbol{u}}_h^n \in
\bar{\boldsymbol{V}}_h^n(\boldsymbol{w})$ to the trace of the velocity
on element boundaries, $p_h \in Q_h^n$ to the pressure $p$, and
$\bar{p}_h^n \in \bar{Q}_h^n$ to the trace of the pressure on element
boundaries. Here the `velocity' finite element spaces are defined as
\begin{subequations}
  \label{eq:velocity_spc}
  \begin{align}
    \label{eq:velocity_spc_a}
    \boldsymbol{V}_h^n &:= \cbr[1]{ \boldsymbol{v}_h \in \sbr[0]{L^2(\mathcal{E}_h^n)}^2 \, | \,
      \boldsymbol{v}_h \in \sbr[0]{P_k(\mathcal{K})}^2\ \forall \mathcal{K} \in \mathcal{E}^n_h},
    \\
    \label{eq:velocity_spc_b}
    \bar{\boldsymbol{V}}^n_h(\boldsymbol{w})
    &:= \big\{ \bar{\boldsymbol{v}}_h \in \sbr[0]{L^2(\mathcal{F}^n)}^2 \, | \,
      \bar{\boldsymbol{v}}_h \in \sbr[0]{P_k(F)}^2
      \ \forall F \in \mathcal{F}^n,\ \bar{\boldsymbol{v}}_h = \boldsymbol{w} \
      \text{on } \partial \mathcal{E}^D \cap \partial \mathcal{E}^n,
    \\ \nonumber
    & \hspace{21em} \bar{\boldsymbol{v}}_h \cdot \boldsymbol{n} = 0 \ \text{on } \partial\mathcal{E}^F \cap \partial\mathcal{E}^n \big \} \cap C(\Gamma^n),
  \end{align}
\end{subequations}
where $P_l(D)$ denotes the space of polynomials of degree $l \ge 0$ on
a domain $D$. The `pressure' finite element spaces are defined as
\begin{subequations}
  \label{eq:pressure_spc}
  \begin{align}
    \label{eq:pressure_spc_a}
    Q_h^n &:= \cbr[1]{ q_h \in L^2(\mathcal{E}_h^n) \, | \, q_h \in P_{k-1}(\mathcal{K}) \
      \forall \mathcal{K} \in \mathcal{E}^n_h},
    \\
    \label{eq:pressure_spc_b}
    \bar{Q}^n_h
    &:= \cbr[1]{ \bar{q}_h \in L^2(\mathcal{F}) \, | \,
      \bar{q}_h \in P_k(F) \ \forall F\in\mathcal{F}^n}.
  \end{align}
\end{subequations}
Note that the space-time EHDG method approximates the trace of the
velocity using a continuous facet space but that the trace of the
pressure is approximated in a discontinuous facet space. We remark,
however, that with minimal changes to the space-time EHDG
discretization it is possible to obtain also a space-time HDG (in
which both facet spaces are discontinuous) and a space-time EDG (in
which both facet spaces are continuous) discretization of the
fluid-rigid body interaction problem. See \cite{Horvath:2020} for a
detailed comparison of the three methods.

The space-time EHDG discretization of the momentum and continuity
equations \cref{eq:navierstokes} in each space-time slab
$\mathcal{E}_h^n$, $n = 0, \hdots, N-1$, is now given by
\cite{Horvath:2019, Horvath:2020}: find
$(\boldsymbol{u}_h, \bar{\boldsymbol{u}}_h) \in \boldsymbol{V}_h^n
\times \bar{\boldsymbol{V}}_h^n(\boldsymbol{w})$ and
$(p_h, \bar{p}_h) \in Q_h^n \times \bar{Q}_h^n$ such that
\begin{subequations}
  \label{eq:ST-EHDG}
  \begin{equation}
    \label{eq:momentumdisc}
    \begin{split}
      &-\sum_{\mathcal{K} \in \mathcal{E}_h^n} \int_{\mathcal{K}}
      \big( \boldsymbol{u}_h \partial_t \boldsymbol{v}_h
      + \boldsymbol{\sigma}_h:\nabla\boldsymbol{v}_h \big) \dif\boldsymbol{x}\dif t
      + \sum_{\mathcal{K}\in\mathcal{E}^n_h} \int_{\mathcal{Q}_{\mathcal{K}}}
      \widehat{\boldsymbol{\sigma}}_h \cdot (\boldsymbol{v}_h - \bar{\boldsymbol{v}}_h) \dif s
      + \int_{\Omega(t_{n+1})}\boldsymbol{u}_h \cdot \boldsymbol{v}_h \dif \boldsymbol{x}
      \\
      &
      - \sum_{\mathcal{K}\in\mathcal{E}^n_h} \int_{\mathcal{Q}_{\mathcal{K}}} 2\nu\boldsymbol{\varepsilon}(\boldsymbol{v}_h) :
      (\boldsymbol{u}_h - \bar{\boldsymbol{u}}_h)\otimes\boldsymbol{n} \dif s
      + \int_{\partial\mathcal{E}^N \cap I_n} \del{n_t + \bar{\boldsymbol{u}}_h\cdot \boldsymbol{n}}
      \bar{\boldsymbol{u}}_h \cdot \bar{\boldsymbol{v}}_h \dif s  
      \\
      =&
      \sum_{\mathcal{K}\in\mathcal{E}^n_h} \int_{\mathcal{K}} \boldsymbol{f} \cdot \boldsymbol{v}_h \dif \boldsymbol{x} \dif t
      - \int_{\partial\mathcal{E}^N \cap I^n} \boldsymbol{g} \cdot \bar{\boldsymbol{v}}_h \dif s 
      + \int_{\Omega(t_n)} \boldsymbol{u}_h^- \cdot \boldsymbol{v}_h \dif \boldsymbol{x},      
    \end{split}
  \end{equation}
  for all
  $(\boldsymbol{v}_h, \bar{\boldsymbol{v}}_h) \in \boldsymbol{V}_h^n
  \times \bar{\boldsymbol{V}}_h^n(\boldsymbol{0})$ and
  \begin{equation}
    \label{eq:contdisc}
    - \sum_{\mathcal{K}\in\mathcal{E}^n_h} \int_{\mathcal{K}} q_h \nabla \cdot \boldsymbol{u}_h \dif \boldsymbol{x} \dif t 
    + \sum_{\mathcal{K}\in\mathcal{E}^n_h} \int_{\mathcal{Q}_{\mathcal{K}}} (\boldsymbol{u}_h - \bar{\boldsymbol{u}}_h)
    \cdot \boldsymbol{n} \bar{q}_h \dif s = 0,
  \end{equation}
\end{subequations}
for all $(q_h, \bar{q}_h) \in Q_h^n \times \bar{Q}_h^n$. When $n=0$,
$\boldsymbol{u}_h^-$ is the projection of the initial condition
$\boldsymbol{u}_0$ into $\boldsymbol{V}_h^0 \cap H(\text{div})$ such
that $\nabla \cdot \boldsymbol{u}_h^- = 0$ element-wise. When $n>0$,
$\boldsymbol{u}_h^- = \lim_{\varepsilon\to 0}
\boldsymbol{u}_h(t^n-\varepsilon)$. Furthermore,
$\widehat{\boldsymbol{\sigma}}_h$ is an approximation to the normal
component of the momentum flux on element boundaries, which is given
by
\begin{equation}
  \label{eq:numflux}
  \widehat{\boldsymbol{\sigma}}_h :=
  (n_t + \boldsymbol{u}_h \cdot \boldsymbol{n})(\boldsymbol{u}_h + \lambda(\bar{\boldsymbol{u}}_h - \boldsymbol{u}_h))
  + (\bar{p}_h\boldsymbol{I} - 2\nu\boldsymbol{\varepsilon}(\boldsymbol{u}_h))\boldsymbol{n}
  + \frac{2\nu\alpha}{h_{\mathcal{K}}}(\boldsymbol{u}_h - \bar{\boldsymbol{u}}_h),
\end{equation}
where $\lambda=1$ when
$(n_t + \boldsymbol{u}_h \cdot \boldsymbol{n}) < 0$ and $\lambda = 0$
when $(n_t + \boldsymbol{u}_h \cdot \boldsymbol{n}) \ge 0$,
$\alpha > 0$ is a penalty parameter that needs to be sufficiently
large to ensure stability, and $h_{\mathcal{K}}$ is the characteristic
length of the element $\mathcal{K}$. We remark that the advective part
of this numerical flux is a space-time upwind flux while the diffusive
part is an interior penalty flux.

In \cite{Horvath:2019, Horvath:2020} we proved that the space-time
EHDG discretization \cref{eq:ST-EHDG} results in an approximate
velocity field $\boldsymbol{u}_h$ that is exacty divergence-free and
$H(\text{div})$-conforming on space-time elements, i.e.,
$\nabla \cdot \boldsymbol{u}_h = 0$ for
$(t, \boldsymbol{x}) \in \mathcal{K}$ for all
$\mathcal{K} \in \mathcal{E}_h^n$,
$\boldsymbol{u}_h \cdot \boldsymbol{n}$ is single-valued on interior
faces and
$\boldsymbol{u}_h \cdot \boldsymbol{n} = \bar{\boldsymbol{u}}_h \cdot
\boldsymbol{n}$ on boundary faces. A consequence of these results is
that the discretization is exactly mass conserving and that the
discretization is both locally conservative and energy stable, even on
deforming domains.

The space-time EHDG discretization \cref{eq:ST-EHDG} in space-time
slab $\mathcal{E}_h^n$, $n=0,\hdots,N-1$, results in a nonlinear
system of algebraic equations. To solve this nonlinear problem we use
Picard iterations: given
$(\boldsymbol{u}_h^k, \bar{\boldsymbol{u}}_h^k) \in \boldsymbol{V}_h^n
\times \bar{\boldsymbol{V}}_h^n(\boldsymbol{w})$ we seek
$(\boldsymbol{u}_h^{k+1}, \bar{\boldsymbol{u}}_h^{k+1}) \in
\boldsymbol{V}_h^n \times \bar{\boldsymbol{V}}_h^n(\boldsymbol{w})$
and $(p_h^{k+1}, \bar{p}_h^{k+1}) \in Q_h^n \times \bar{Q}_h^n$ such
that
\begin{subequations}
  \label{eq:Picard}
  \begin{equation}
    \label{eq:momentum_Picard}
    \begin{split}
      &-\sum_{\mathcal{K} \in \mathcal{E}_h^n} \int_{\mathcal{K}}
      \big( \boldsymbol{u}_h^{k+1} \partial_t \boldsymbol{v}_h
      + \boldsymbol{\sigma}_h^k:\nabla\boldsymbol{v}_h \big) \dif\boldsymbol{x}\dif t
      + \sum_{\mathcal{K}\in\mathcal{E}^n_h} \int_{\mathcal{Q}_{\mathcal{K}}}
      \widehat{\boldsymbol{\sigma}}_h^k \cdot (\boldsymbol{v}_h - \bar{\boldsymbol{v}}_h) \dif s
      + \int_{\Omega(t_{n+1})}\boldsymbol{u}_h^{k+1} \cdot \boldsymbol{v}_h \dif \boldsymbol{x}
      \\
      &
      - \sum_{\mathcal{K}\in\mathcal{E}^n_h} \int_{\mathcal{Q}_{\mathcal{K}}} 2\nu\boldsymbol{\varepsilon}(\boldsymbol{v}_h) :
      (\boldsymbol{u}_h^{k+1} - \bar{\boldsymbol{u}}_h^{k+1})\otimes\boldsymbol{n} \dif s
      + \int_{\partial\mathcal{E}^N \cap I_n} \del{n_t + \bar{\boldsymbol{u}}_h^k\cdot \boldsymbol{n}}
      \bar{\boldsymbol{u}}_h^{k+1} \cdot \bar{\boldsymbol{v}}_h \dif s  
      \\
      =&
      \sum_{\mathcal{K}\in\mathcal{E}^n_h} \int_{\mathcal{K}} \boldsymbol{f} \cdot \boldsymbol{v}_h \dif \boldsymbol{x} \dif t
      - \int_{\partial\mathcal{E}^N \cap I^n} \boldsymbol{g} \cdot \bar{\boldsymbol{v}}_h \dif s 
      + \int_{\Omega(t_n)} \boldsymbol{u}_h^- \cdot \boldsymbol{v}_h \dif \boldsymbol{x},      
    \end{split}
  \end{equation}
  for all
  $(\boldsymbol{v}_h, \bar{\boldsymbol{v}}_h) \in \boldsymbol{V}_h^n
  \times \bar{\boldsymbol{V}}_h^n(\boldsymbol{0})$ and
  \begin{equation}
    \label{eq:cont_Picard}
    - \sum_{\mathcal{K}\in\mathcal{E}^n_h} \int_{\mathcal{K}} q_h \nabla \cdot \boldsymbol{u}_h^{k+1} \dif \boldsymbol{x} \dif t 
    + \sum_{\mathcal{K}\in\mathcal{E}^n_h} \int_{\mathcal{Q}_{\mathcal{K}}} (\boldsymbol{u}_h^{k+1} - \bar{\boldsymbol{u}}_h^{k+1})
    \cdot \boldsymbol{n} \bar{q}_h \dif s = 0,
  \end{equation}
\end{subequations}
for all $(q_h, \bar{q}_h) \in Q_h^n \times \bar{Q}_h^n$ for
$k=0,1,2,\hdots$. Here $\boldsymbol{\sigma}_h^k$ and
$\widehat{\boldsymbol{\sigma}}_h^k$ are the following modifications
of, respectively, \cref{eq:momflux,eq:numflux}:
\begin{align*}
  \boldsymbol{\sigma}_h^k
  & = \boldsymbol{u}_h^k \otimes \boldsymbol{u}_h^{k+1}
    + p_h^{k+1}\boldsymbol{I} - 2\nu\boldsymbol{\varepsilon}(\boldsymbol{u}_h^{k+1}),
  \\
  \widehat{\boldsymbol{\sigma}}_h^k
  &=   (n_t + \boldsymbol{u}_h^k \cdot \boldsymbol{n})(\boldsymbol{u}_h^{k+1} + \lambda(\bar{\boldsymbol{u}}_h^{k+1} - \boldsymbol{u}_h^{k+1}))
    + (\bar{p}_h^{k+1}\boldsymbol{I} - 2\nu\boldsymbol{\varepsilon}(\boldsymbol{u}_h^{k+1})\boldsymbol{n}
    + \frac{2\nu\alpha}{h_{\mathcal{K}}}(\boldsymbol{u}_h^{k+1} - \bar{\boldsymbol{u}}_h^{k+1}).
\end{align*}
To initialize the Picard iteration we set $\boldsymbol{u}_h^0 = 0$ on
all elements in $\mathcal{E}_h^n$,
$\bar{\boldsymbol{u}}_h^0 = \boldsymbol{w}$ on Dirichlet boundary
facets, and $\bar{\boldsymbol{u}}_h^0 = 0$ on all remaining facets in
$\mathcal{F}^n$. As stopping criterion we use
\begin{equation}
  \label{eq:stopping_criterion}
  \max \cbr{
    \frac{\norm[0]{\boldsymbol{u}_h^{k+1} - \boldsymbol{u}_h^{k}}_{\infty}}{\norm[0]{\boldsymbol{u}_h^{k+1} - \boldsymbol{u}_h^0}_{\infty}},
    \frac{\norm[0]{p_h^{k+1} - p_h^{k}}_{\infty}}{\norm[0]{p_h^{k+1} - p_h^0}_{\infty}} }
  < \delta_{\text{NS}},
\end{equation}
where $\norm[0]{\cdot}_\infty$ is the discrete maximum-norm and
$\delta_{\text{NS}}$ is the desired tolerance. Once the convergence
criterion is met we set
$(\boldsymbol{u}_h, \boldsymbol{p}_h) = (\boldsymbol{u}_h^{k+1},
\boldsymbol{p}_h^{k+1})$.

\subsection{Rigid body time integration}
\label{ss:discretization_ODE}

In this section we discretize the equations that describe the motion
of the rigid body \cref{eq:linearMotion}. Since the discretization for
\cref{eq:linearMotion_y} and \cref{eq:linearMotion_theta} are similar,
we present here only the details for \cref{eq:linearMotion_y}.

We write the second-order ordinary differential equation
\cref{eq:linearMotion_y} as the following first-order system:
\begin{equation}
  \label{eq:firstordersystem}
  \dot{d} = b, \qquad \dot{b} = (-k_yd - c_yb + F_y)/m.
\end{equation}
We then discretize \cref{eq:firstordersystem} by a predictor-corrector
method. In particular, starting from the initial condition
$(d_0, b_0)$, we generate the values $(d_{n+1}, b_{n+1})$ for
$n = 0, 1, \ldots, N-1$ using the following algorithm. In the
predictor step we first calculate an initial guess
$(d_{n+1}^0, b_{n+1}^0)$ using Euler's method:
\begin{subequations}
  \label{eq:predictor}
  \begin{align}
    d_{n+1}^0 &= \Delta t^n b_n + d_n,
    \\
    b_{n+1}^0 &= \Delta t^n \del{ -k_y d_n - c_y b_n + F_y(t_n) } / m + b_n.
  \end{align}
\end{subequations}
We then update the initial guess in the corrector step to
$(d_{n+1}^l, b_{n+1}^l)$ for $l = 1, 2, \ldots$ until a stopping
criterion is met. These updates are given by a BDF2 discretization:
\begin{subequations}
  \label{eq:corrector}
  \begin{align}
    d_{n+1}^l &= \frac{2}{3} \Delta t^n b_{n+1}^{l-1} + \frac{4}{3} d_n - \frac{1}{3} d_{n-1},
    \\
    b_{n+1}^l &= \frac{2}{3} \Delta t^n \del{-k_y d_{n+1}^{l-1} - c_y b_{n+1}^{l-1} + F_y(t_{n+1}) } / m
            + \frac{4}{3} b_n - \frac{1}{3} b_{n-1}.
  \end{align}
\end{subequations}
We use the following stopping criterion:
\begin{equation}
  \label{eq:stopping_crit_rb}
  \sqrt{(d_{n+1}^l - d_{n+1}^{l-1})^2 + (b_{n+1}^l - b_{n+1}^{l-1})^2} < \delta_{\text{rb}},
\end{equation}
where $\delta_{\text{rb}}$ is the desired tolerance. When the stopping
criterion is met, we define $(d_{n+1}, b_{n+1})$ as the last corrector
step.

The interface velocity between the fluid and the body,
$\boldsymbol{u}_B$, is updated at every iteration using
$(d_{n+1}^l, b_{n+1}^l)$ and $(\theta_{n+1}^l, \kappa_{n+1}^l)$, where
$(\theta_{n+1}^l, \kappa_{n+1}^l)$ is an approximation to
$(\theta(t^{n+1}), \dot{\theta}(t^{n+1}))$.

\subsection{Staggered coupling for fluid-rigid body interaction}
\label{ss:discretization_coupling}

Given the discretization of the fluid problem and the equations
describing the motion of the rigid body in, respectively,
\cref{ss:discretization_NS,ss:discretization_ODE}, we now briefly
discuss the coupling of these discretizations when advancing from
$t^n$ to $t^{n+1}$.

Given the velocity and pressure solution at time $t=t^n$, we first
estimate the rigid body movement by \cref{eq:predictor} after which we
update the flow domain and the mesh $\mathcal{E}_h^n$. Given the new
mesh, we solve the flow equations by Picard iteration \cref{eq:Picard}
until \cref{eq:stopping_criterion} is satisfied. The updated velocity
and pressure solutions are then used to calculate the forces
\cref{eq:forcemoment} after which we apply the corrector step
\cref{eq:corrector}. If we satisfy \cref{eq:stopping_crit_rb}, we move
to the next time slab; otherwise, we again update the flow domain and
solve the fluid problem. See \cref{alg:frbi}.

\begin{algorithm}[!h!t]
  \caption{Staggered coupling for fluid-rigid body solver in space-time slab $\mathcal{E}_h^n$. \label{alg:frbi}}
  \begin{algorithmic}[1]
    \State Predictor step \cref{eq:predictor} to obtain an initial guess for the rigid body position
    \While {Rigid body stopping criterion \cref{eq:stopping_crit_rb} not satisfied}
    \State Update the flow domain and mesh $\mathcal{E}_h^n$
    \While {Picard stopping criterion \cref{eq:stopping_criterion} not satisfied}
    \State Solve \cref{eq:Picard} to obtain the flow solution
           $(\boldsymbol{u}_h^{k+1}, \bar{\boldsymbol{u}}_h^{k+1}, p_h^{k+1},\bar{p}_h^{k+1})$
    \EndWhile
    \State Set $(\boldsymbol{u}_h, \bar{\boldsymbol{u}}_h, p_h,\bar{p}_h)
           =(\boldsymbol{u}_h^{k+1}, \bar{\boldsymbol{u}}_h^{k+1}, p_h^{k+1},\bar{p}_h^{k+1})$
    \State Corrector step \cref{eq:corrector} to update the rigid body position
    \EndWhile
  \end{algorithmic}
\end{algorithm}

\section{A conforming sliding mesh technique}
\label{s:sliding_grid}

Due to the fluid-rigid body interaction, an update of the flow
domain/mesh $\mathcal{E}_h^n$ is required in each space-time slab (see
line 3 of \cref{alg:frbi}). However, for some problems the rotational
motion and displacement of the rigid body can be large which may
result in mesh entanglement if the flow domain is not re-meshed. To
account for large deformation of the domain we therefore consider here
a conforming sliding mesh technique.

In this section we will first discuss the extension of a triangular
spatial mesh to a tetrahedral space-time mesh after which we explain
the necessary edge swapping to obtain a conforming sliding mesh. We
will assume that the number of spatial elements $K^n$ is constant for
all $n=0,\hdots,N$ and we will denote the (global) identifiers of the
three vertices of an element $K^n$ by $p_i^n$ with
$i \in \cbr{1,\hdots, N_v}$ with $N_v$ the total number of vertices in the
spatial mesh.

\subsection{Tetrahedral space-time mesh generation}
\label{ss:prism2tet}

Let $\Omega_h(t^n) := \cup_{j=1}^JK_j^n$ be the discretization of the
spatial domain $\Omega(t)$ at time level $t=t^n$, with $J$ the total
number of elements in the spatial mesh, and let
$\boldsymbol{x}_i^{j,n}$, $i=1,2,3$, be the spatial coordinates of the
vertices of a spatial element $K_j^n$. Furthermore, let
$\Phi^n:\Omega_h(t^n) \to \Omega_h(t)\ :\ \boldsymbol{x} \mapsto
\Phi^n(\boldsymbol{x})$ be the mapping that represents the evolution
of the spatial domain due to the fluid-rigid body interaction from
$t^n$ to $t^{n+1}$. Following \cite{Vegt:2002}, the spatial elements
$K_j^{n+1}$ at time level $t=t^{n+1}$ are obtained from $K_j^n$ by
moving the vertices $\boldsymbol{x}_i^{j,n}$ to their new position at
time level $t=t^{n+1}$ using the mapping $\Phi^n$. We remark that if
the vertex $\boldsymbol{x}_i^{j,n}$ has global identifier $p_j^n$,
then the vertex $\boldsymbol{x}_i^{j,n+1}=\Phi^n(\boldsymbol{x}_i^n)$
has global identifier $p_j^{n+1} = p_j^n + N_v$. A space-time element
can then be obtained by connecting $K_j^n$ and $K_j^{n+1}$ via linear
interpolation in time. For ease of discussion we call this space-time
element a space-time `prism' but remark that, depending on the mapping
$\Phi^n$, this prism may not have flat quadrilateral sides. In our
discretization, however, we never create a space-time mesh consisting
of prisms. Instead, the space-time prisms are used as an intermediate
step to create a space-time mesh consisting only of tetrahedra. This
is because our space-time EHDG discretization \cref{eq:ST-EHDG} is
exactly mass conserving, locally momentum conservative, and
energy-stable on simplicial meshes (see
\cite{Horvath:2019,Horvath:2020}). A further consequence of cutting
the prisms into tetrahedra is that all boundary sides of the
space-time tetrahedra are flat surfaces.

\begin{figure}[tbp]
  \begin{center}
    \includegraphics[width=.9\linewidth]{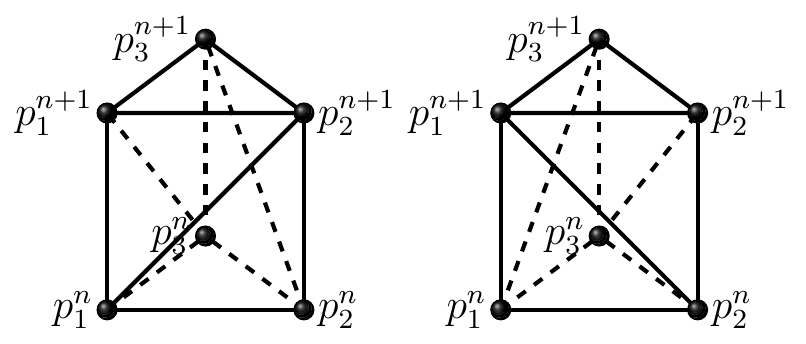} 
    \caption{Invalid cuts of a prism.}
    \label{fig:bad_prison_cut}
  \end{center}
\end{figure}

There are six different ways to cut a prism into three tetrahedra and
these can be identified by the way the quadrilateral sides of a prism
are cut (see \cite{Wang:2015}). To identify these different cuts, let
us call a set of diagonal cuts of the quadrilateral sides
\textit{valid} if they generate three tetrahedra. The three
quadrilateral sides of a prism can be cut along either of the two
diagonals resulting in $2^3$ possible ways to cut a prism. Only six of
these cuts, however, are valid and these are the cuts where two of the
quadrilateral diagonals share a vertex. The two invalid cuts are shown
in \cref{fig:bad_prison_cut}.

\begin{figure}[tbp]
  \begin{center}
    \includegraphics[width=\linewidth]{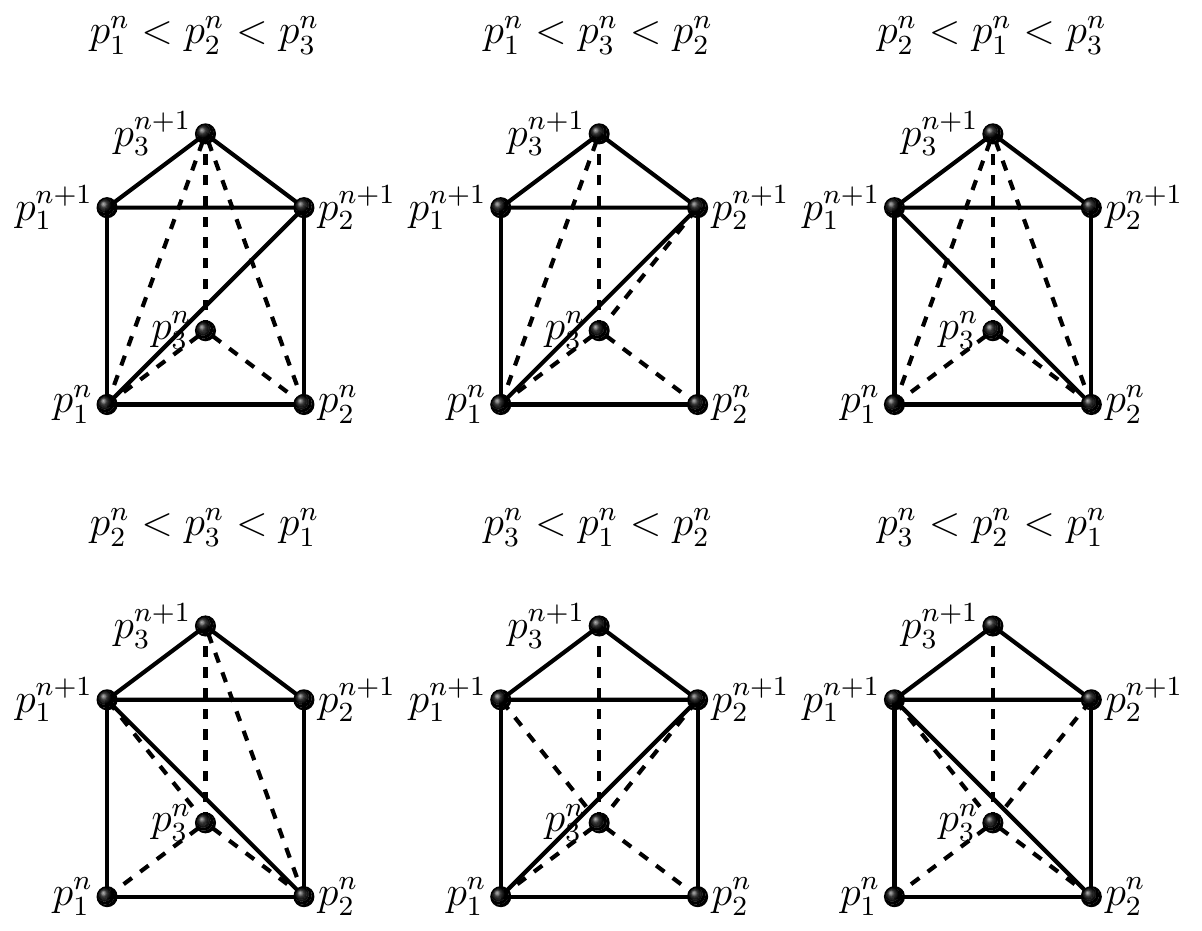} 
    \caption{Valid cuts of a prism.}
    \label{fig:proper_prison_cut}
  \end{center}
\end{figure}

A valid cut is guaranteed by using the quadrilateral diagonal starting
at the vertex with the smallest identifier on every quadrilateral
side; this vertex is included in two quadrilateral cuts (see
\cref{fig:proper_prison_cut}). This approach generates a conforming
tetrahedral mesh since the quadrilateral side shared by two
neighboring prisms is divided into two triangles for both prisms in
the same way.

\subsection{Edge swapping}

\begin{figure}[tbp]
  \begin{center}
    \includegraphics[width=\linewidth]{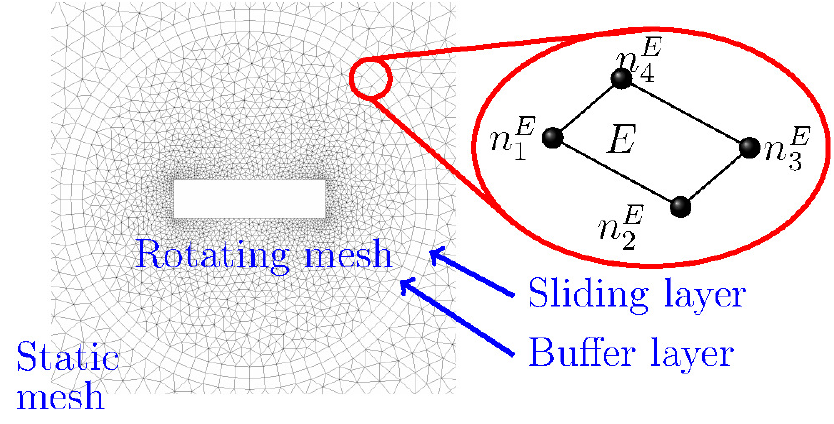} 
    \caption{The initial mixed spatial mesh and the numbering of the
      quadrilateral nodes.}
    \label{fig:mixed_mesh}
  \end{center}
\end{figure}

In this paper, we consider sliding meshes to account for large
deformation of the flow domain. Following for example
\cite{Ferrer:2012}, we consider a mesh where an inner mesh rotates
with respect to an outer static mesh (see
\cref{fig:mixed_mesh}). Different from \cite{Ferrer:2012}, however, is
that we consider edge swapping in the annulus separating the static
and rotating mesh to avoid hanging nodes. This preserves exact mass
conservation of our EHDG discretization and furthermore has the
advantage that we do not require internal curved edges. In this
section we explain the edge swapping to obtain a conforming space-time
sliding mesh in the time interval $[t^n , t^{n+1}]$.

The spatial mesh at time level $t=t^n$ is fixed since it is the upper
boundary of the space-time mesh $\mathcal{E}_h^{n-1}$ from the
previous time interval $[t^{n-1}, t^n]$. Consider therefore the
spatial mesh at time $t=t^{n+1}$. Initially, we consider this to be a
mixed mesh consisting of quadrilaterals in the annulus that separates
the inner and outer mesh, and triangles elsewhere. The annulus itself
is separated into two layers of elements resulting in an inner
annulus, which we refer to as the buffer layer, and an outer annulus,
which we refer to as the sliding layer (see
\cref{fig:mixed_mesh}). When the inner mesh rotates, all vertices in
the annulus move except the vertices shared also by the outer static
mesh. In the initial mesh we will require an even number of
quadrilaterals in the buffer layer and that the number of
quadrilaterals in the buffer and sliding layers are equal.

\begin{figure}[tbp]
  \begin{center}
    \includegraphics[width=\linewidth]{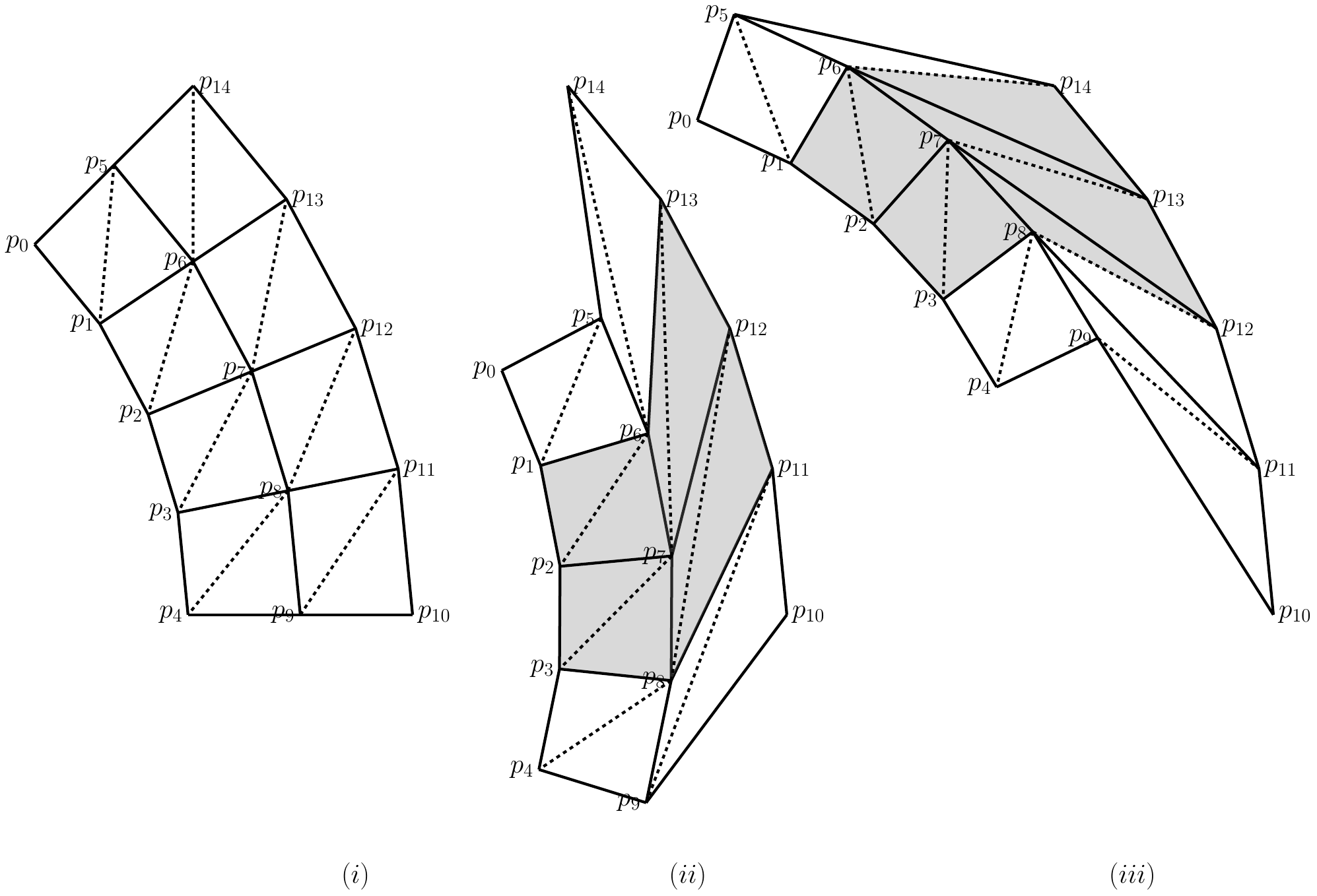}
    \caption{Here we show part of the triangular mesh in the sliding
      layer. In $(i)$ we show part of this mesh before rotation. In
      $(ii)$ and $(iii)$ we show part of this mesh once the inner mesh
      has started to rotate. In $(ii)$ edge swapping is necessary in
      the quadrilaterals of the initial quadrilateral mesh. In $(iii)$
      edge swapping is necessary in the quadrilaterals of the
      secondary quadrilateral mesh.}
    \label{fig:mixed_ghost}
  \end{center}
\end{figure}

Consider now a quadrilateral $E$ in the spatial mesh at time
$t=t^{n+1}$ with vertices $n_i^E$, $i=1, \hdots, 4$. We introduce the
following numbering convention: $n_1^E$ and $n_2^E$ are the vertices
closest to the center of rotation and ordered clockwise on the annulus
while the remaining two vertices $n_3^E$ and $n_4^E$ are ordered
anticlockwise with respect to the center of the quadrilateral (see
\cref{fig:mixed_mesh}). Since we are interested in generating a
tetrahedral space-time mesh on the whole domain (see
\cref{ss:prism2tet}), we next cut the quadrilaterals along their
$n_2^E - n_4^E$ diagonals to obtain a triangular mesh on the whole
spatial domain. These diagonal cuts now form a secondary quadrilateral
mesh in the sliding layer. For example, the quadrilaterals
$\cbr[0]{p_6^{n+1}, p_7^{n+1}, p_{13}^{n+1}, p_{14}^{n+1}}$ and
$\cbr[0]{p_7^{n+1}, p_8^{n+1}, p_{12}^{n+1}, p_{13}^{n+1}}$ in
\cref{fig:mixed_ghost}$(i)$ are elements in a secondary quadrilateral
mesh. Note that cutting quadrilaterals in the secondary mesh along the
$n_1^E - n_3^E$ diagonal results in the same simplicial mesh as
cutting quadrilaterals in the primary mesh along the $n_2^E - n_4^E$
diagonal.

Rotation of the inner mesh may result in triangles in the sliding
layer of which their inscribed circles have a small radius thereby
resulting in a low quality spatial mesh at time $t=t^{n+1}$. This can
be avoided, however, by edge swapping. The edge swapping can occur on
the initial quadrilateral mesh or on the secondary quadrilateral
mesh. Consider, for example, \cref{fig:mixed_ghost}$(ii)$. Here we
would swap the diagonals of the initial quadrilateral mesh, i.e., we
would replace $p_6^{n+1} - p_{14}^{n+1}$ by
$p_5^{n+1} - p_{13}^{n+1}$, and similarly for the other
quadrilaterals. In the case of \cref{fig:mixed_ghost}$(iii)$, however,
we would swap the diagonal of the secondary quadrilateral mesh, i.e.,
we would replace $p_6^{n+1} - p_{13}^{n+1}$ by
$p_7^{n+1} - p_{14}^{n+1}$, and similarly for the remaining
quadrilaterals of the secondary mesh. Note that the same diagonal cut
is used in all elements in the sliding layer. To determine whether
edge swapping is necessary in the sliding layer, we therefore need to
compute the diagonals only of one quadrilateral element in the primary
mesh and of one quadrilateral in the secondary mesh. We always choose
the mesh built from quadrilaterals with the shortest diagonal in the
sliding layer; if the shortest diagonal appears in the current mesh,
no edge swapping is necessary, otherwise it is. There is no edge
swapping in the buffer layer.

Given the spatial mesh at times $t=t^n$ and $t=t^{n+1}$, we next
construct the space-time tetrahedral mesh. At this point we remark
that $p_j^n$ is always connected to $p_j^{n+1}$ by linear
interpolation in time (see \cref{ss:prism2tet}) and that the
connectivity in the inner rotating mesh and in the outer static mesh
do not change throughout the simulation. If no edge swapping takes
place at $t=t^{n+1}$ we can directly use the algorithm from
\cref{ss:prism2tet}. Otherwise, we note that there are only four
different connectivity sets for the space-time mesh of the annulus
when edge swapping does take place and these can be calculated before
the start of a simulation. During the simulation, in each space-time
slab, we then only need to choose one of the four connectivity
sets. We discuss this next.

\begin{figure}[tbp]
  \begin{center}
    \includegraphics[width=.4\linewidth]{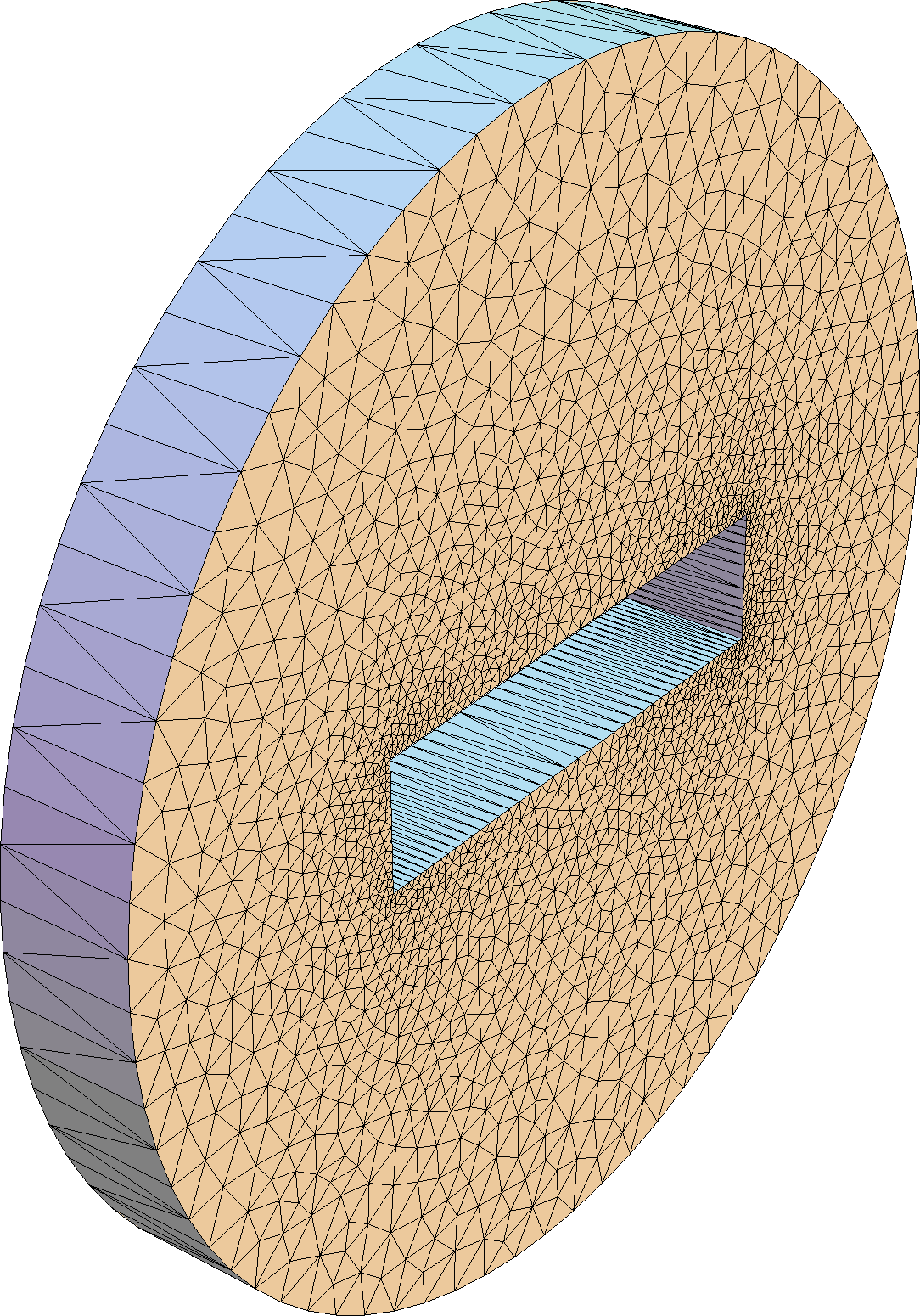} 
    \includegraphics[width=.4\linewidth]{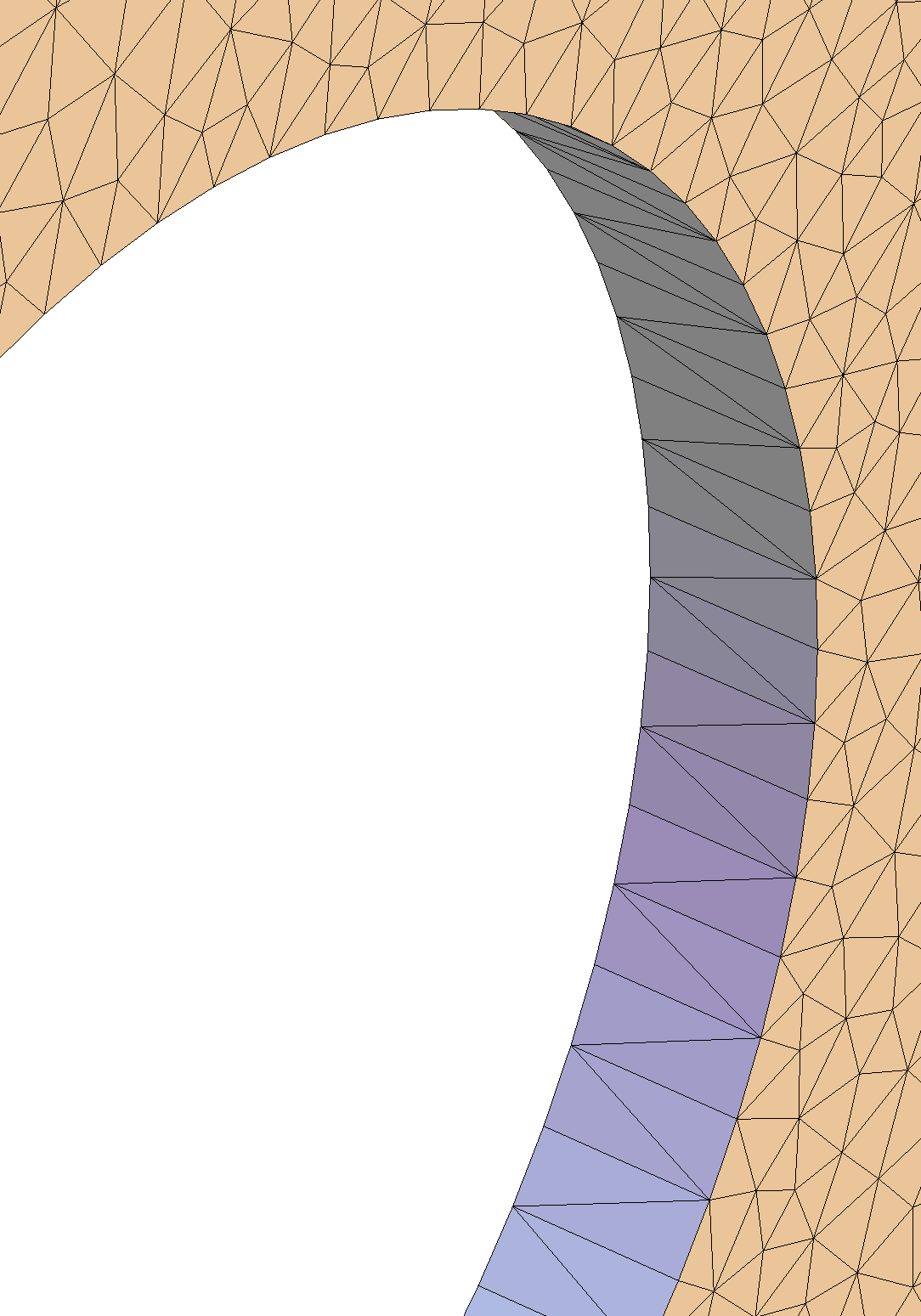} 
    \caption{On the outer boundary of the rotating space-time inner
      mesh (left) and on the inner boundary the outer space-time
      static mesh (right) we observe a `chainsaw' pattern.}
    \label{fig:chainsaw}
  \end{center}
\end{figure}

To construct the space-time mesh, consider again the annulus in
\cref{fig:mixed_mesh}. The boundary of this annulus consists of the
inner boundary circle and the outer boundary circle. For ease of
discussion, we now introduce the following numbering of the vertices
on the inner boundary circle: no three adjacent vertices are in an
increasing or decreasing order, i.e., for any three adjacent vertices,
the identifier of the middle vertex is either smaller or larger than
its two neighboring vertices. We use the same numbering rule for the
vertices on the outer boundary circle. This numbering of the vertices
combined with the tetrahedral space-time mesh generation algorithm of
\cref{ss:prism2tet} for the rotating space-time inner mesh and the
outer space-time static mesh, results in a `chainsaw' pattern on the
outer boundary of the inner mesh and on the inner boundary of the
outer mesh, see \cref{fig:chainsaw}.

The space-time mesh generated in the annulus now has to match these
chainsaw patterns so that the space-time facets of the outer boundary
of the rotating space-time inner mesh match with the space-time facets
of the inner boundary of the space-time annulus. Likewise, the
space-time facets of the inner boundary of the outer space-time static
mesh have to match with the space-time facets of the outer boundary of
the space-time annulus. The remaining space-time facets of the
tetrahedral space-time mesh in the interior of the annulus now provide
enough freedom to construct a conforming space-time mesh, as described
next.

\begin{figure}[tbp]
  \centering
  \subfloat[Cut 1.]{\includegraphics[width=.22\linewidth]{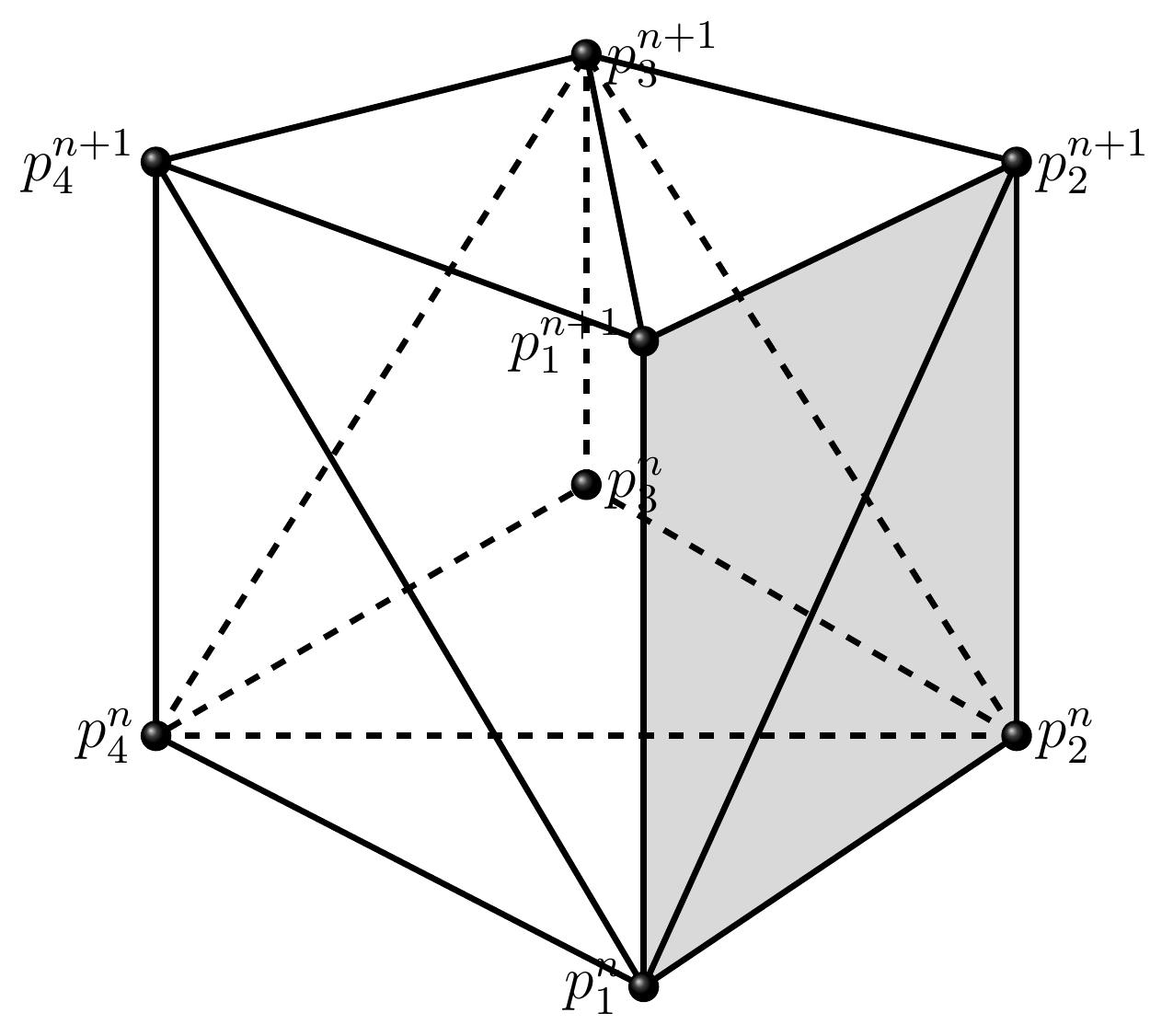}}
  \subfloat[Cut 2.]{\includegraphics[width=.22\linewidth]{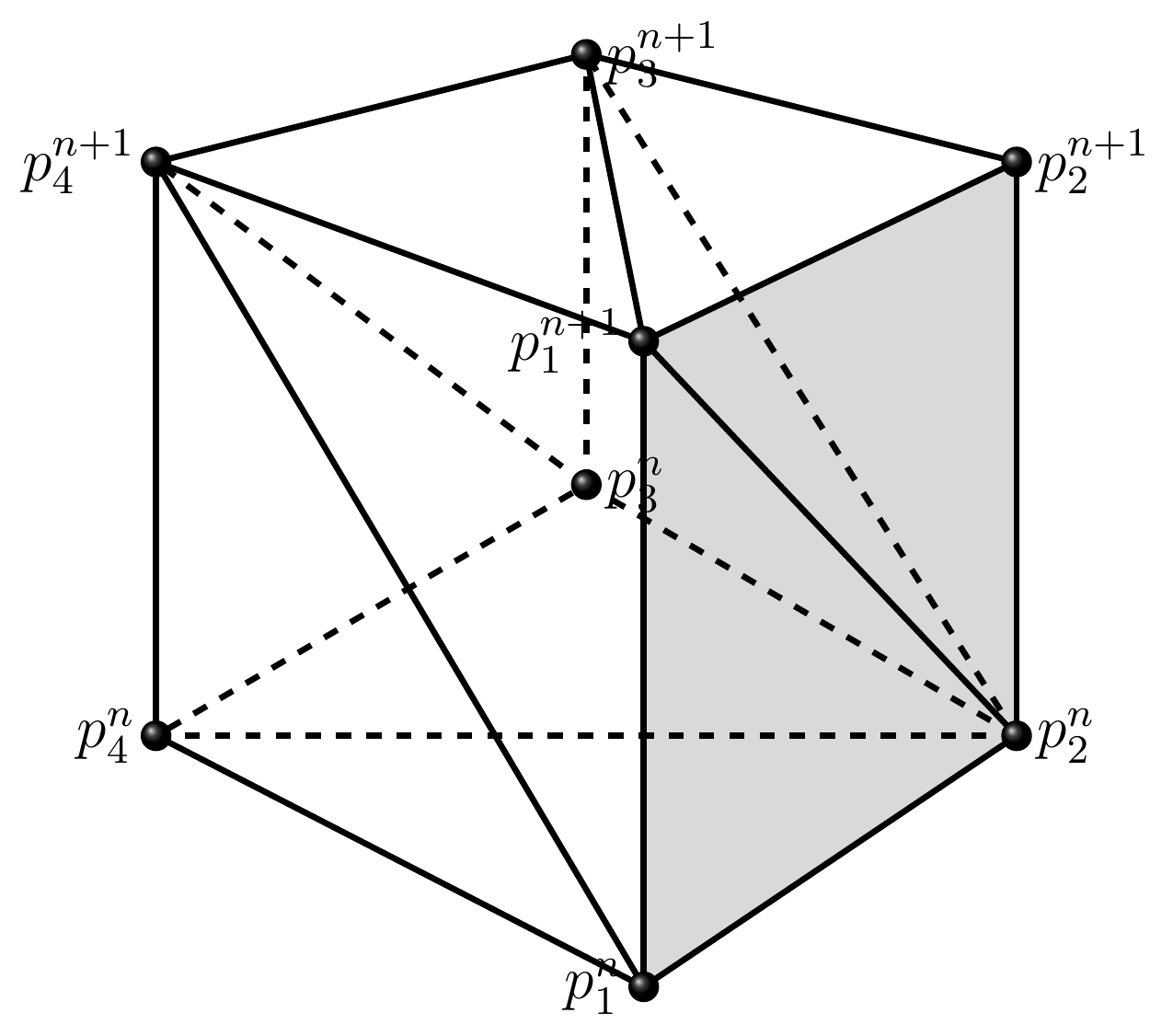}}
  \subfloat[Cut 5.]{\includegraphics[width=.22\linewidth]{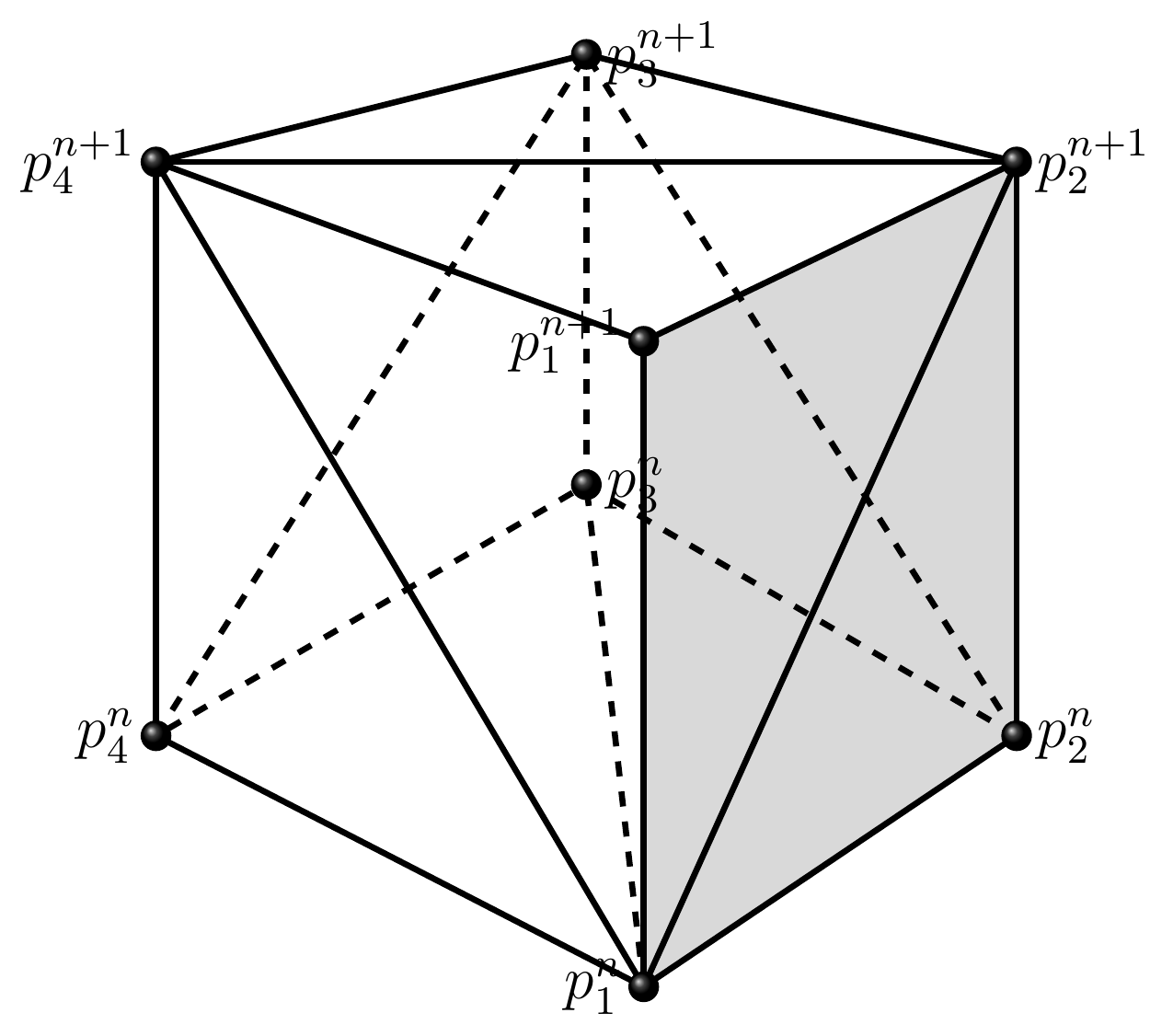}}
  \subfloat[Cut 6.]{\includegraphics[width=.22\linewidth]{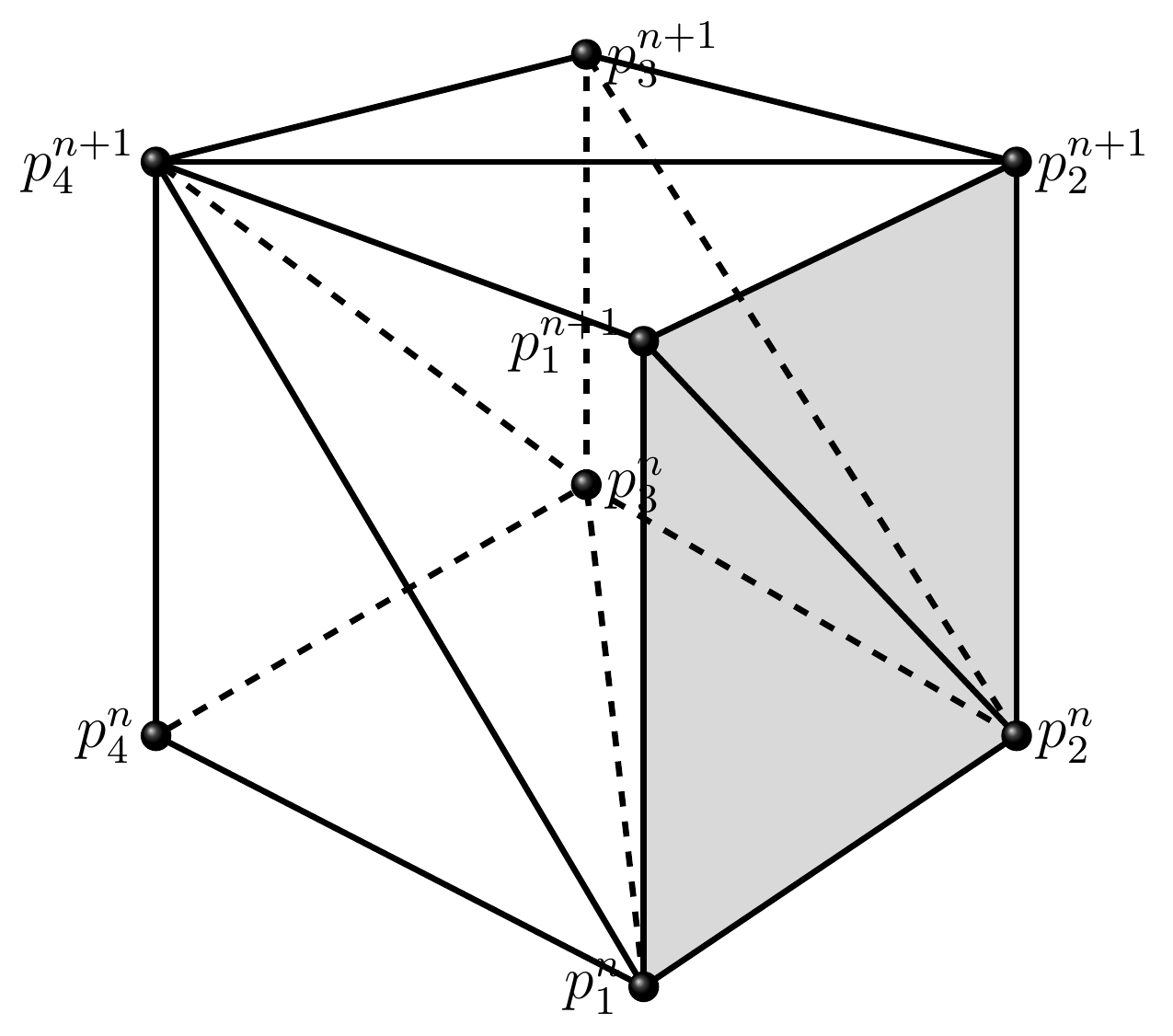}}
  \\
  \subfloat[Cut 3.]{\includegraphics[width=.22\linewidth]{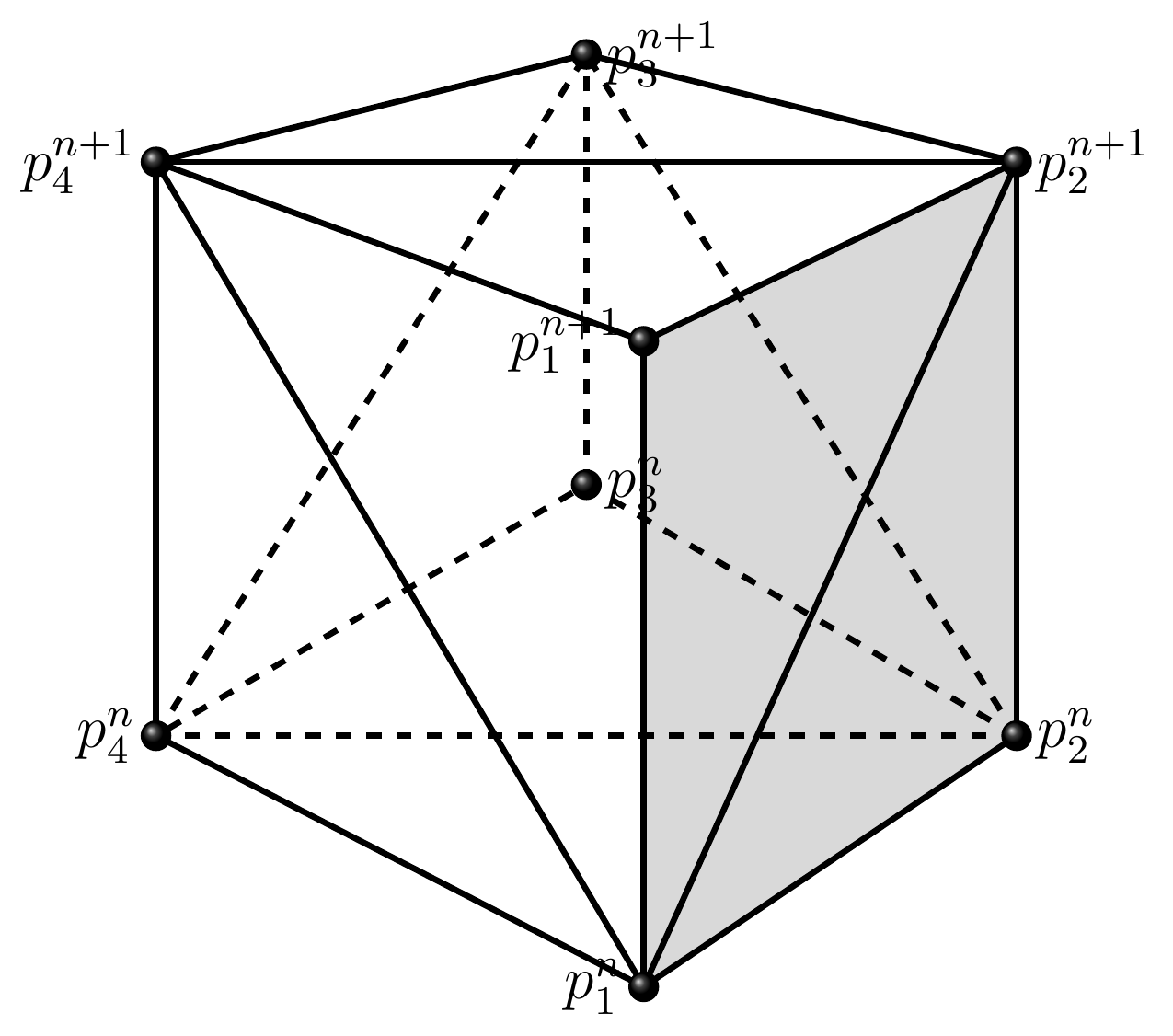}}
  \subfloat[Cut 4.]{\includegraphics[width=.22\linewidth]{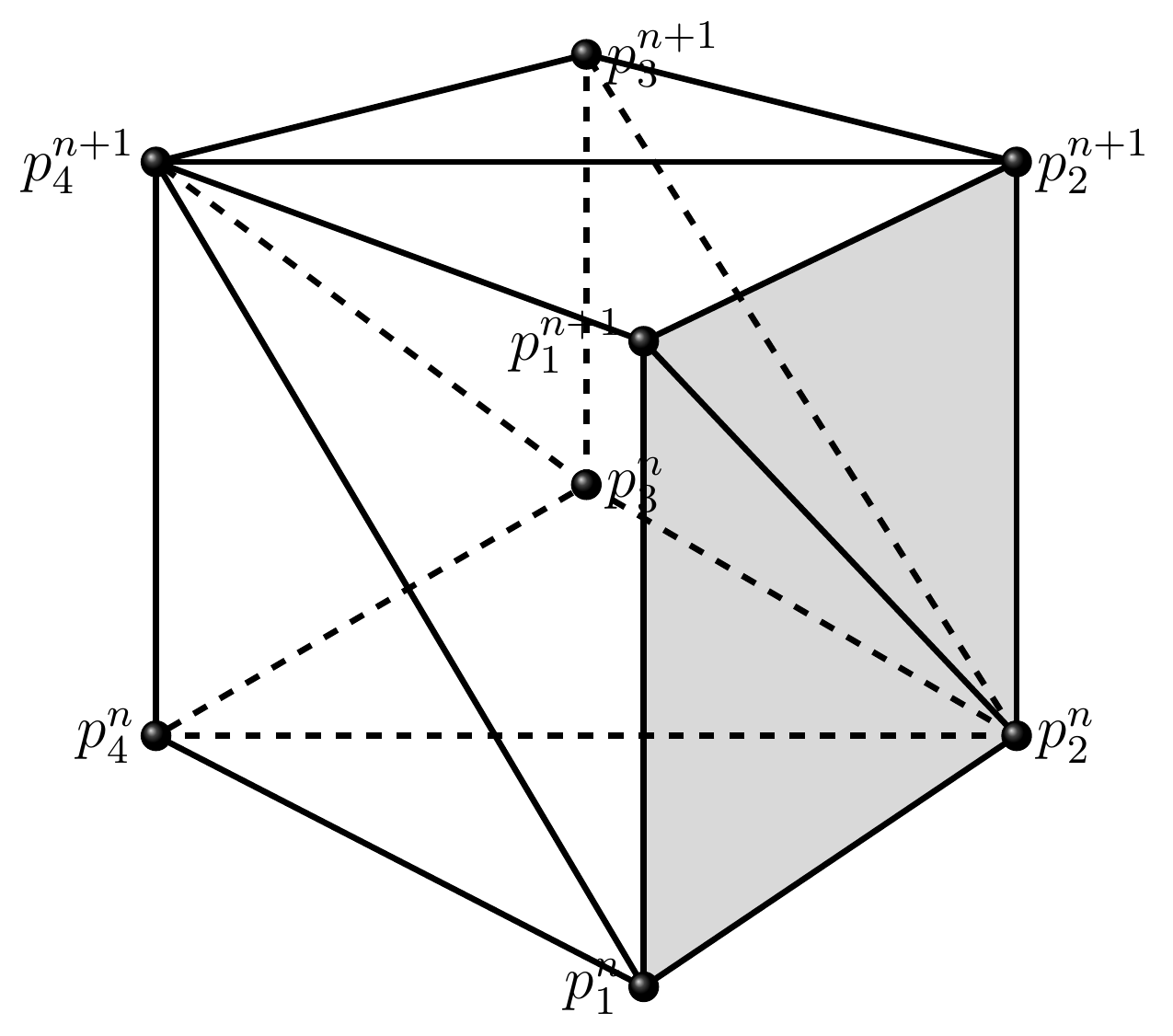}}
  \subfloat[Cut 7.]{\includegraphics[width=.22\linewidth]{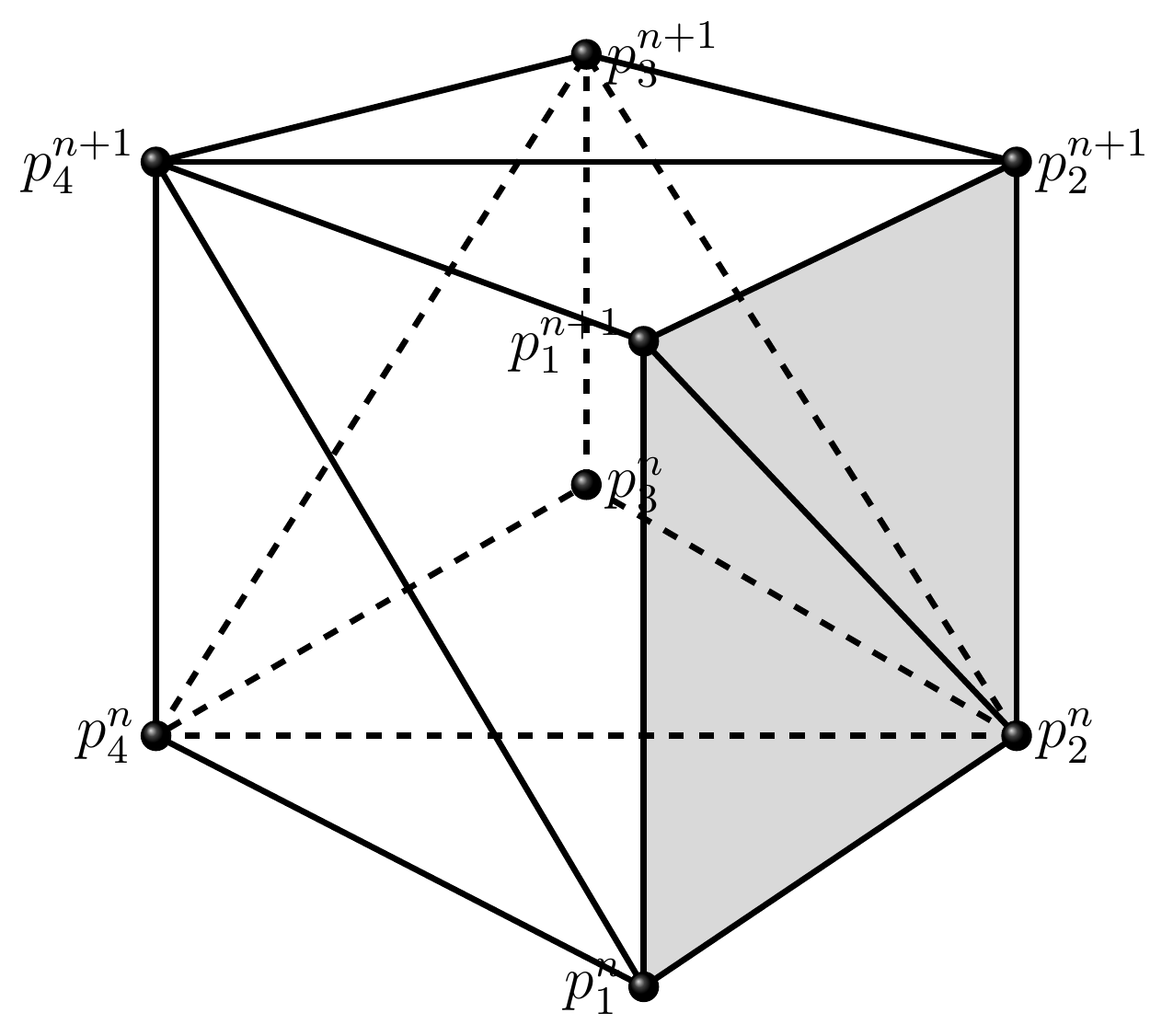}}
  \subfloat[Cut 8.]{\includegraphics[width=.22\linewidth]{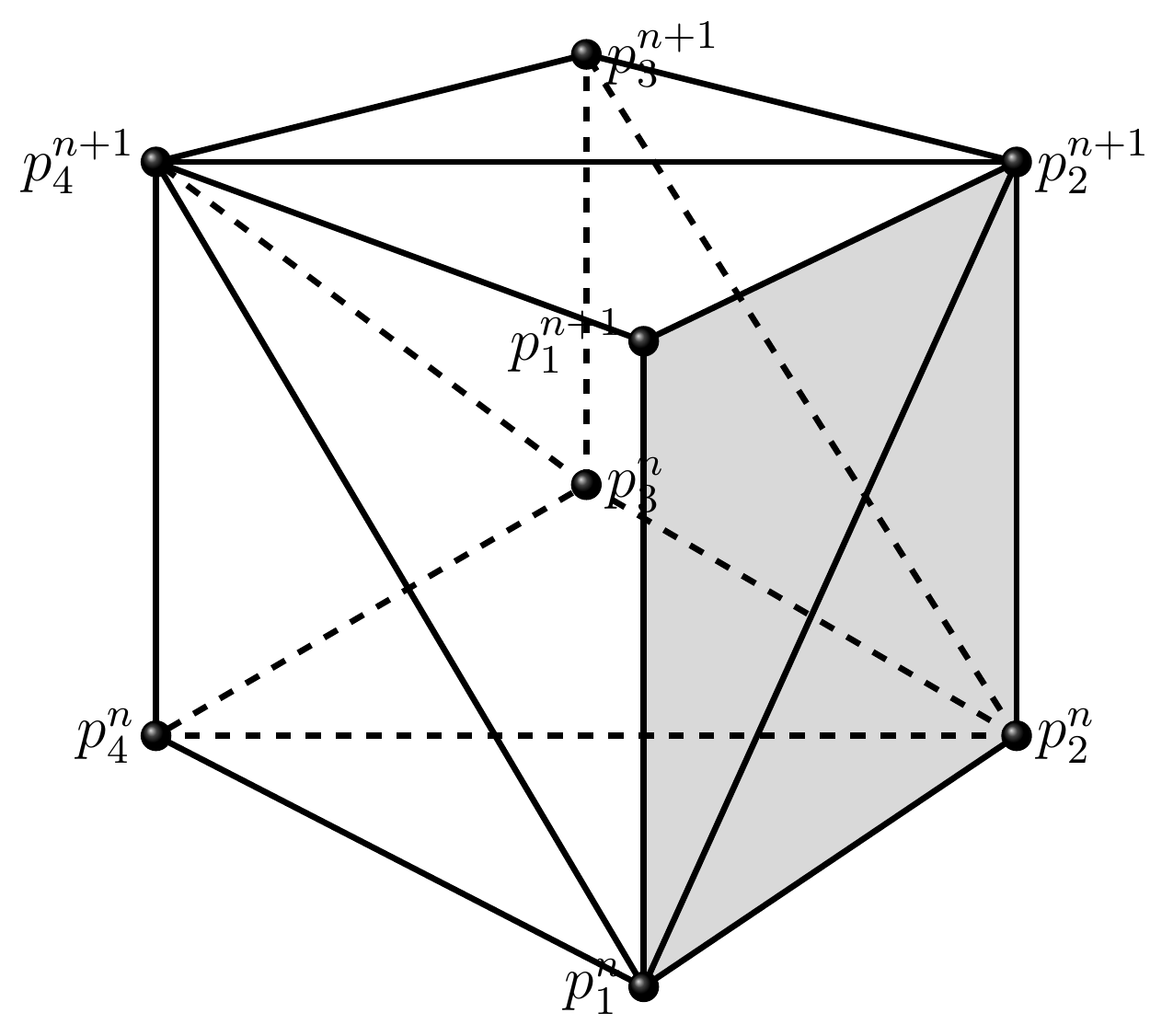}}
  \caption{The gray shaded sides of the above hexahedra are the sides
    closest to the body. Top row: The four different cuttings of a
    hexahedron in the sliding layer when an edge swapping has taken
    place at $t=t^{n+1}$. Bottom row: The four different cuttings of a
    hexahedron in the buffer layer. There is no edge swapping in the
    buffer layer at $t=t^{n+1}$, but the diagonal in the quadrilateral
    side shared by the inner mesh (shaded gray) is either the same as
    the diagonal of the quadrilateral side shared by the sliding layer
    (cuts 3 and 4), or opposite (cuts 7 and 8).}
  \label{fig:quad_cuts}
\end{figure}

Consider the space-time extension from $t^n$ to $t^{n+1}$ of a
$2 \times 2$ block of quadrilaterals in the annulus of which two are
adjacent elements in the sliding layer and the other two are their
neighbors in the buffer layer. Note that this space-time extension may
not result in a space-time hexahedra as not all faces may be flat
space-time quadrilaterals. However, similar to the algorithm described
in \cref{ss:prism2tet}, we never create a space-time hexahedra; it is
used only as an intermediate step to create space-time tetrahedra. To
simplify the discussion we therefore refer to this space-time
extension as a space-time hexahedra with space-time quadrilateral
sides. To continue, let us denote the space-time quadrilateral side of
a space-time hexahedron in the buffer layer that is shared by the
inner mesh by $S_{in}$ (shaded gray in the bottom row of
\cref{fig:quad_cuts}) and the space-time quadrilateral side of a
space-time hexahedron that is shared by the sliding layer by
$S_{sl}$. There are now four different cuts for each space-time
hexahedron in the sliding layer and four different cuts for each
space-time hexahedron in the buffer layer:
\begin{enumerate}
\item In the sliding layer at $t=t^n$ the quadrilateral is cut along
  the $n_2^E - n_4^E$ diagonal while it is cut along the
  $n_1^E - n_3^E$ diagonal at $t=t^{n+1}$, resulting in cuts 1 or 2 in
  \cref{fig:quad_cuts}. Or, in the sliding layer at $t=t^n$ the
  quadrilateral is cut along the $n_1^E - n_3^E$ diagonal while it is
  cut along the $n_2^E - n_4^E$ diagonal at $t=t^{n+1}$, resulting in
  cuts 5 or 6 in \cref{fig:quad_cuts}.
\item In the buffer layer, the diagonal cut of $S_{in}$ is the same as
  the diagonal cut of $S_{sl}$, resulting in cuts 3 or 4 in
  \cref{fig:quad_cuts}. Or, in the buffer layer, the diagonal cut of
  $S_{in}$ is opposite to the diagonal cut of $S_{sl}$, resulting in
  cuts 7 or 8 in \cref{fig:quad_cuts}.
\end{enumerate}
To obtain a conforming mesh we note that all diagonals of the
space-time quadrilateral sides of a space-time hexahedron have to
match with the diagonals of the space-time quadrilateral sides of an
adjacent space-time hexahedron. Furthermore, the diagonals of the
space-time quadrilateral sides of a space-time hexahedron have to
match with the `chainsaw' pattern of the outer boundary of the
rotating space-time inner mesh and on the inner boundary the outer
space-time static mesh. This results in four configurations of the
$2\times 2$ block using the above eight cuts, see
\cref{tab:quad_cuts}. Consider, for example, the first configuration
of \cref{tab:quad_cuts}. Then, in the sliding layer, Cut 1 is adjacent
to Cut 2. Cut 1 and Cut 2 are furthermore adjacent to Cuts 3 and 4,
respectively, in the buffer layer. Configurations 2, 3, and 4 are
defined similarly.

\begin{table}[tbp]
  \centering
  \begin{tabular}{c|ccc|ccc|ccc|c}
    Cut 1 & Cut 2 & & Cut 5 & Cut 6 & & Cut 1  & Cut 2 & & Cut 5  & Cut 6  \\\cline{1-2}\cline{4-5}\cline{7-8}\cline{10-11}
    Cut 3 & Cut 4 & & Cut 3 & Cut 4 & & Cut 7 & Cut 8 & & Cut 7 & Cut 8 \\
    \multicolumn{11}{c}{}\\
    \multicolumn{2}{c}{Configuration 1} && \multicolumn{2}{c}{Configuration 2} && \multicolumn{2}{c}{Configuration 3} && \multicolumn{2}{c}{Configuration 4}
  \end{tabular}
  \caption{The four different hexahedral cut configurations to create
    a conforming tetrahedral mesh over the sliding and buffer
    layers. See \cref{fig:quad_cuts} for the different cuts.}
  \label{tab:quad_cuts}
\end{table}

\captionsetup[subfigure]{labelformat=empty}

\section{Numerical examples}
\label{s:examples}

In this section we consider various test cases of increasing
complexity and compare, where possible, our results to those found in
the literature. For this we added the fluid-rigid body interaction and
conforming sliding mesh technique algorithms to our space-time EHDG
implementation of the Navier--Stokes equations as described previously
in \cite{Horvath:2020}. In particular, all simulations in this section
have been implemented in the Modular Finite Element Method (MFEM)
library \cite{mfem-library}. We use the direct solver of
MUMPS~\cite{MUMPS:1,MUMPS:2} through
PETSc~\cite{petsc-web-page,petsc-user-ref,petsc-efficient} to solve
the linear systems resulting from the Picard iteration
\cref{eq:Picard}. Unless specified differently, we choose as
tolerances for the Picard iteration \cref{eq:stopping_criterion} and
the rigid body stopping criterion \cref{eq:stopping_crit_rb},
respectively, as $\delta_{\text{NS}} = 10^{-6}$ and
$\delta_{\text{rb}} = 10^{-5}$. Finally, the penalty parameter
$\alpha$ in \cref{eq:numflux} is chosen as $\alpha = 6k^2$ where $k$
is the polynomial degree used in
\cref{eq:velocity_spc,eq:pressure_spc} (see \cite{Riviere:book} for
more information on the penalty parameter). In all simulations we take
$k=2$.

\subsection{Vortex induced oscillations and galloping}
\label{ss:transv_gall}

\begin{figure}[tbp]
  \centering
  \subfloat[Domain setup.]{\includegraphics[width=.5\linewidth]{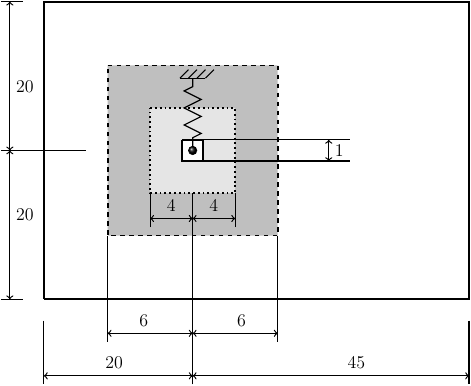}}
  \qquad
  \subfloat[Zoom of the spatial mesh near the body.]{\includegraphics[width=.4\linewidth]{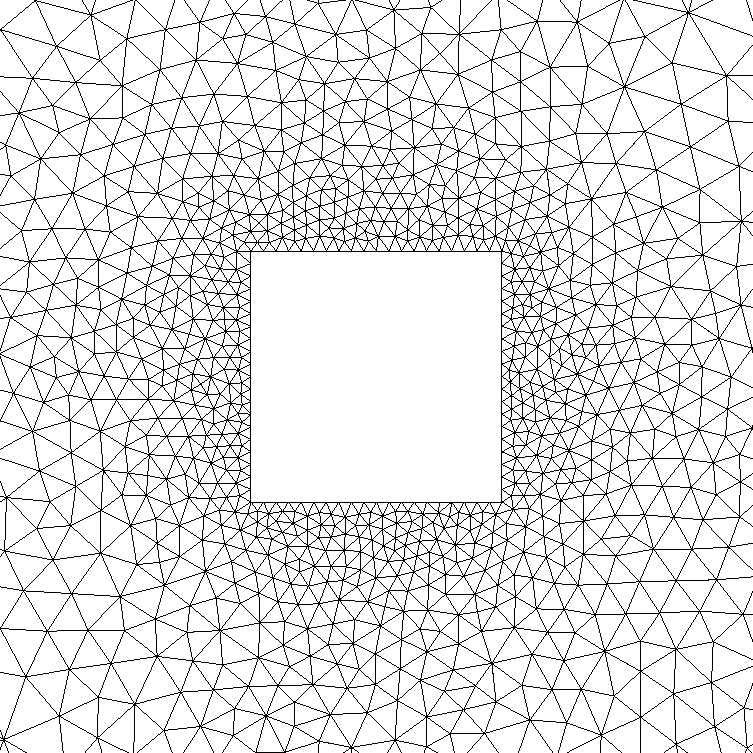}}
  \caption{Left: depiction of the problem setup as described
    in~\cref{ss:transv_gall}. Right: a zoom of the spatial mesh near
    the body.}
  \label{fig:transv_gall_geo}
\end{figure}

In this example, we test the coupling algorithm between the fluid and
the rigid body, but we do not yet consider the sliding mesh technique
of \cref{s:sliding_grid}. To test this coupling we consider vertical
galloping only, i.e., we only solve for the displacement $d$ in
\cref{eq:linearMotion}. We follow the setup of this test case as
described previously in \cite{Robertson:2003, Dettmer:2006,
  Kadapa:2017}. For this, consider a square object with sides of
length 1 in a rectangular fluid domain with outer boundary sides of
length 65 and 40. At time $t=0$ we place the center of the object at a
distance of 20 from the left and bottom outer fluid domain boundary
(see \cref{fig:transv_gall_geo}).

The coefficients of the rigid-body spring-system are chosen as
$m = 20$, $c_y = 0.00581195$, $k_y = 3.08425$. In the Navier--Stokes
equations \cref{eq:navierstokes} we choose the viscosity as
$\nu = 0.01$, while the density in \cref{eq:forcemoment} is chosen as
$\rho = 1$. We impose $\boldsymbol{u}=(2.5, 0)$ on the left outer
boundary of the fluid domain, free-slip on the top and bottom outer
boundaries, and a homogeneous Neumann boundary condition on the right
outer boundary. The spatial mesh consists of 8510 triangles with 80
elements on the boundary of the body (see \cref{fig:transv_gall_geo}
for a zoom of the mesh near the body). We set the time step equal to
$\Delta t = 0.1$.

Since we consider the motion of a rigid body without rotation, we
update our mesh in a similar manner as we used previously in
\cite[Section 5.4]{Horvath:2020}: nodes inside the spatial box
$\Omega^{\text{in}}(t)=[-4,\ 4] \times [-4 + d(t),\ 4+d(t)]$ (the
light grey area in \cref{fig:transv_gall_geo}) move with the body (the
position of each node in $\Omega^{\text{in}}(t)$ does not change with
respect to the body); the movement of the nodes in
$\Omega^{\text{mid}}(t)=([-6,\ 6] \times [-6,\ 6]) \backslash
\Omega^{\text{in}}(t)$ (the dark grey area in
\cref{fig:transv_gall_geo}) decreases linearly with distance from the
body; and the remaining nodes are fixed.

\begin{figure}[tbp]
  \begin{center}
    \includegraphics[width=\linewidth]{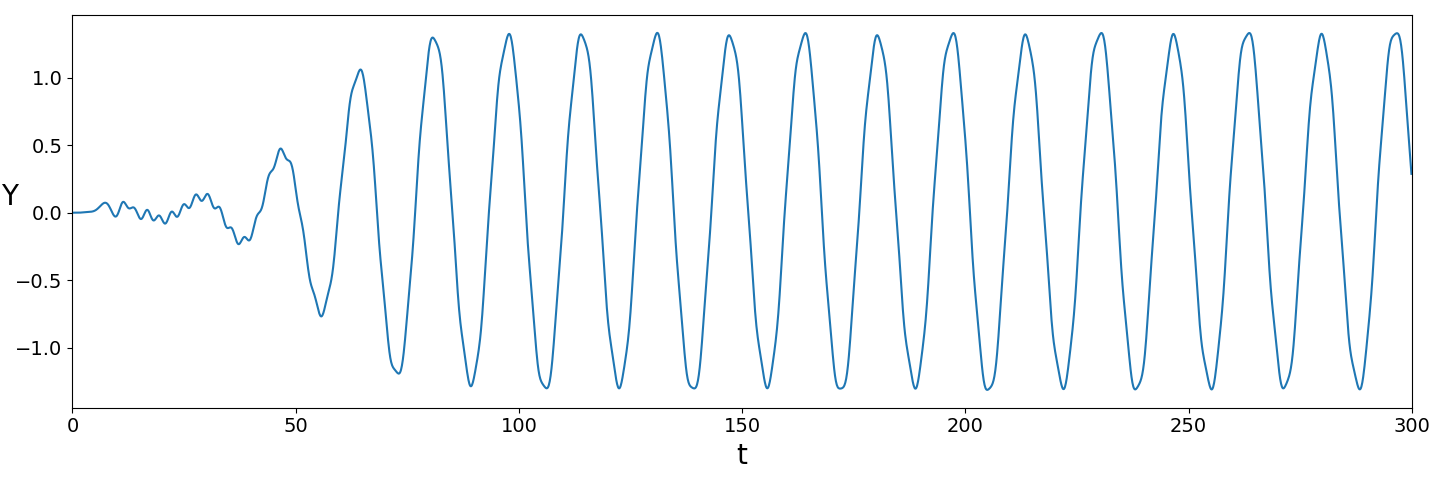} 
    \caption{Evolution of the amplitude for the problem described
      in~\cref{ss:transv_gall}}
    \label{fig:transv_gall_amplitude}
  \end{center}
\end{figure}

\begin{figure}[tbp]
  \begin{center}
    \includegraphics[width=0.3\linewidth]{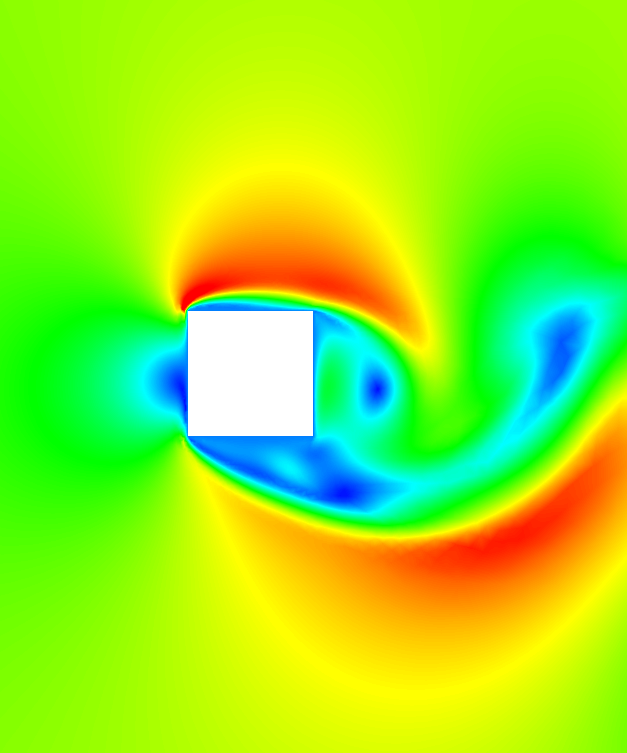} 
    \includegraphics[width=0.3\linewidth]{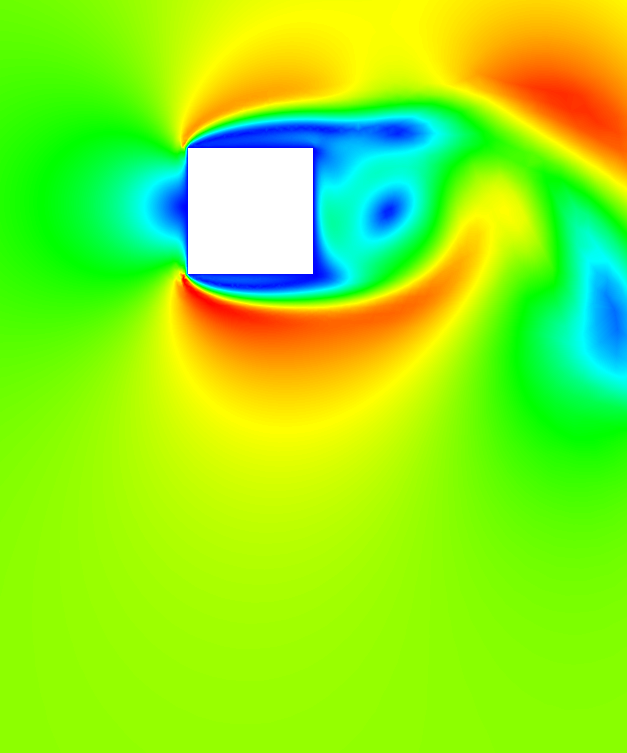} 
    \includegraphics[width=0.3\linewidth]{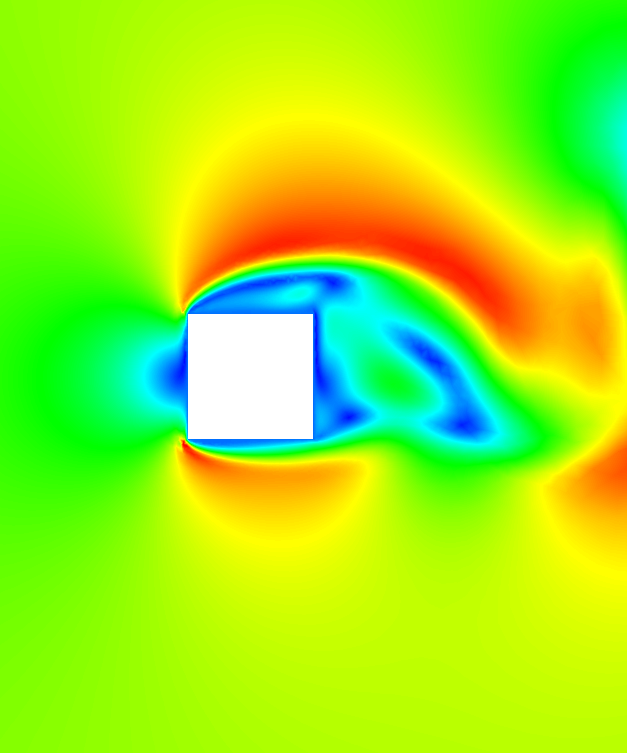}\\
    \includegraphics[width=0.3\linewidth]{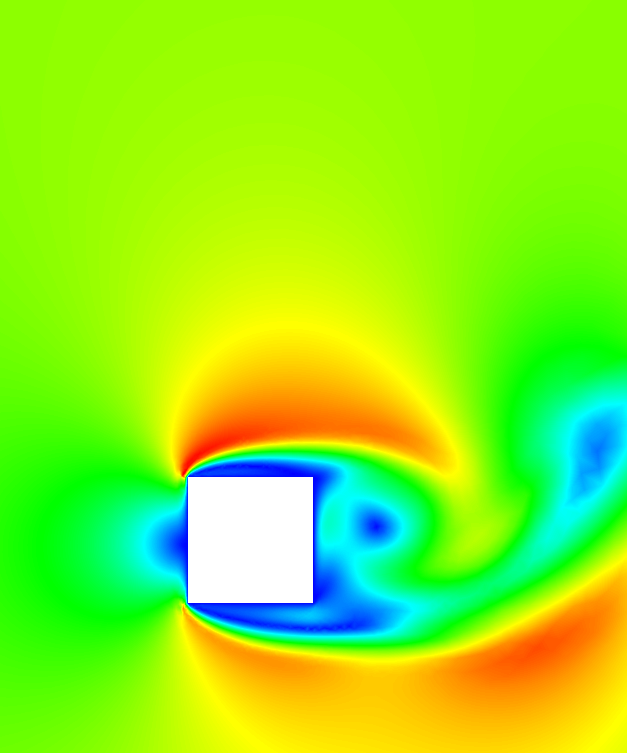} 
    \includegraphics[width=0.3\linewidth]{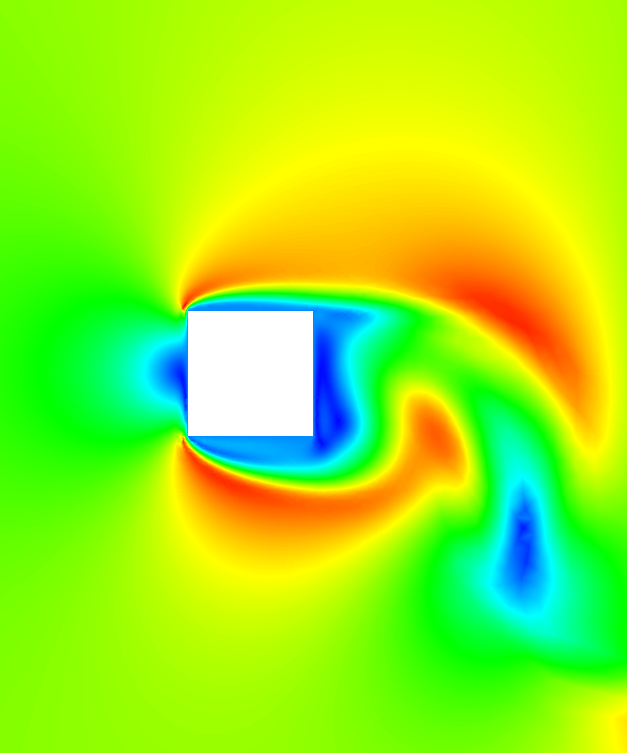} 
    \includegraphics[width=0.3\linewidth]{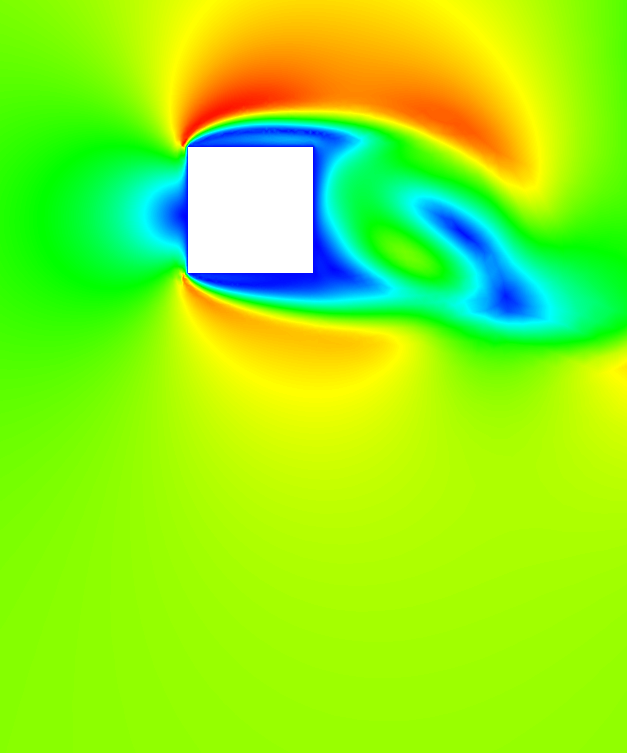}\\
    \includegraphics[width=0.3\linewidth]{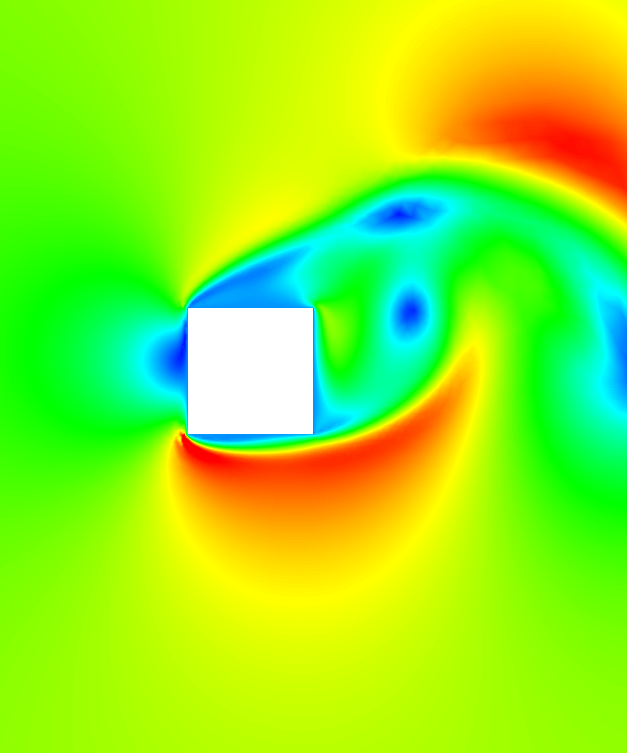} 
    \includegraphics[width=0.3\linewidth]{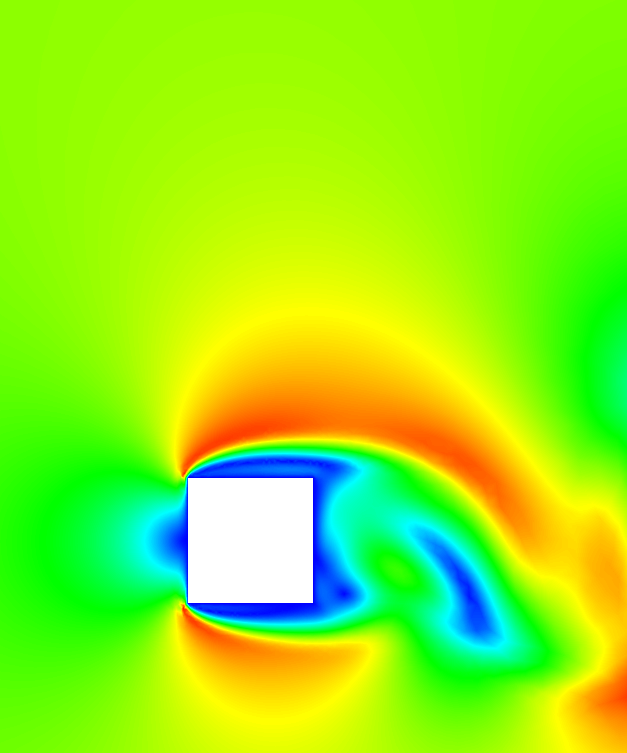} 
    \includegraphics[width=0.3\linewidth]{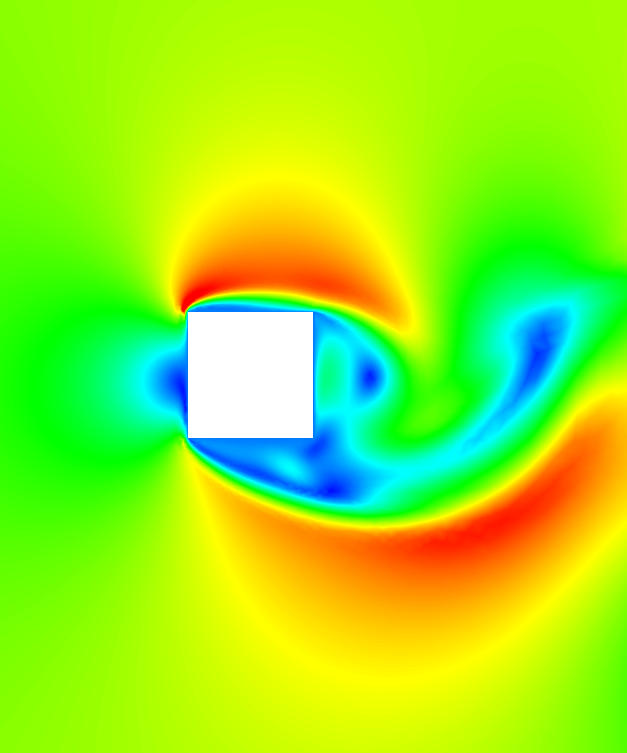} 
    \caption{The velocity magnitude at different, evenly spaced,
      points in time during two up-down periods for the problem
      described in~\cref{ss:transv_gall}.}
    \label{fig:transv_gall_velocity_evolution}
  \end{center}
\end{figure}

We plot the evolution of the amplitude of the motion of the rigid body
in \cref{fig:transv_gall_amplitude} for $t \in [0, 300]$. We find a
maximum amplitude of $\max|y| = 1.3$ and a frequency of
$0.62$. Compared to results found in the literature, in which the
maximum amplitude lies in the interval $[1.08, 1.2]$ and the frequency
lies in the interval $[0.058, 0.069]$ (see, e.g.,
\cite{Robertson:2003, Dettmer:2006, Kadapa:2017}), we note that we
slightly over-predict the maximum amplitude, but that the frequency
lies in the same range as other works. The velocity magnitude, at
different evenly spaced points in time during two up-down periods, is
shown in \cref{fig:transv_gall_velocity_evolution}.

\subsection{Rotational galloping}
\label{ss:rectangle}

In this section we combine the conforming sliding mesh technique
described in~\cref{s:sliding_grid} with fluid-rigid body interaction
by simulating the rotational galloping of a rectangle. The rotational
galloping test case, in which there is no vertical motion, has also
been considered in, for example, \cite{Robertson:2003,Dettmer:2006}.

\begin{figure}[tbp]
  \begin{center}
    \includegraphics[width=.6\linewidth]{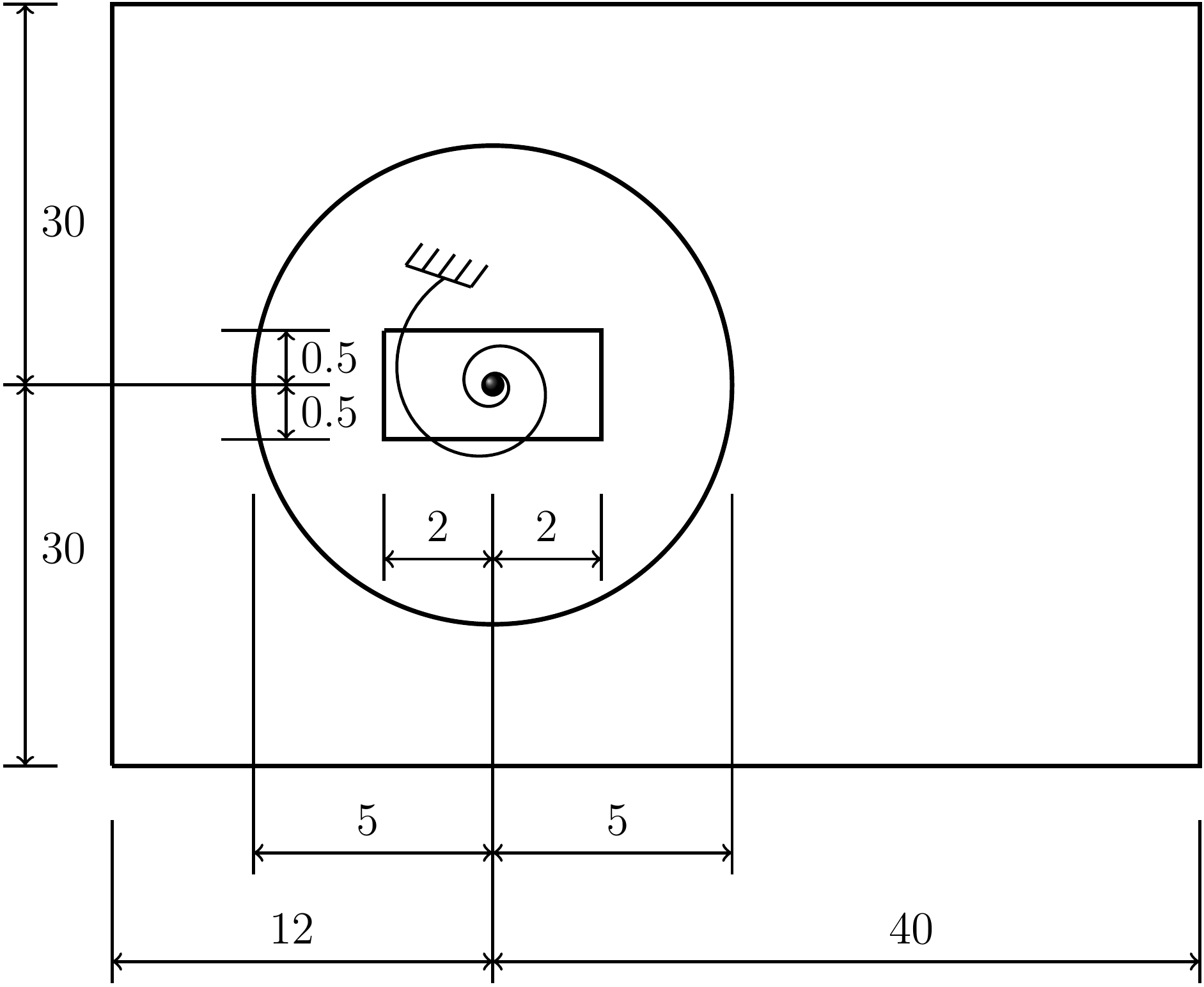} 
    \caption{Set up of the problem of~\cref{ss:rectangle}.}
    \label{fig:galloping_geo}
  \end{center}
\end{figure}

We consider a rectangular fluid domain of dimensions 52 by 60. The
center of a rectangular rigid body, with dimensions 4 by 1, is placed
at a distance 30 from the outer bottom wall of the fluid domain and a
distance 12 from the outer left wall. See \cref{fig:galloping_geo} for
a setup of the geometry. The coefficients of the rigid body spring
system \cref{eq:linearMotion_theta} are $I_\theta = 400$,
$c_\theta = 78.54$, and $k_\theta = 61.685$. For the fluid problem we
again choose the viscosity and density as, respectively, $\nu = 0.01$
and $\rho = 1$. We impose $\boldsymbol{u} = (10, 0)$ on the left outer
boundary of the fluid domain, free-slip on the top and bottom outer
boundaries, and a homogeneous Neumann boundary condition on the right
outer boundary. We set the time step equal to $\Delta t = 0.15$.

\begin{figure}[tbp]
  \begin{center}
    \includegraphics[width=0.49\linewidth]{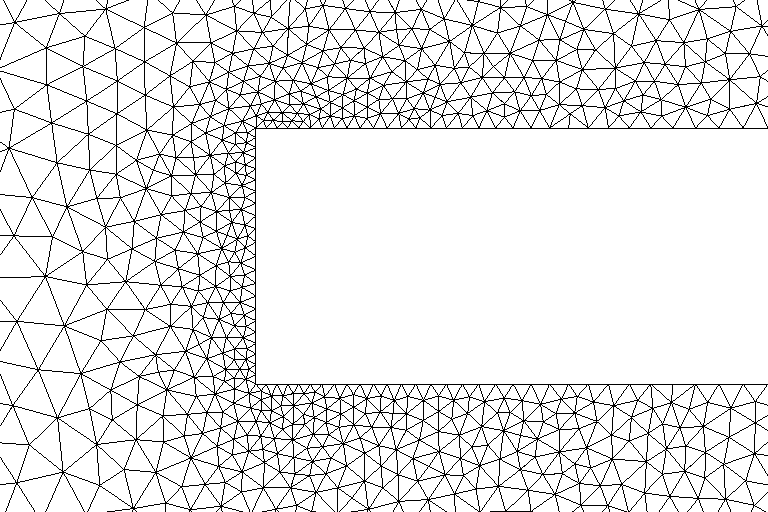} 
    \includegraphics[width=0.49\linewidth]{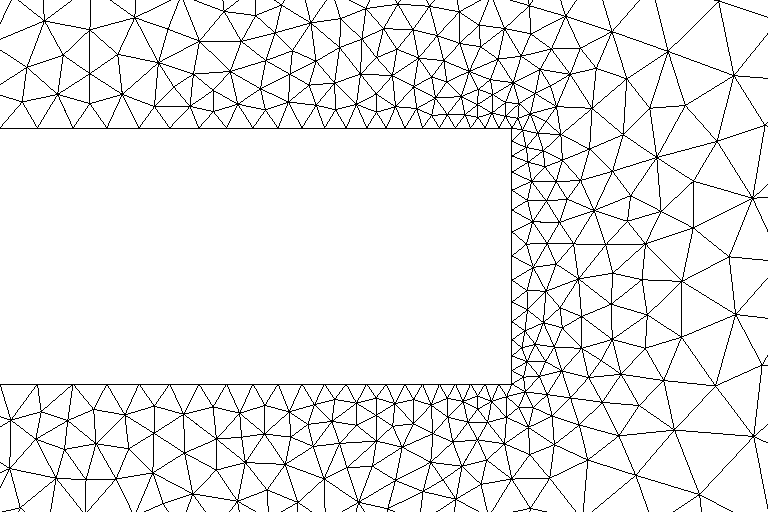} 
    \caption{Fine (left) and coarse (right) mesh near the body for the
      problem described in~\cref{ss:rectangle}.}
    \label{fig:galloping_mesh}
  \end{center}
\end{figure}

In this simulation we compute results on a coarse mesh, consisting of
6030 triangles with 116 elements on the boundary of the body, and on a
fine mesh with \num{13722} triangles with 168 body boundary elements
(see \cref{fig:galloping_mesh} for a zoom of the mesh near the
body). The inner and outer radii of the annulus consisting of the
sliding and buffer layer are placed at a distance of, respectively,
4.4 and 5 from the center of the rigid body. This annulus is initially
(before creating the space-time tetrahedral mesh) tessellated into 100
quadrilaterals for the coarse mesh and 152 for the fine mesh.

\begin{figure}[tbp]
  \begin{center}
    \includegraphics[width=\linewidth]{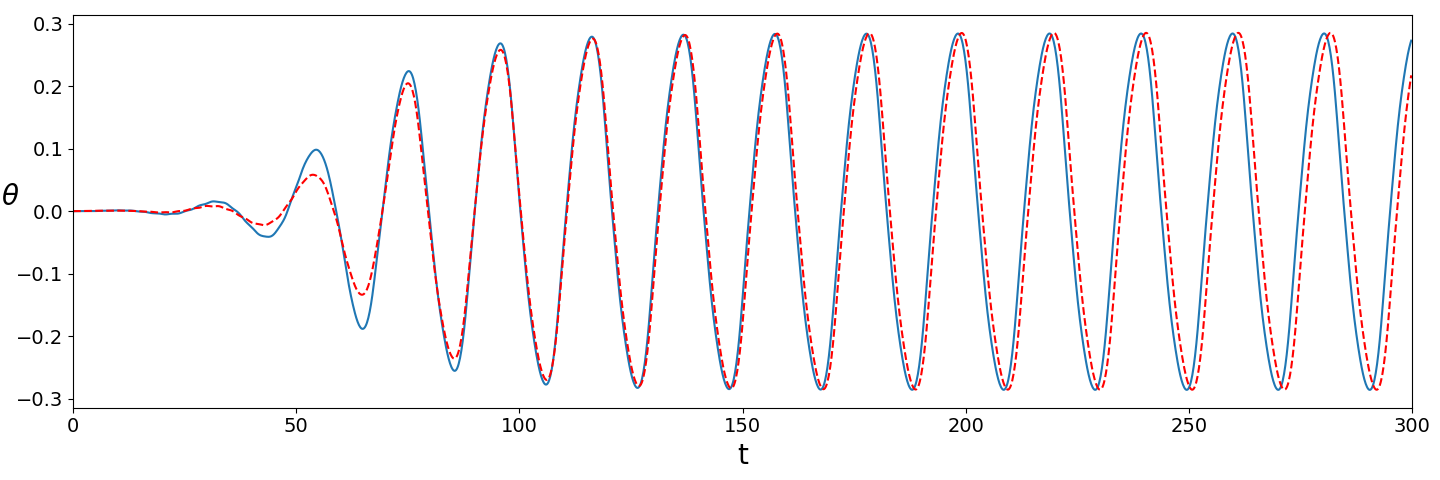} 
    \caption{Evolution of the rotational angle computed on the coarse
      mesh (blue) and the fine mesh (red). See \cref{ss:rectangle}.}
    \label{fig:galloping_angle}
  \end{center}
\end{figure}

The evolution of the rotational angle for $t \in [0, 300]$ is shown in
\cref{fig:galloping_angle}. On the coarse mesh we find a maximum
amplitude of $\max|y| = 0.286$ and a frequency of $0.048$. On the fine
mesh we find a maximum amplitude of $\max|y| = 0.285$ and a frequency
of $0.049$. These results compare well to results found in the
literature. Indeed, \cite{Robertson:2003, Dettmer:2006, Kadapa:2017}
predict the maximum amplitude to lie in the interval $[0.223, 0.283]$
and the frequency to lie in the interval $[0.048, 0.052]$.

\begin{figure}[tbp]
  \begin{center}
    \includegraphics[width=0.45\linewidth]{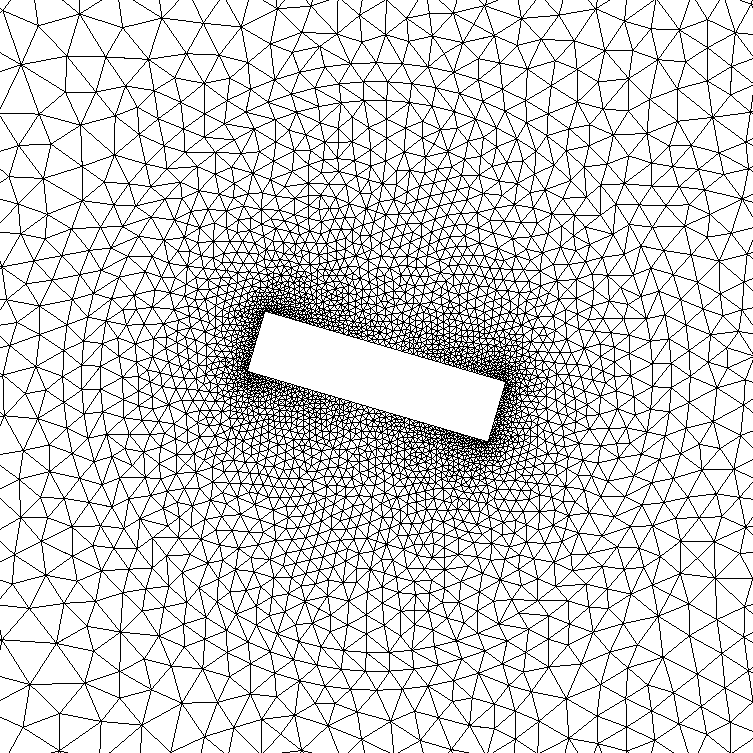} 
    \includegraphics[width=0.45\linewidth]{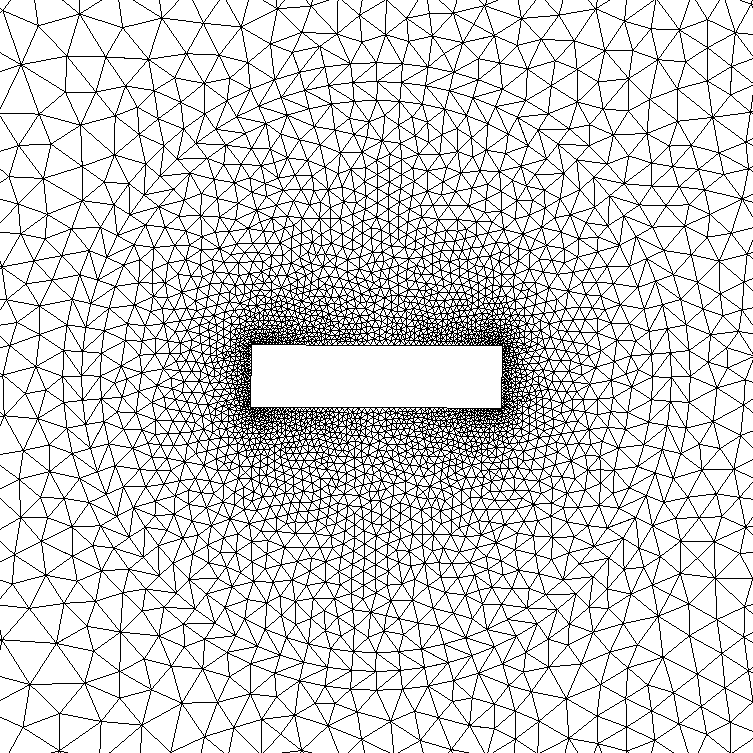} \\
    \includegraphics[width=0.45\linewidth]{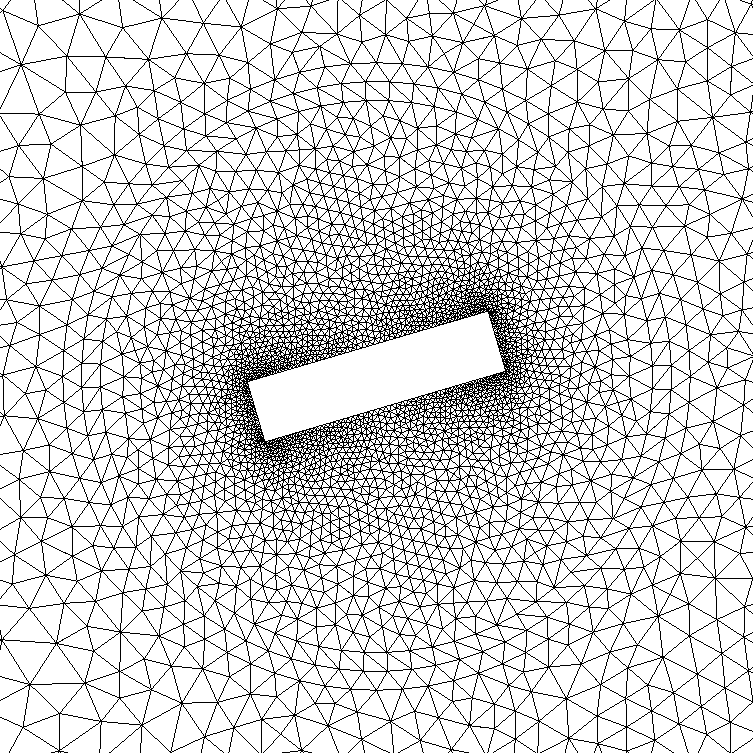} 
    \caption{The mesh at the maximum, zero, and minimum angle of the
      body for the problem described in~\cref{ss:rectangle}.}
    \label{fig:rot_gal_mesh_evolution}
  \end{center}
\end{figure}

\begin{figure}[tbp]
  \begin{center}
    \includegraphics[width=0.45\linewidth]{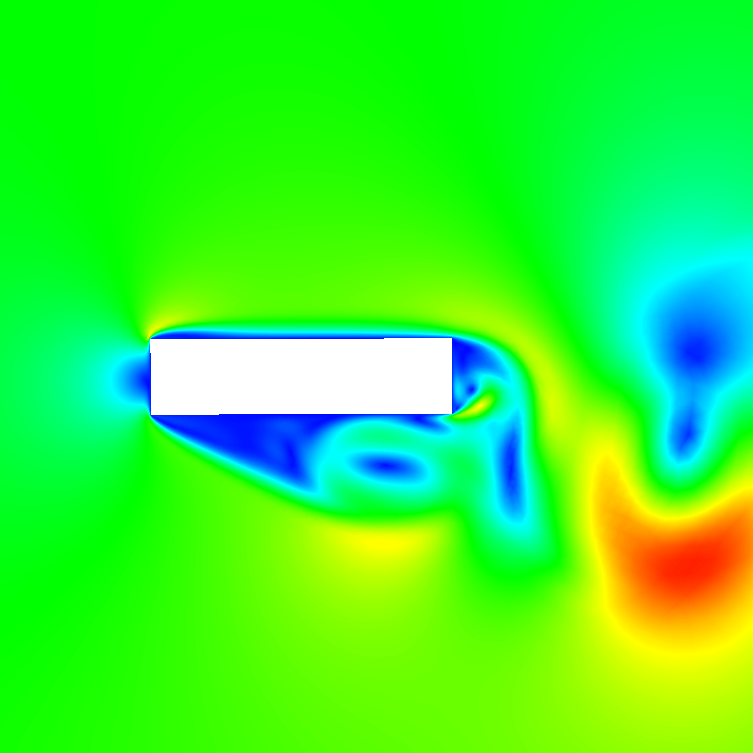} 
    \includegraphics[width=0.45\linewidth]{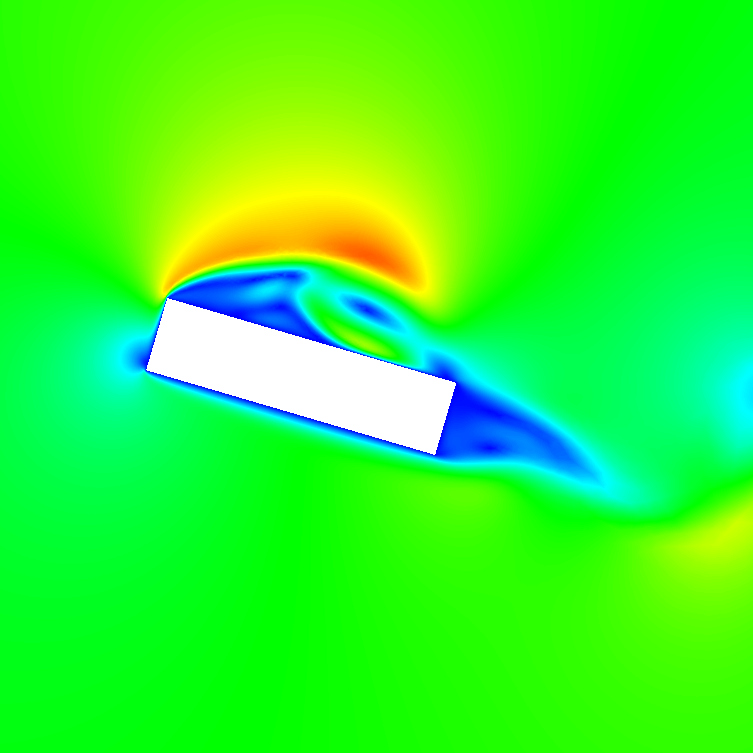} \\
    \includegraphics[width=0.45\linewidth]{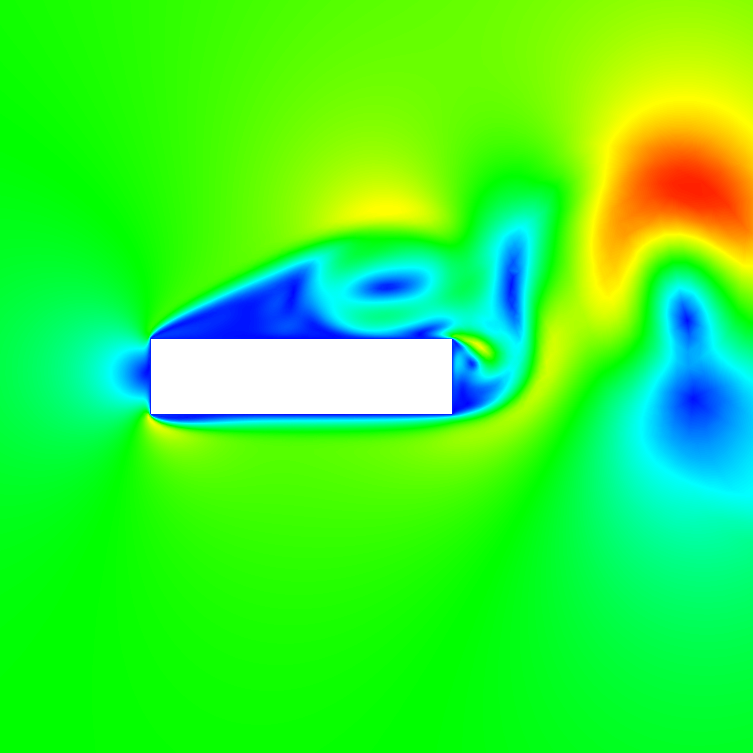}
    \includegraphics[width=0.45\linewidth]{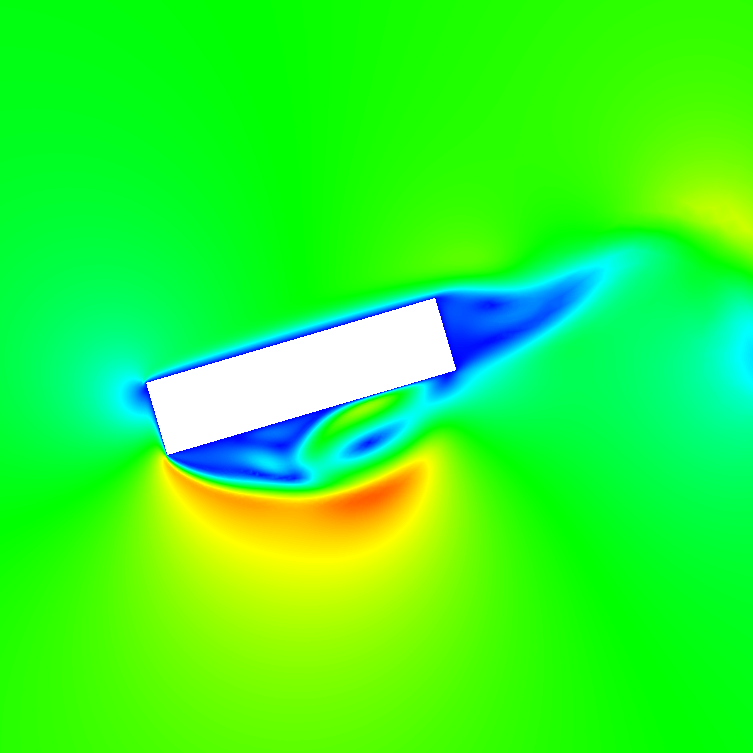} \\
    \includegraphics[width=0.45\linewidth]{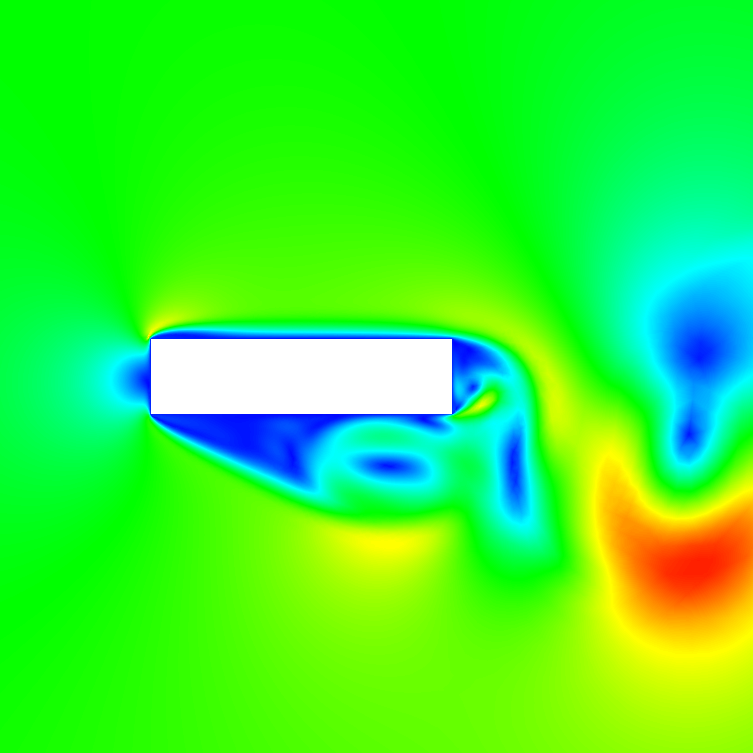} 
    \caption{The velocity magnitude computed on the fine mesh at
      different, evenly spaced, points in time during one period for
      the problem described in~\cref{ss:rectangle}.}
    \label{fig:rot_gal_velocity_evolution}
  \end{center}
\end{figure}

A plot of the fine mesh at maximum, zero, and minimum angle of the
body is given in \cref{fig:rot_gal_mesh_evolution}. Observe that the
connectivity changes only in the sliding layer. Finally, we plot in
\cref{fig:rot_gal_velocity_evolution} the velocity magnitude computed
on the fine mesh at different, evenly spaced, points in time during
one period.

\subsection{H-shaped bridge}
\label{ss:bridge}

We now consider the fluttering H-shaped bridge example of \cite[Section
8.4]{Dettmer:2006}. In particular, we consider a fluid-rigid body
interaction simulation in which both rotational and vertical
displacement of the body are computed for.

\begin{figure}[tbp]
  \begin{center}
    \includegraphics[width=\linewidth]{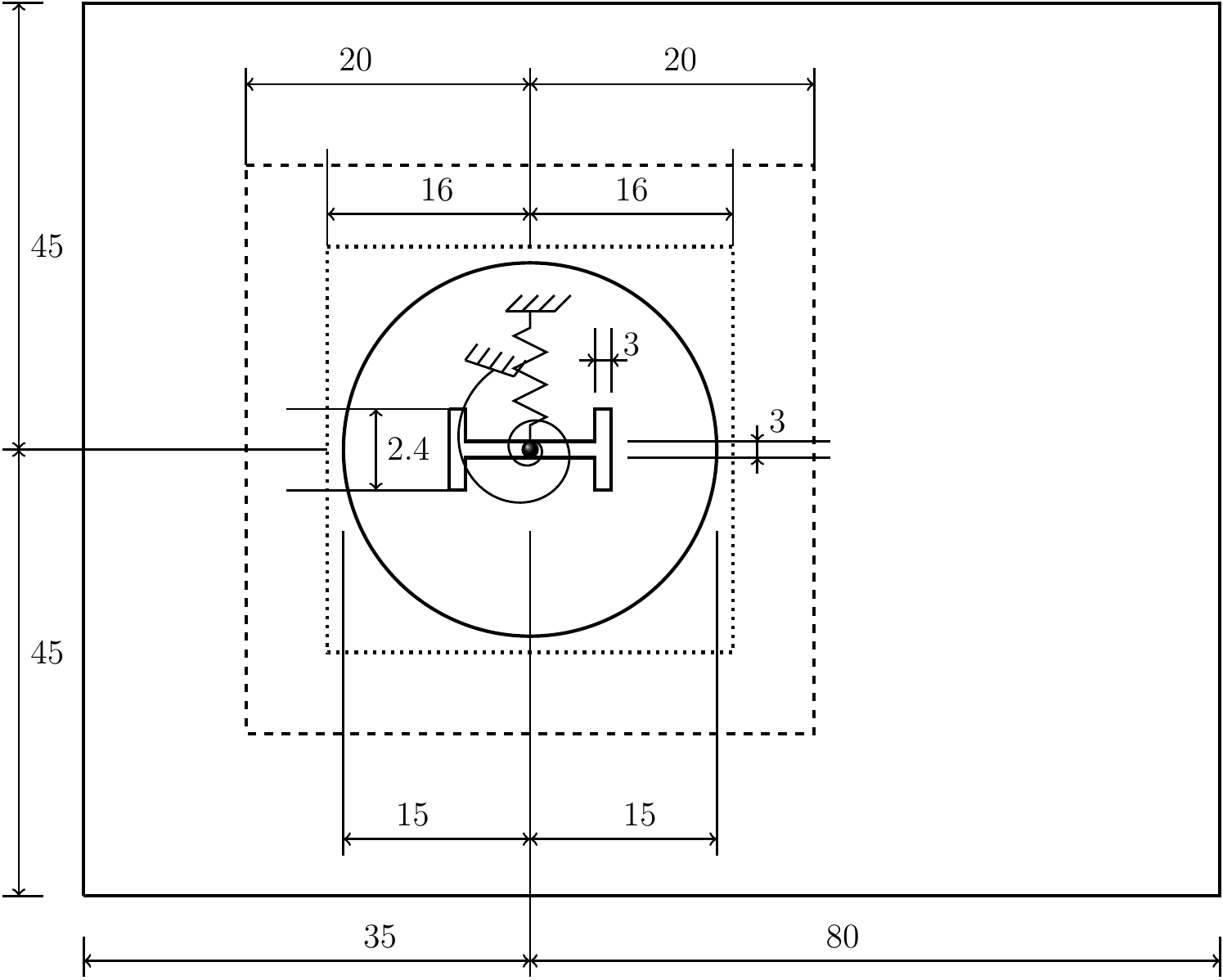} 
    \caption{Set up of the fluttering H-shaped bridge example
      of~\cref{ss:bridge}.}
    \label{fig:fluttering_geo}
  \end{center}
\end{figure}

The setup of the geometry is depicted
in~\cref{fig:fluttering_geo}. The coefficients of the rigid body
spring system \cref{eq:linearMotion} are chosen as
$I_\theta = \num{25300}$, $c_\theta = 0$, $k_\theta = \num{40000}$,
$m = \num{3000}$, $c_y = 0$, and $k_y = \num{2000}$. For the fluid
problem we choose $\nu = 0.1$ and $\rho = 1.25$. We impose
$\boldsymbol{u} = (10,0)$ on the left boundary, free-slip on the top
and the bottom boundaries, and a homogeneous Neumann boundary
condition on the right boundary. As time step we choose
$\Delta t = 0.05$.

\begin{figure}[tbp]
  \begin{center}
    \includegraphics[width=0.49\linewidth]{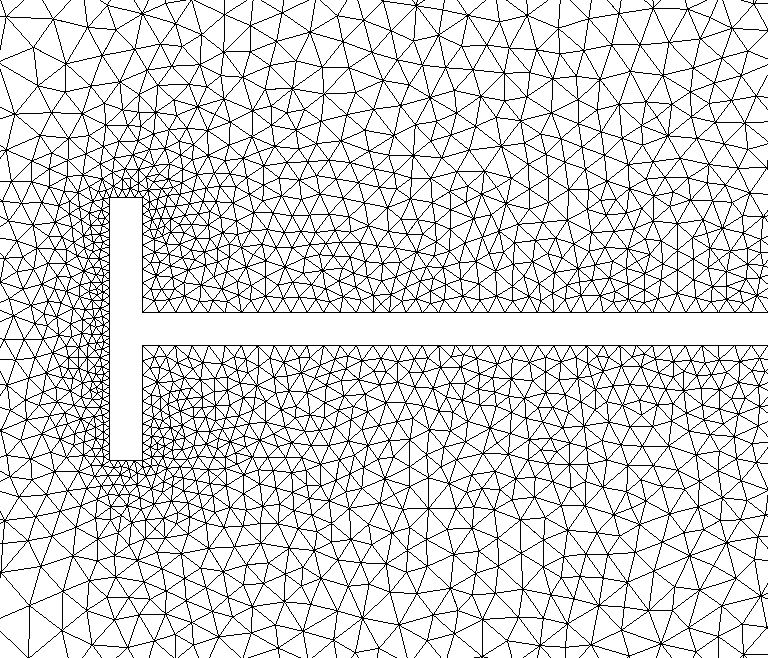} 
    \includegraphics[width=0.49\linewidth]{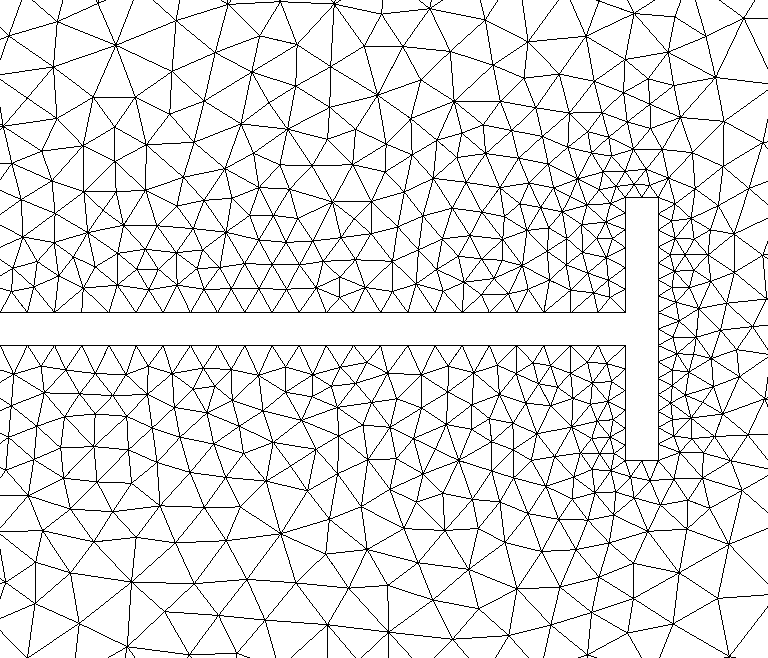} 
    \caption{A zoom of the mesh near the body for the problem
      described in~\cref{ss:bridge}.  Left: fine mesh. Right: coarse
      mesh.}
    \label{fig:fluttering_mesh}
  \end{center}
\end{figure}

We again consider two different meshes: a coarse mesh containing 5901
triangles with 156 body boundary elements and a fine mesh containing
\num{12163} triangles with 272 body boundary elements (see
\cref{fig:fluttering_mesh} for a zoom of the mesh near the body). We
place the inner and outer radii of the annulus consisting of the
sliding and buffer layer at a distance of, respectively, 14 and 15
from the center of the H-shaped bridge. This annulus is initially
partitioned into 100 quadrilaterals when computing on the coarse mesh
and 152 when computing on the fine mesh.

The vertical displacement is handled similarly as in
\cref{ss:transv_gall}. The nodes inside the spatial box
$\Omega^{\text{in}}(t)=[-16,\ 16] \times [-16 + d(t),\ 16 + d(t)]$
move vertically with the body; the vertical movement of the nodes in
$\Omega^{\text{mid}}(t)=([-20,\ 20] \times [-20,\ 20]) \backslash
\Omega^{\text{in}}(t)$ decreases linearly with distance from
$\Omega^{\text{in}}$; and the remaining nodes are fixed. Note that the
rotating and sliding meshes are completely contained in
$\Omega^{\text{in}}(t)$.

\begin{figure}[tbp]
  \begin{center}
    \includegraphics[width=\linewidth]{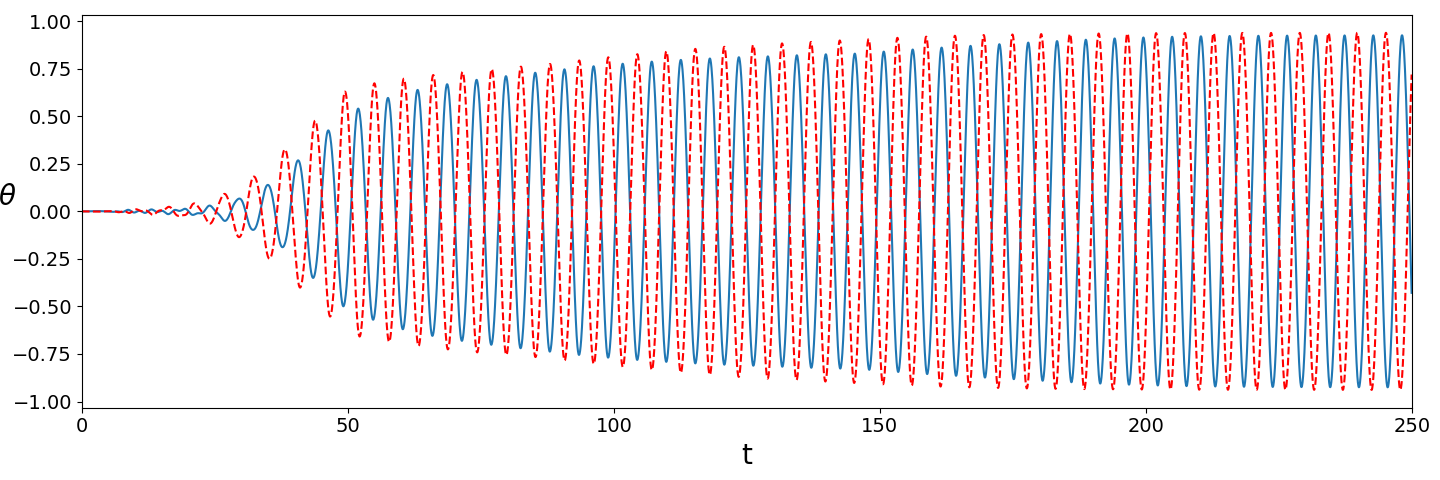} \\
    \includegraphics[width=\linewidth]{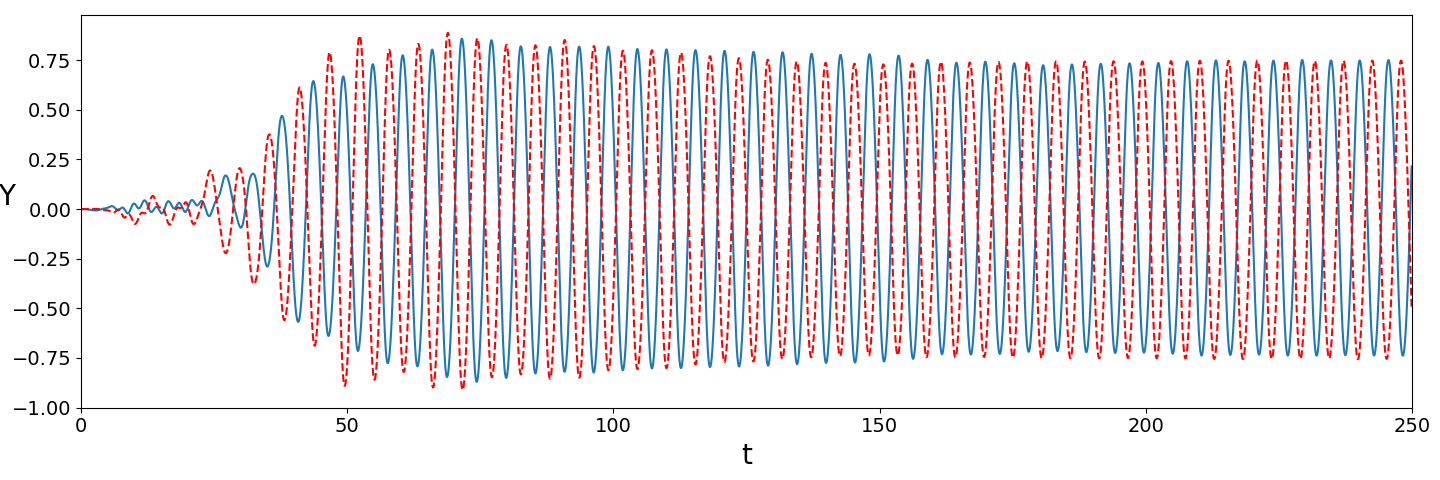} 
    \caption{Evolution of the rotational angle (top) and vertical
      displacement (bottom) computed on the coarse mesh (blue) and the
      fine mesh (red). See \cref{ss:bridge} for more details.}
    \label{fig:fluttering_out}
  \end{center}
\end{figure}

\begin{figure}[tbp]
  \begin{center}
    \includegraphics[width=0.45\linewidth]{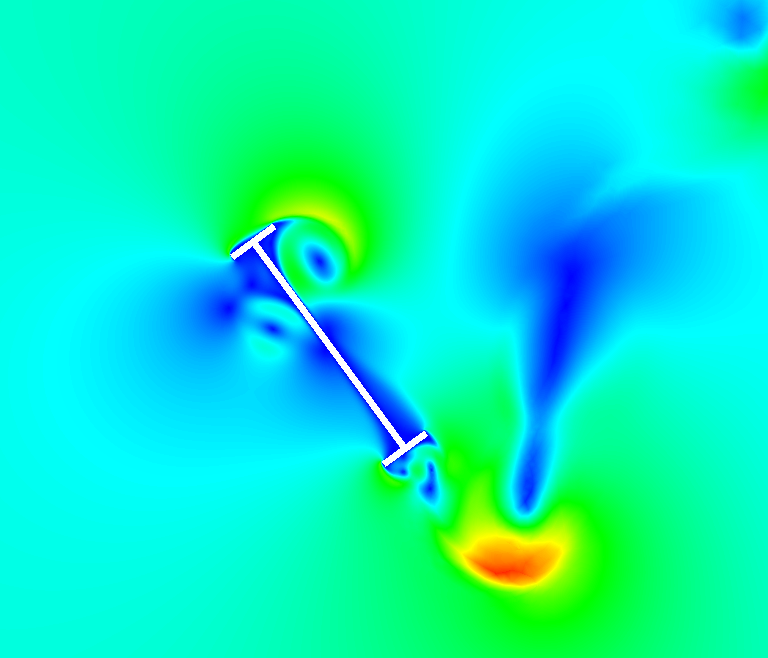} 
    \includegraphics[width=0.45\linewidth]{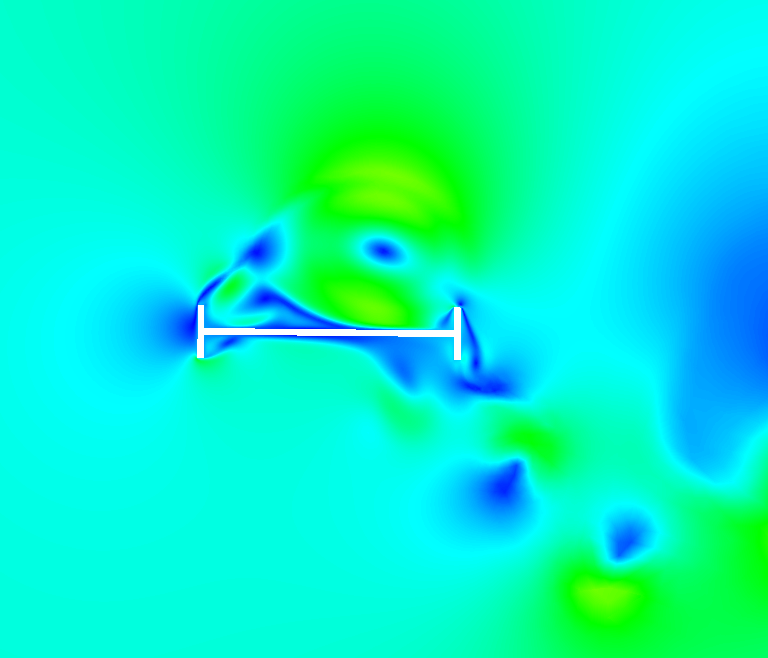} \\
    \includegraphics[width=0.45\linewidth]{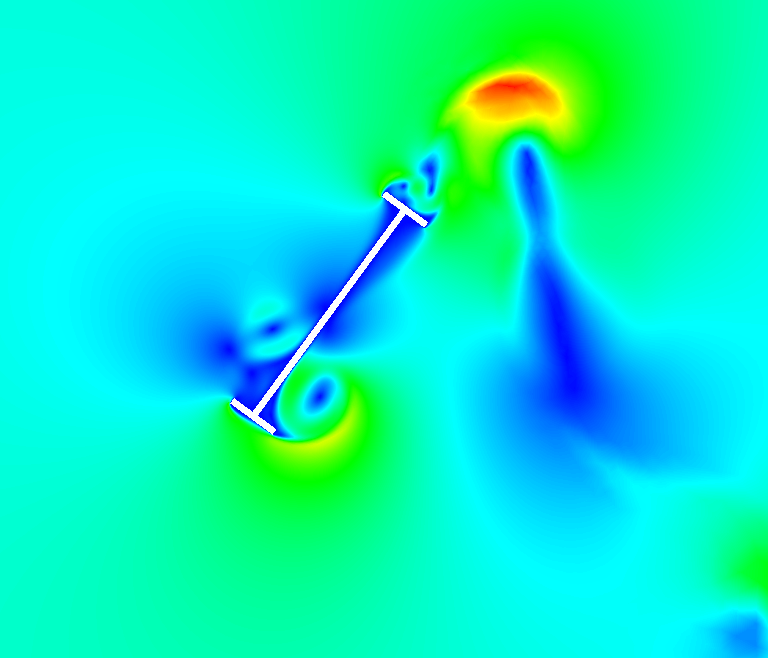} 
    \includegraphics[width=0.45\linewidth]{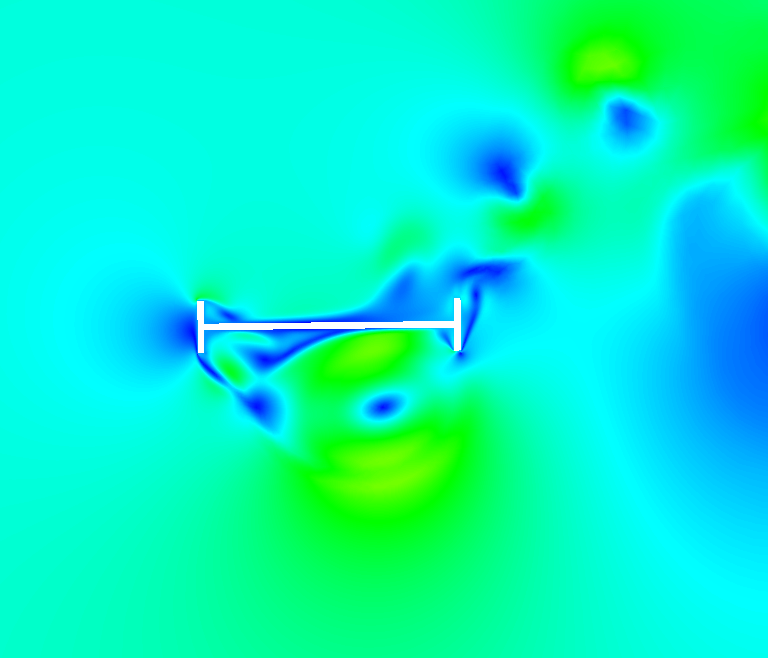} \\
    \includegraphics[width=0.45\linewidth]{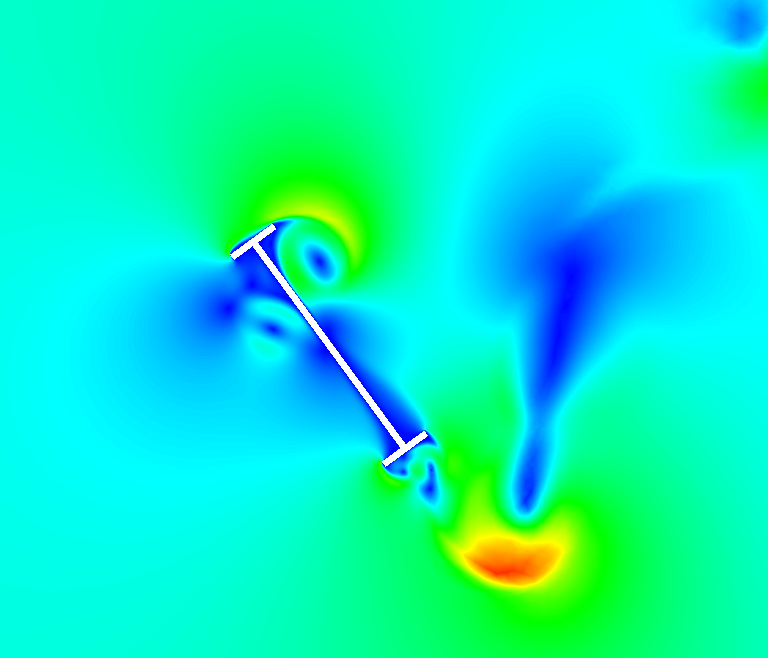} 
    \caption{The velocity magnitude at different, evenly spaced,
      points in time during a full period ($t>190$) for the problem
      described in~\cref{ss:bridge}.}
    \label{fig:bridge_velocity_evolution}
  \end{center}
\end{figure}

We plot, for $t \in [0, 250]$, the evolution of the vertical
displacement and the rotational angle of the rigid body in
\cref{fig:fluttering_out}. We observe that the oscillations develop
into a stable time-periodic pattern at approximately $t = 190$.  For
$t > 190$ we find $\max|y| = 0.74$. For $t > 190$ we also observe that
$\max|\theta| = 0.925$ when computed on the coarse mesh and
$\max|\theta| = 0.935$ when computed on the fine mesh. The frequency
of both the displacement and the rotation is $0.185$. These results
compare well with \cite{Dettmer:2006}, in which
$\max|y|\in [0.75, 0.85]$, $\max|\theta|\in [0.9, 1]$, and the
frequency is $0.186$. To conclude this section, we plot the velocity
magnitude at different, evenly spaced, points in time during a full
period ($t > 190$) in \cref{fig:bridge_velocity_evolution}.

\subsection{Rotating turbine}
\label{ss:turbine}

\begin{figure}[tbp]
  \begin{center}
    \includegraphics[width=\linewidth]{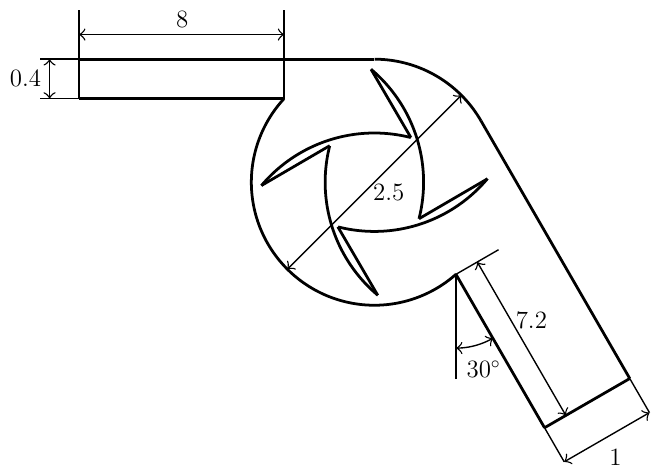} 
    \caption{Set up of the geometry of the rotating turbine example
      of~\cref{ss:turbine}.}
    \label{fig:turbine_geo}
  \end{center}
\end{figure}

In this final test case, we consider a model turbine simulation
inspired by \cite[Section 10.1]{Saksono:2007}. We refer to
\cref{fig:turbine_geo} for the setup of the geometry. For this example
we set $I_\theta = \num{1.5}$, $c_\theta = k_\theta = 0$ in
\cref{eq:linearMotion_theta} (note that we do not have vertical motion
in this example) and choose $\nu = 0.1$ and $\rho = 1$ for the fluid
problem. We impose $\boldsymbol{u} = (1250(1.25-y)(y-0.85), 0)$ on the
left inflow boundary, $\boldsymbol{u} = \boldsymbol{0}$ on the walls,
and a homogeneous Neumann boundary condition on the bottom right
outflow boundary. As time step we choose $\Delta t = 0.001$. For this
test we choose $\delta_{\text{NS}} = 10^{-5}$ in
\cref{eq:stopping_criterion} and $\delta_{\text{rb}} = 10^{-4}$ in
\cref{eq:stopping_crit_rb}.

\begin{figure}[tbp]
  \begin{center}
    \includegraphics[width=\linewidth]{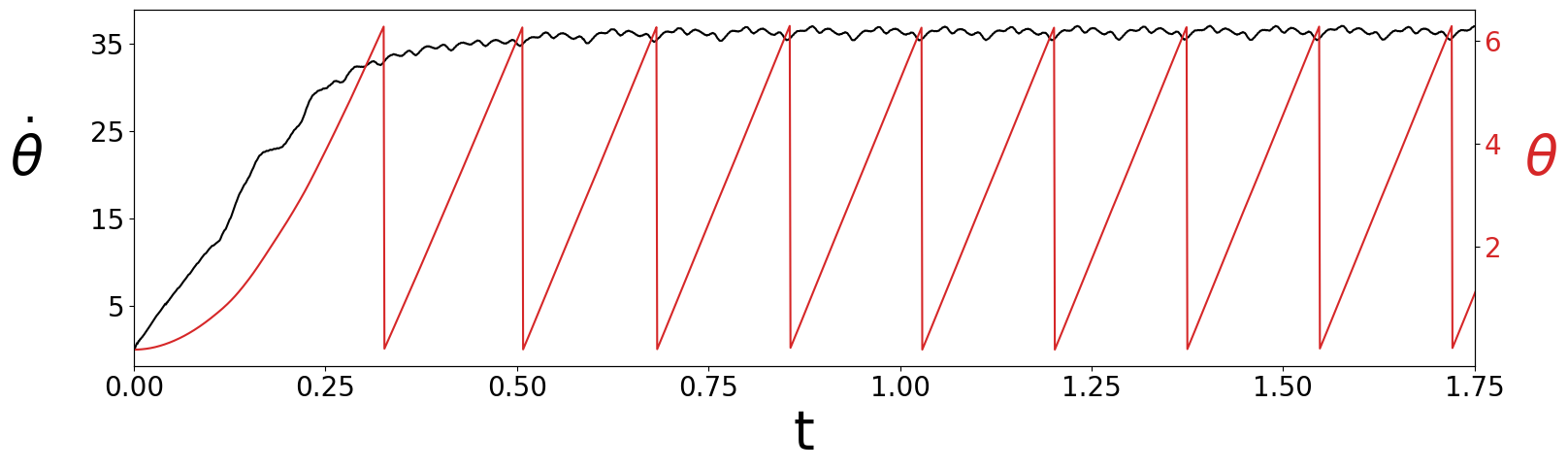}
    \caption{Evolution of the angular velocity $\dot{\theta}$ and
      rotational angle $\theta$ for the rotating turbine problem
      (\cref{ss:turbine}).}
    \label{fig:turbine}
  \end{center}
\end{figure}

\begin{figure}[tbp]
  \begin{center}
    \includegraphics[width=0.45\linewidth]{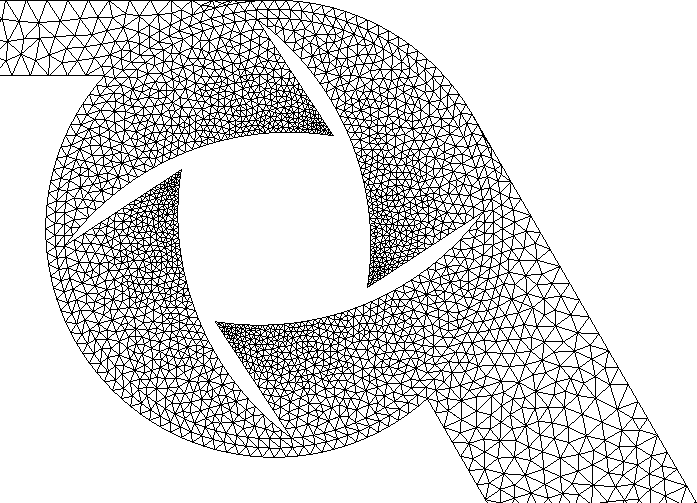} 
    \includegraphics[width=0.45\linewidth]{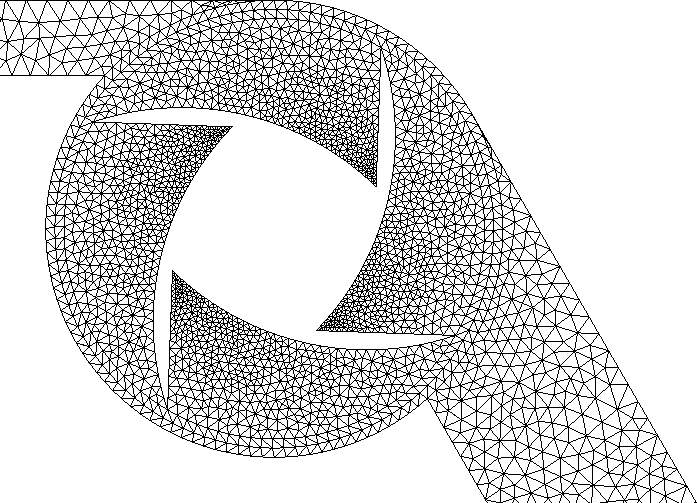} \\
    \includegraphics[width=0.45\linewidth]{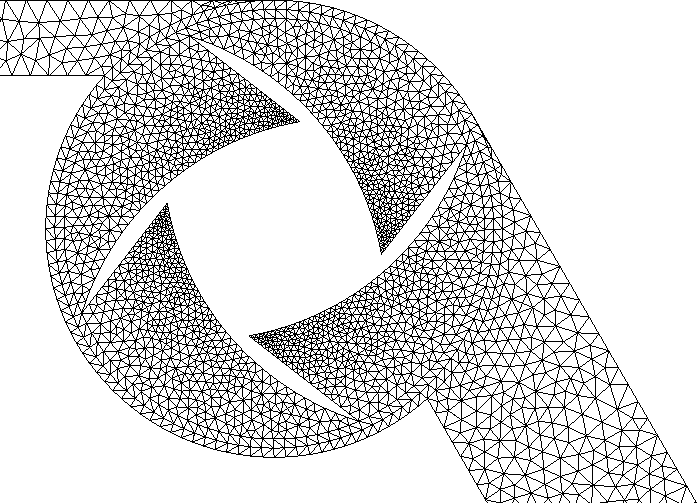}
    \includegraphics[width=0.45\linewidth]{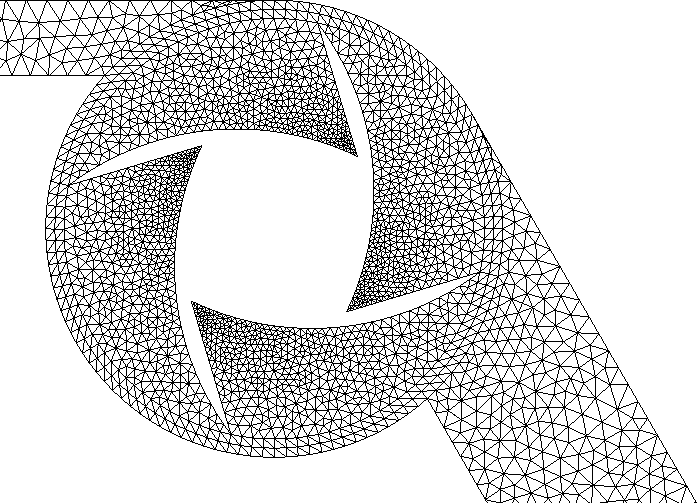} \\
    \includegraphics[width=0.45\linewidth]{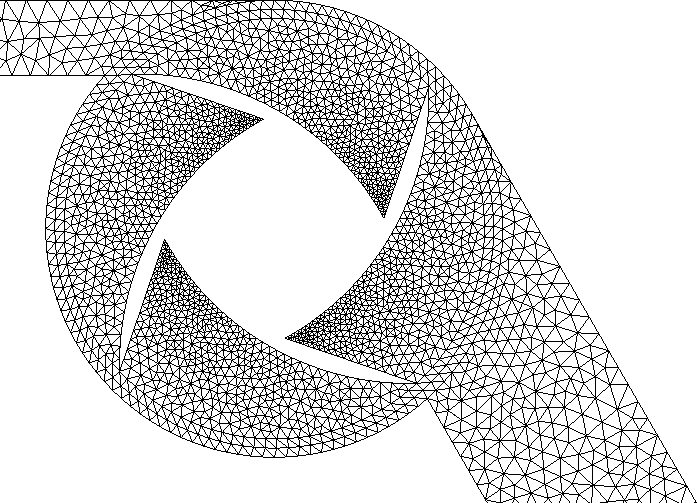} 
    \includegraphics[width=0.45\linewidth]{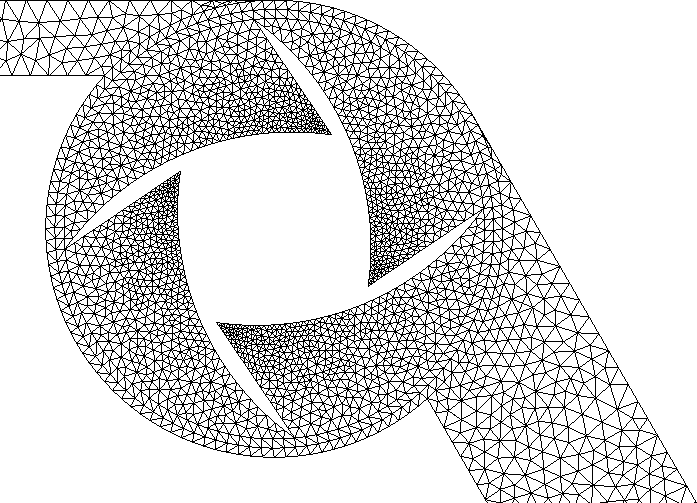} 
    \caption{The meshes at different, evenly spaced, points in time
      during half a rotation of the turbine for the problem described
      in~\cref{ss:turbine}.}
    \label{fig:turbine_mesh_evolution}
  \end{center}
\end{figure}

\begin{figure}[tbp]
  \begin{center}
    \includegraphics[width=0.45\linewidth]{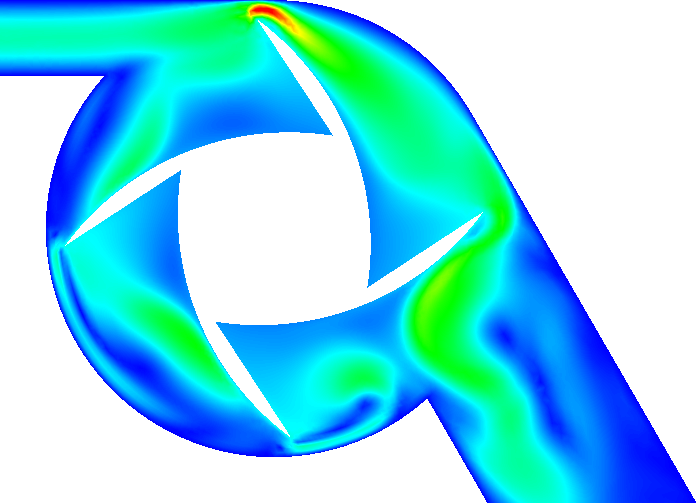} 
    \includegraphics[width=0.45\linewidth]{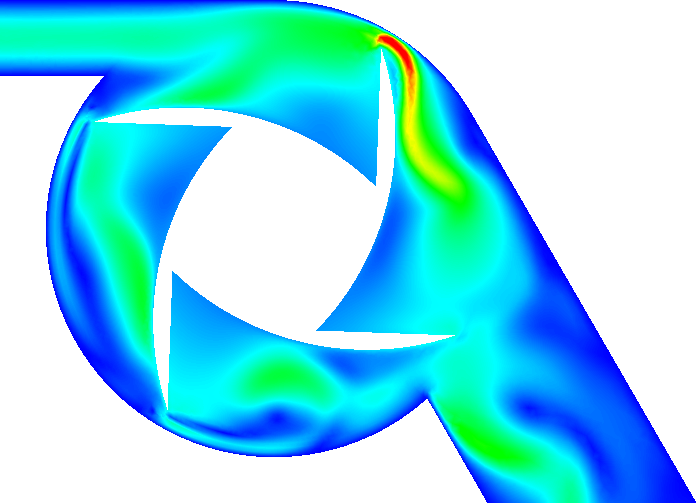} \\
    \includegraphics[width=0.45\linewidth]{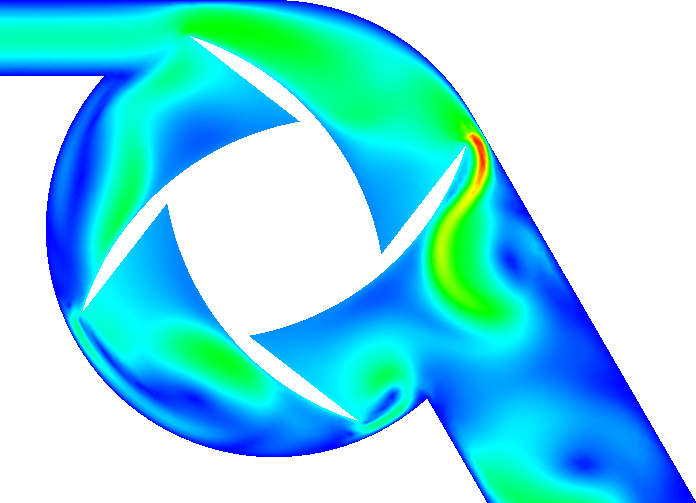}
    \includegraphics[width=0.45\linewidth]{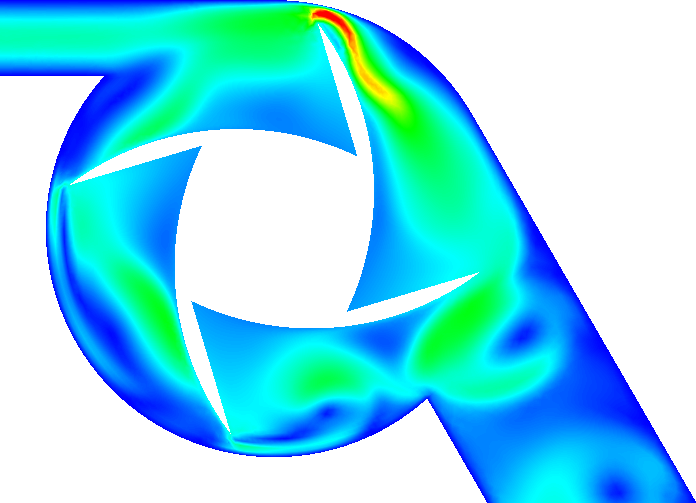} \\
    \includegraphics[width=0.45\linewidth]{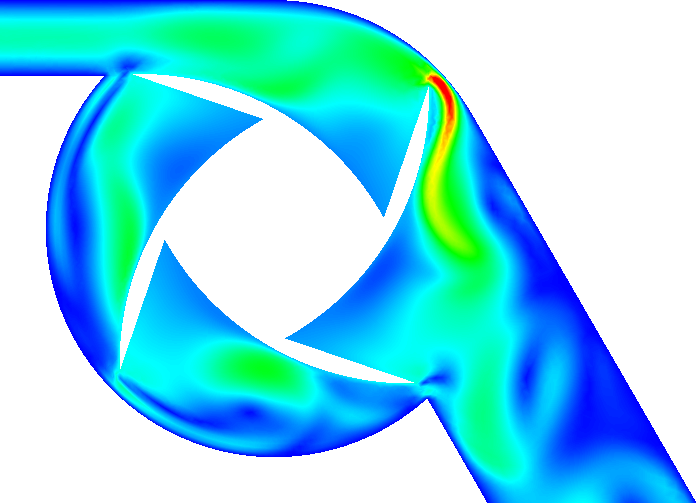} 
    \includegraphics[width=0.45\linewidth]{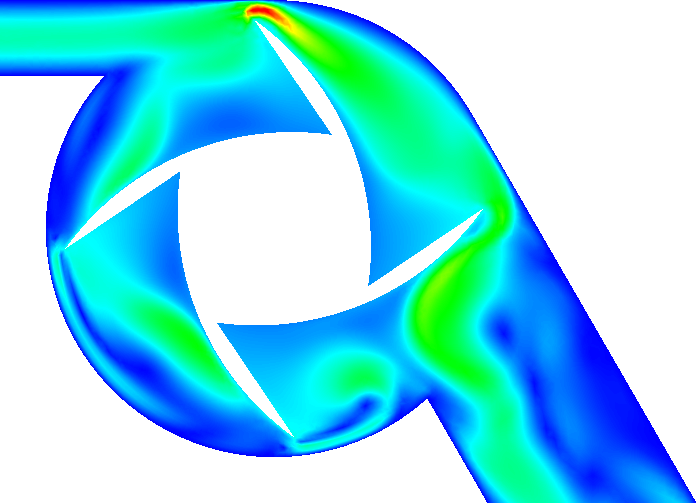} 
    \caption{The velocity magnitude computed on the meshes of
      \cref{fig:turbine_mesh_evolution}.}
    \label{fig:turbine_velocity_evolution}
  \end{center}
\end{figure}

In \cref{fig:turbine} we plot the evolution of the rotational angle
and the angular velocity of the turbine for $t \in [0, 1.75]$. We
observe that the rotational velocity increases until $t = 0.5$ and
then slightly oscillates around $\dot{\theta}=35$. We furthermore
observe from \cref{fig:turbine} that the turbine completes multiple
full rotations throughout the simulation. We show the mesh at
different, evenly spaced, points in time during half a rotation
($t > 0.5$) in \cref{fig:turbine_mesh_evolution}. Finally, we plot the
velocity magnitude on these meshes in
\cref{fig:turbine_velocity_evolution}.

\section{Conclusions}
\label{s:conclusion}

In this paper we presented an exactly mass conserving space-time
embedded-hybridizable discontinuous Galerkin finite element
discretization for fluid-rigid body interaction problems. Furthermore,
we introduced a conforming sliding mesh technique, using edge
swapping, for simulations where the rotational motion of the rigid
body is large. By noting that our edge swapping algorithm results in
one of only four different mesh connectivity sets, and that each of
these sets can be determined before the start of a simulation,
additional computational cost due to edge swapping is kept low.

\subsubsection*{Acknowledgements}

TH would like to thank Oakland University for the URC Faculty Research
Fellowship Award. SR gratefully acknowledges support from the Natural
Sciences and Engineering Research Council of Canada through the
Discovery Grant program (RGPIN-05606-2015).

\bibliographystyle{elsarticle-num-names}
\bibliography{references}

\begin{thebibliography}{45}
\expandafter\ifx\csname natexlab\endcsname\relax\def\natexlab#1{#1}\fi
\providecommand{\url}[1]{\texttt{#1}}
\providecommand{\href}[2]{#2}
\providecommand{\path}[1]{#1}
\providecommand{\DOIprefix}{doi:}
\providecommand{\ArXivprefix}{arXiv:}
\providecommand{\URLprefix}{URL: }
\providecommand{\Pubmedprefix}{pmid:}
\providecommand{\doi}[1]{\href{http://dx.doi.org/#1}{\path{#1}}}
\providecommand{\Pubmed}[1]{\href{pmid:#1}{\path{#1}}}
\providecommand{\bibinfo}[2]{#2}
\ifx\xfnm\relax \def\xfnm[#1]{\unskip,\space#1}\fi
\bibitem[{Robertson et~al.(2003)Robertson, Li, Sherwin, and
  Bearman}]{Robertson:2003}
\bibinfo{author}{I.~Robertson}, \bibinfo{author}{L.~Li}, \bibinfo{author}{S.~J.
  Sherwin}, \bibinfo{author}{P.~W. Bearman},
\newblock \bibinfo{title}{A numerical study of rotational and transverse
  galloping rectangular bodies},
\newblock \bibinfo{journal}{J. Fluid. Struct} \bibinfo{volume}{17}
  (\bibinfo{year}{2003}) \bibinfo{pages}{681--699}.
  \DOIprefix\doi{10.1016/S0889-9746(03)00008-2}.
\bibitem[{Dettmer and Peri{\'c}(2006)}]{Dettmer:2006}
\bibinfo{author}{W.~Dettmer}, \bibinfo{author}{D.~Peri{\'c}},
\newblock \bibinfo{title}{A computational framework for fluid-rigid body
  interaction: finite element formulation and applications},
\newblock \bibinfo{journal}{Comput. Methods. Appl. Mech. Engrg.}
  \bibinfo{volume}{195} (\bibinfo{year}{2006}) \bibinfo{pages}{1633--1666}.
  \DOIprefix\doi{10.1016/j.cma.2005.05.033}.
\bibitem[{Saksono et~al.(2007)Saksono, Dettmer, and Peri{\'c}}]{Saksono:2007}
\bibinfo{author}{P.~H. Saksono}, \bibinfo{author}{W.~G. Dettmer},
  \bibinfo{author}{D.~Peri{\'c}},
\newblock \bibinfo{title}{An adaptive remeshing strategy for flows with moving
  boundaries and fluid-structure interaction},
\newblock \bibinfo{journal}{Int. J. Numer. Meth. Eng.} \bibinfo{volume}{71}
  (\bibinfo{year}{2007}) \bibinfo{pages}{1009--1050}.
  \DOIprefix\doi{10.1002/nme.1971}.
\bibitem[{{\v{C}}ani{\'c}(2021)}]{Canic:2021}
\bibinfo{author}{S.~{\v{C}}ani{\'c}},
\newblock \bibinfo{title}{Moving boundary problems},
\newblock \bibinfo{journal}{B. Am. Math. Soc.} \bibinfo{volume}{58}
  (\bibinfo{year}{2021}) \bibinfo{pages}{79--106}.
  \DOIprefix\doi{10.1090/bull/1703}.
\bibitem[{Buka{\v{c}} et~al.(2013)Buka{\v{c}}, {\v{C}}ani{\'c}, Glowinski,
  Tamba{\v{c}}a, and Quaini}]{Bukac:2013}
\bibinfo{author}{M.~Buka{\v{c}}}, \bibinfo{author}{S.~{\v{C}}ani{\'c}},
  \bibinfo{author}{R.~Glowinski}, \bibinfo{author}{J.~Tamba{\v{c}}a},
  \bibinfo{author}{A.~Quaini},
\newblock \bibinfo{title}{Fluid-structure interaction in blood flow capturing
  non-zero longitudinal structure displacement},
\newblock \bibinfo{journal}{J. Comput. Phys.} \bibinfo{volume}{235}
  (\bibinfo{year}{2013}) \bibinfo{pages}{515--541}.
  \DOIprefix\doi{10.1016/j.jcp.2012.08.033}.
\bibitem[{Buka{\v{c}} et~al.(2019)Buka{\v{c}}, {\v{C}}ani{\'c}, Tamba{\v{c}}a,
  and Wang}]{Bukac:2019}
\bibinfo{author}{M.~Buka{\v{c}}}, \bibinfo{author}{S.~{\v{C}}ani{\'c}},
  \bibinfo{author}{J.~Tamba{\v{c}}a}, \bibinfo{author}{Y.~Wang},
\newblock \bibinfo{title}{Fluid-structure interaction between pulsatile blood
  flow and a curved stented coronary artery on a beating heart: a four stent
  computational study},
\newblock \bibinfo{journal}{Comput. Methods. Appl. Mech. Engrg.}
  \bibinfo{volume}{350} (\bibinfo{year}{2019}) \bibinfo{pages}{679--700}.
  \DOIprefix\doi{10.1016/j.cma.2019.03.034}.
\bibitem[{Neunteufel and Sch\"oberl(2021)}]{Neunteufel:2021}
\bibinfo{author}{M.~Neunteufel}, \bibinfo{author}{J.~Sch\"oberl},
\newblock \bibinfo{title}{Fluid-structure interaction with {H}(div)-conforming
  finite elements},
\newblock \bibinfo{journal}{Comput. Struct.} \bibinfo{volume}{243}
  (\bibinfo{year}{2021}) \bibinfo{pages}{106402}.
  \DOIprefix\doi{10.1016/j.compstruc.2020.106402}.
\bibitem[{Fu and Kuang(2020)}]{Fu:2020}
\bibinfo{author}{G.~Fu}, \bibinfo{author}{W.~Kuang},
\newblock \bibinfo{title}{A monolithic divergence-conforming {HDG} scheme for a
  linear fluid-structure interaction model},
\newblock \bibinfo{journal}{arXiv preprint arXiv:2012.00128}
  (\bibinfo{year}{2020}).
\bibitem[{Labeur and Wells(2009)}]{Labeur:2009}
\bibinfo{author}{R.~J. Labeur}, \bibinfo{author}{G.~N. Wells},
\newblock \bibinfo{title}{Interface stabilised finite element method for moving
  domains and free surface flows},
\newblock \bibinfo{journal}{Comput. Methods Appl. Mech. Engrg.}
  \bibinfo{volume}{198} (\bibinfo{year}{2009}) \bibinfo{pages}{615--630}.
  \DOIprefix\doi{10.1016/j.cma.2008.09.014}.
\bibitem[{Hauke and Hughes(1994)}]{Hauke:1994}
\bibinfo{author}{G.~Hauke}, \bibinfo{author}{T.~J.~R. Hughes},
\newblock \bibinfo{title}{A unified approach to compressible and incompressible
  flows},
\newblock \bibinfo{journal}{Comput. Methods. Appl. Mech. Engrg.}
  \bibinfo{volume}{113} (\bibinfo{year}{1994}) \bibinfo{pages}{389--395}.
  \DOIprefix\doi{10.1016/0045-7825(94)90055-8}.
\bibitem[{Johnson and Tezduyar(1994)}]{Johnson:1994}
\bibinfo{author}{A.~A. Johnson}, \bibinfo{author}{T.~E. Tezduyar},
\newblock \bibinfo{title}{Mesh update strategies in parallel finite element
  computations of flow problems with moving boundaries and interfaces},
\newblock \bibinfo{journal}{Comput. Methods. Appl. Mech. Engrg.}
  \bibinfo{volume}{119} (\bibinfo{year}{1994}) \bibinfo{pages}{73--94}.
  \DOIprefix\doi{10.1016/0045-7825(94)00077-8}.
\bibitem[{N'dri et~al.(2001)N'dri, Garon, and Fortin}]{Ndri:2001}
\bibinfo{author}{D.~N'dri}, \bibinfo{author}{A.~Garon},
  \bibinfo{author}{A.~Fortin},
\newblock \bibinfo{title}{A new stable space-time formulation for
  two-dimensional and three-dimensional incompressible viscous flow},
\newblock \bibinfo{journal}{Int. J. Numer. Meth. Fluids} \bibinfo{volume}{37}
  (\bibinfo{year}{2001}) \bibinfo{pages}{865--884}.
  \DOIprefix\doi{10.1002/fld.174}.
\bibitem[{N'dri et~al.(2002)N'dri, Garon, and Fortin}]{Ndri:2002}
\bibinfo{author}{D.~N'dri}, \bibinfo{author}{A.~Garon},
  \bibinfo{author}{A.~Fortin},
\newblock \bibinfo{title}{Incompressible {N}avier--{S}tokes computations with
  stable and stabilized space-time formulations: a comparative study},
\newblock \bibinfo{journal}{Commun. Numer. Meth. Engng.} \bibinfo{volume}{18}
  (\bibinfo{year}{2002}) \bibinfo{pages}{495--512}.
  \DOIprefix\doi{10.1002/cnm.507}.
\bibitem[{Rhebergen et~al.(2013)Rhebergen, Cockburn, and van~der
  Vegt}]{Rhebergen:2013b}
\bibinfo{author}{S.~Rhebergen}, \bibinfo{author}{B.~Cockburn},
  \bibinfo{author}{J.~van~der Vegt},
\newblock \bibinfo{title}{A space-time discontinuous {G}alerkin method for the
  incompressible {N}avier--{S}tokes equations},
\newblock \bibinfo{journal}{J. Comput. Phys.} \bibinfo{volume}{233}
  (\bibinfo{year}{2013}) \bibinfo{pages}{339--358}.
  \DOIprefix\doi{10.1016/j.jcp.2012.08.052}.
\bibitem[{Tavelli and Dumbser(2015)}]{Tavelli:2015}
\bibinfo{author}{M.~Tavelli}, \bibinfo{author}{M.~Dumbser},
\newblock \bibinfo{title}{A staggered space-time discontinuous {G}alerkin
  method for the incompressible {N}avier--{S}tokes equations on two-dimensional
  triangular meshes},
\newblock \bibinfo{journal}{Comput. Fluids} \bibinfo{volume}{119}
  (\bibinfo{year}{2015}) \bibinfo{pages}{235--249}.
  \DOIprefix\doi{10.1016/j.compfluid.2015.07.003}.
\bibitem[{Tavelli and Dumbser(2016)}]{Tavelli:2016}
\bibinfo{author}{M.~Tavelli}, \bibinfo{author}{M.~Dumbser},
\newblock \bibinfo{title}{A staggered space-time discontinuous {G}alerkin
  method for the three-dimensional incompressible {N}avier--{S}tokes equations
  on unstructured tetrahedral meshes},
\newblock \bibinfo{journal}{J. Comput. Phys.} \bibinfo{volume}{319}
  (\bibinfo{year}{2016}) \bibinfo{pages}{294--323}.
  \DOIprefix\doi{10.1016/j.jcp.2016.05.009}.
\bibitem[{van~der Vegt and Sudirham(2008)}]{Vegt:2008}
\bibinfo{author}{J.~van~der Vegt}, \bibinfo{author}{J.~Sudirham},
\newblock \bibinfo{title}{A space-time discontinuous {G}alerkin method for the
  time-dependent {O}seen equations},
\newblock \bibinfo{journal}{Appl. Numer. Math} \bibinfo{volume}{58}
  (\bibinfo{year}{2008}) \bibinfo{pages}{1892--1917}.
  \DOIprefix\doi{10.1016/j.apnum.2007.11.010}.
\bibitem[{Tezduyar et~al.(1992{\natexlab{a}})Tezduyar, Behr, and
  Liou}]{Tezduyar:1992a}
\bibinfo{author}{T.~E. Tezduyar}, \bibinfo{author}{M.~Behr},
  \bibinfo{author}{J.~Liou},
\newblock \bibinfo{title}{A new strategy for finite element computations
  involving moving boundaries and interfaces--the {DSD/ST} procedure: {I}.
  {T}he concept and the preliminary numerical tests},
\newblock \bibinfo{journal}{Comput. Methods. Appl. Mech. Engrg.}
  \bibinfo{volume}{94} (\bibinfo{year}{1992}{\natexlab{a}})
  \bibinfo{pages}{339--351}. \DOIprefix\doi{10.1016/0045-7825(92)90059-S}.
\bibitem[{Tezduyar et~al.(1992{\natexlab{b}})Tezduyar, Behr, Mittal, and
  Liou}]{Tezduyar:1992b}
\bibinfo{author}{T.~E. Tezduyar}, \bibinfo{author}{M.~Behr},
  \bibinfo{author}{S.~Mittal}, \bibinfo{author}{J.~Liou},
\newblock \bibinfo{title}{A new strategy for finite element computations
  involving moving boundaries and interfaces--{T}he
  deforming-spatial-domain/space-time procedure: {II}. {C}omputation of
  free-surface flows, two-liquid flows, and flows with drifting cylinders},
\newblock \bibinfo{journal}{Comput. Methods. Appl. Mech. Engrg.}
  \bibinfo{volume}{94} (\bibinfo{year}{1992}{\natexlab{b}})
  \bibinfo{pages}{353--371}. \DOIprefix\doi{10.1016/0045-7825(92)90060-W}.
\bibitem[{Tezduyar et~al.(2006)Tezduyar, Sathe, Keedy, and
  Stein}]{Tezduyar:2006}
\bibinfo{author}{T.~E. Tezduyar}, \bibinfo{author}{S.~Sathe},
  \bibinfo{author}{R.~Keedy}, \bibinfo{author}{K.~Stein},
\newblock \bibinfo{title}{Space-time finite element techniques for computation
  of fluid-structure interactions},
\newblock \bibinfo{journal}{Comput. Methods Appl. Mech. Engrg.}
  \bibinfo{volume}{195} (\bibinfo{year}{2006}) \bibinfo{pages}{2002--2027}.
  \DOIprefix\doi{10.1016/j.cma.2004.09.014}.
\bibitem[{H{\"u}bner et~al.(2004)H{\"u}bner, Walhorn, and
  Dinkler}]{Hubner:2004}
\bibinfo{author}{B.~H{\"u}bner}, \bibinfo{author}{E.~Walhorn},
  \bibinfo{author}{D.~Dinkler},
\newblock \bibinfo{title}{A monolithic approach to fluid-structure interaction
  using space-time finite elements},
\newblock \bibinfo{journal}{Comput. Methods Appl. Mech. Engrg.}
  \bibinfo{volume}{193} (\bibinfo{year}{2004}) \bibinfo{pages}{2087--2104}.
  \DOIprefix\doi{10.1016/j.cma.2004.01.024}.
\bibitem[{Horvath and Rhebergen(2020)}]{Horvath:2020}
\bibinfo{author}{T.~L. Horvath}, \bibinfo{author}{S.~Rhebergen},
\newblock \bibinfo{title}{An exactly mass conserving space-time
  embedded-hybridized discontinuous {G}alerkin method for the
  {N}avier--{S}tokes equations on moving domains},
\newblock \bibinfo{journal}{J. Comput. Phys.} \bibinfo{volume}{417}
  (\bibinfo{year}{2020}). \DOIprefix\doi{10.1016/j.jcp.2020.109577}.
\bibitem[{Rhebergen and Wells(2017)}]{Rhebergen:2017}
\bibinfo{author}{S.~Rhebergen}, \bibinfo{author}{G.~N. Wells},
\newblock \bibinfo{title}{Analysis of a hybridized/interface stabilized finite
  element method for the {S}tokes equations},
\newblock \bibinfo{journal}{SIAM J. Numer. Anal.} \bibinfo{volume}{55}
  (\bibinfo{year}{2017}) \bibinfo{pages}{1982--2003}.
  \DOIprefix\doi{10.1137/16M1083839}.
\bibitem[{Rhebergen and Wells(2018)}]{Rhebergen:2018}
\bibinfo{author}{S.~Rhebergen}, \bibinfo{author}{G.~N. Wells},
\newblock \bibinfo{title}{A hybridizable discontinuous {G}alerkin method for
  the {N}avier--{S}tokes equations with pointwise divergence-free velocity
  field},
\newblock \bibinfo{journal}{J. Sci. Comput.} \bibinfo{volume}{76}
  (\bibinfo{year}{2018}) \bibinfo{pages}{1484--1501}.
  \DOIprefix\doi{10.1007/s10915-018-0671-4}.
\bibitem[{Rhebergen and Wells(2020)}]{Rhebergen:2020}
\bibinfo{author}{S.~Rhebergen}, \bibinfo{author}{G.~N. Wells},
\newblock \bibinfo{title}{An embedded-hybridized discontinuous {G}alerkin
  finite element method for the {S}tokes equations},
\newblock \bibinfo{journal}{Comput. Methods Appl. Mech. Engrg.}
  \bibinfo{volume}{367} (\bibinfo{year}{2020}).
  \DOIprefix\doi{10.1016/j.cam.2019.112476}.
\bibitem[{Horvath and Rhebergen(2019)}]{Horvath:2019}
\bibinfo{author}{T.~L. Horvath}, \bibinfo{author}{S.~Rhebergen},
\newblock \bibinfo{title}{A locally conservative and energy-stable finite
  element method for the {N}avier--{S}tokes problem on time-dependent domains},
\newblock \bibinfo{journal}{Int. J. Numer. Meth. Fluids} \bibinfo{volume}{89}
  (\bibinfo{year}{2019}) \bibinfo{pages}{519--532}.
  \DOIprefix\doi{10.1002/fld.4707}.
\bibitem[{Kirk et~al.(2021)Kirk, Horvath, and Rhebergen}]{Kirk:2021}
\bibinfo{author}{K.~L.~A. Kirk}, \bibinfo{author}{T.~L. Horvath},
  \bibinfo{author}{S.~Rhebergen},
\newblock \bibinfo{title}{Analysis of an exactly mass conserving space-time
  hybridized discontinuous {G}alerkin method for the time-dependent
  {N}avier--{S}tokes equations},
\newblock \bibinfo{journal}{arXiv preprint arXiv:2103.13492}
  (\bibinfo{year}{2021}).
\bibitem[{Cockburn et~al.(2009)Cockburn, Gopalakrishnan, and
  Lazarov}]{Cockburn:2009a}
\bibinfo{author}{B.~Cockburn}, \bibinfo{author}{J.~Gopalakrishnan},
  \bibinfo{author}{R.~Lazarov},
\newblock \bibinfo{title}{Unified hybridization of discontinuous {G}alerkin,
  mixed, and continuous {G}alerkin methods for second order elliptic problems},
\newblock \bibinfo{journal}{SIAM J. Numer. Anal.} \bibinfo{volume}{47}
  (\bibinfo{year}{2009}) \bibinfo{pages}{1319--1365}.
  \DOIprefix\doi{10.1137/070706616}.
\bibitem[{Ferrer and Willden(2012)}]{Ferrer:2012}
\bibinfo{author}{E.~Ferrer}, \bibinfo{author}{R.~H.~J. Willden},
\newblock \bibinfo{title}{A high order discontinuous {G}alerkin -- {F}ourier
  incompressible 3{D} {N}avier--{S}tokes solver with rotating sliding meshes},
\newblock \bibinfo{journal}{J. Comput. Phys.} \bibinfo{volume}{231}
  (\bibinfo{year}{2012}) \bibinfo{pages}{7037--7056}.
  \DOIprefix\doi{10.1016/j.jcp.2012.04.039}.
\bibitem[{Cottrell et~al.(2009)Cottrell, Hughes, and Bazilevs}]{Cottrell:2009}
\bibinfo{author}{J.~A. Cottrell}, \bibinfo{author}{T.~J.~R. Hughes},
  \bibinfo{author}{Y.~Bazilevs}, \bibinfo{title}{Isogeometric Analysis: Toward
  Integration of CAD and FEA}, \bibinfo{publisher}{Wiley},
  \bibinfo{address}{New York}, \bibinfo{year}{2009}.
\bibitem[{Anagnostou(1990)}]{Anagnostou:1990}
\bibinfo{author}{G.~Anagnostou}, \bibinfo{title}{Nonconforming sliding spectral
  element methods for the unsteady incompressible Navier--Stokes equations},
  Ph.D. thesis, Massachusetts Institute of Technology, \bibinfo{year}{1990}.
\bibitem[{Wang and Persson(2015)}]{Wang:2015}
\bibinfo{author}{L.~Wang}, \bibinfo{author}{P.-O. Persson},
\newblock \bibinfo{title}{A high-order discontinuous {G}alerkin method with
  unstructured space-time meshes for two-dimensional compressible flows on
  domains with large deformations},
\newblock \bibinfo{journal}{Comput. Fluids} \bibinfo{volume}{118}
  (\bibinfo{year}{2015}) \bibinfo{pages}{53--68}.
  \DOIprefix\doi{10.1016/j.compfluid.2015.05.026}.
\bibitem[{Calderer and Masud(2010)}]{Calderer:2010}
\bibinfo{author}{R.~Calderer}, \bibinfo{author}{A.~Masud},
\newblock \bibinfo{title}{A multiscale stabilized {ALE} formulation for
  incompressible flows with moving boundaries},
\newblock \bibinfo{journal}{Comput. Mech.} \bibinfo{volume}{46}
  (\bibinfo{year}{2010}) \bibinfo{pages}{185--197}.
  \DOIprefix\doi{10.1007/s00466-010-0487-z}.
\bibitem[{Farhat et~al.(2006)Farhat, der Zee, and Geuzaine}]{Farhat:2006}
\bibinfo{author}{C.~Farhat}, \bibinfo{author}{K.~G.~V. der Zee},
  \bibinfo{author}{P.~Geuzaine},
\newblock \bibinfo{title}{Provably second-order time-accurate loosely-coupled
  solution algorithms for transient nonlinear computational aeroelasticity},
\newblock \bibinfo{journal}{Comput. Methods Appl. Mech. Engrg.}
  \bibinfo{volume}{195} (\bibinfo{year}{2006}) \bibinfo{pages}{1973--2001}.
  \DOIprefix\doi{10.1016/j.cma.2004.11.031}.
\bibitem[{Kadapa et~al.(2017)Kadapa, Dettmer, and Peri{\'c}}]{Kadapa:2017}
\bibinfo{author}{C.~Kadapa}, \bibinfo{author}{W.~G. Dettmer},
  \bibinfo{author}{D.~Peri{\'c}},
\newblock \bibinfo{title}{A stabilised immersed boundary method on hierarchical
  b-spline grids for fluid-rigid body interaction with solid-solid contact},
\newblock \bibinfo{journal}{Comput. Methods Appl. Mech. Engrg.}
  \bibinfo{volume}{318} (\bibinfo{year}{2017}) \bibinfo{pages}{242--269}.
  \DOIprefix\doi{10.1016/j.cma.2017.01.024}.
\bibitem[{Rugonyi and Bathe(2001)}]{Rugonyi:2001}
\bibinfo{author}{S.~Rugonyi}, \bibinfo{author}{K.~J. Bathe},
\newblock \bibinfo{title}{On finite element analysis of fluid flows fully
  coupled with structural interactions},
\newblock \bibinfo{journal}{CMES Comput. Model. Eng. Sci.} \bibinfo{volume}{2}
  (\bibinfo{year}{2001}) \bibinfo{pages}{195--212}.
  \DOIprefix\doi{10.3970/cmes.2001.002.195}.
\bibitem[{Tezduyar and Sathe(2007)}]{Tezduyar:2007}
\bibinfo{author}{T.~E. Tezduyar}, \bibinfo{author}{S.~Sathe},
\newblock \bibinfo{title}{Modelling of fluid-structure interactions with the
  space-time finite elements: {S}olution techniques},
\newblock \bibinfo{journal}{Int. J. Numer. Meth. Fl.} \bibinfo{volume}{54}
  (\bibinfo{year}{2007}) \bibinfo{pages}{855--900}.
  \DOIprefix\doi{10.1002/fld.1430}.
\bibitem[{{van der Vegt} and {van der Ven}(2002)}]{Vegt:2002}
\bibinfo{author}{J.~J.~W. {van der Vegt}}, \bibinfo{author}{H.~{van der Ven}},
\newblock \bibinfo{title}{Space-time discontinuous {G}alerkin finite element
  method with dynamic grid motion for inviscid compressible flow},
\newblock \bibinfo{journal}{J. Comput. Phys.} \bibinfo{volume}{182}
  (\bibinfo{year}{2002}) \bibinfo{pages}{546--585}.
  \DOIprefix\doi{10.1006/jcph.2002.7185}.
\bibitem[{Dobrev et~al.(2018)Dobrev, Kolev et~al.}]{mfem-library}
\bibinfo{author}{V.~A. Dobrev}, \bibinfo{author}{T.~V. Kolev}, et~al.,
  \bibinfo{title}{{MFEM}: Modular finite element methods},
  \bibinfo{howpublished}{\url{http://mfem.org}}, \bibinfo{year}{2018}.
\bibitem[{Amestoy et~al.(2001)Amestoy, Duff, L'Excellent, and Koster}]{MUMPS:1}
\bibinfo{author}{P.~Amestoy}, \bibinfo{author}{I.~Duff}, \bibinfo{author}{J.-Y.
  L'Excellent}, \bibinfo{author}{J.~Koster},
\newblock \bibinfo{title}{A fully asynchronous multifrontal solver using
  distributed dynamic scheduling},
\newblock \bibinfo{journal}{SIAM J. Matrix Anal. \& Appl.} \bibinfo{volume}{23}
  (\bibinfo{year}{2001}) \bibinfo{pages}{15--41}.
  \DOIprefix\doi{10.1137/S0895479899358194}.
\bibitem[{Amestoy et~al.(2006)Amestoy, Guermouche, L'Excellent, and
  Pralet}]{MUMPS:2}
\bibinfo{author}{P.~R. Amestoy}, \bibinfo{author}{A.~Guermouche},
  \bibinfo{author}{J.-Y. L'Excellent}, \bibinfo{author}{S.~Pralet},
\newblock \bibinfo{title}{Hybrid scheduling for the parallel solution of linear
  systems},
\newblock \bibinfo{journal}{Parallel Comput.} \bibinfo{volume}{32}
  (\bibinfo{year}{2006}) \bibinfo{pages}{136--156}.
  \DOIprefix\doi{10.1016/j.parco.2005.07.004}.
\bibitem[{Balay et~al.(2016{\natexlab{a}})Balay, Abhyankar, Adams, Brown,
  Brune, Buschelman, Dalcin, Eijkhout, Gropp, Kaushik, Knepley, {Curfman
  McInnes}, Rupp, Smith, Zampini, Zhang, and Zhang}]{petsc-web-page}
\bibinfo{author}{S.~Balay}, \bibinfo{author}{S.~Abhyankar},
  \bibinfo{author}{M.~F. Adams}, \bibinfo{author}{J.~Brown},
  \bibinfo{author}{P.~Brune}, \bibinfo{author}{K.~Buschelman},
  \bibinfo{author}{L.~Dalcin}, \bibinfo{author}{V.~Eijkhout},
  \bibinfo{author}{W.~D. Gropp}, \bibinfo{author}{D.~Kaushik},
  \bibinfo{author}{M.~G. Knepley}, \bibinfo{author}{L.~{Curfman McInnes}},
  \bibinfo{author}{K.~Rupp}, \bibinfo{author}{B.~F. Smith},
  \bibinfo{author}{S.~Zampini}, \bibinfo{author}{H.~Zhang},
  \bibinfo{author}{H.~Zhang}, \bibinfo{title}{{PETS}c {W}eb page},
  \bibinfo{howpublished}{\url{http://www.mcs.anl.gov/petsc}},
  \bibinfo{year}{2016}{\natexlab{a}}.
\bibitem[{Balay et~al.(2016{\natexlab{b}})Balay, Abhyankar, Adams, Brown,
  Brune, Buschelman, Dalcin, Eijkhout, Gropp, Kaushik, Knepley, {Curfman
  McInnes}, Rupp, Smith, Zampini, Zhang, and Zhang}]{petsc-user-ref}
\bibinfo{author}{S.~Balay}, \bibinfo{author}{S.~Abhyankar},
  \bibinfo{author}{M.~F. Adams}, \bibinfo{author}{J.~Brown},
  \bibinfo{author}{P.~Brune}, \bibinfo{author}{K.~Buschelman},
  \bibinfo{author}{L.~Dalcin}, \bibinfo{author}{V.~Eijkhout},
  \bibinfo{author}{W.~D. Gropp}, \bibinfo{author}{D.~Kaushik},
  \bibinfo{author}{M.~G. Knepley}, \bibinfo{author}{L.~{Curfman McInnes}},
  \bibinfo{author}{K.~Rupp}, \bibinfo{author}{B.~F. Smith},
  \bibinfo{author}{S.~Zampini}, \bibinfo{author}{H.~Zhang},
  \bibinfo{author}{H.~Zhang}, \bibinfo{title}{{PETS}c Users Manual},
  \bibinfo{type}{Technical Report} \bibinfo{number}{ANL-95/11 - Revision 3.7},
  Argonne National Laboratory, \bibinfo{year}{2016}{\natexlab{b}}. \URLprefix
  \url{http://www.mcs.anl.gov/petsc}.
\bibitem[{Balay et~al.(1997)Balay, Gropp, {Curfman McInnes}, and
  Smith}]{petsc-efficient}
\bibinfo{author}{S.~Balay}, \bibinfo{author}{W.~D. Gropp},
  \bibinfo{author}{L.~{Curfman McInnes}}, \bibinfo{author}{B.~F. Smith},
\newblock \bibinfo{title}{Efficient management of parallelism in object
  oriented numerical software libraries},
\newblock in: \bibinfo{editor}{E.~Arge}, \bibinfo{editor}{A.~M. Bruaset},
  \bibinfo{editor}{H.~P. Langtangen} (Eds.), \bibinfo{booktitle}{Modern
  Software Tools in Scientific Computing}, \bibinfo{organization}{$\ $},
  \bibinfo{publisher}{Birkh{\"{a}}user Press}, \bibinfo{year}{1997}, pp.
  \bibinfo{pages}{163--202}.
\bibitem[{Rivi\`ere(2008)}]{Riviere:book}
\bibinfo{author}{B.~Rivi\`ere}, \bibinfo{title}{Discontinuous {G}alerkin
  Methods for Solving Elliptic and Parabolic Equations},
  volume~\bibinfo{volume}{35} of \textit{\bibinfo{series}{Frontiers in Applied
  Mathematics}}, \bibinfo{publisher}{Society for Industrial and Applied
  Mathematics}, \bibinfo{address}{Philadelphia}, \bibinfo{year}{2008}.

\end{thebibliography}
\end{document}